\documentclass[AMA,STIX1COL]{wileyNJD-v2} 

\articletype{Article Type}%

\received{26 April 2016}
\revised{6 June 2016}
\accepted{6 June 2016}

\usepackage{xcolor}
\usepackage{automated_macros}
\usepackage{anyfontsize}
\usepackage{amsmath}
\usepackage{bm}
\usepackage{subcaption}
\usepackage[usestackEOL]{stackengine}
\DeclareMathAlphabet{\pazocal}{OMS}{zplm}{m}{n}

\newcommand{\afterequation}{\vspace{2.5mm}}
\newcommand{\aftertable}{\vspace{-1cm}}

\newlength{\tempdima}
\newcommand{\rownameform}[1]
{\rotatebox{90}{\makebox[\tempdima][c]{\hspace{4cm}{{#1}}}}}

\newif\ifhidden
\hiddentrue 

\graphicspath{{figures/}}

\begin{document}

\title{Solving Forward and Inverse Problems of Contact Mechanics using Physics-Informed Neural Networks}

\author[1]{Tarik Sahin*}

\author[2]{Max von Danwitz}

\author[1,2]{Alexander Popp}

\authormark{SAHIN \textsc{et al}}

\address[1]{\orgdiv{Institute for Mathematics and Computer-Based Simulation (IMCS)}, 
\orgname{ University of the Bundeswehr Munich}, \orgaddress{\country{Germany}}}

\address[2]{\orgdiv{Institute for the Protection of Terrestrial Infrastructures}, 
\orgname{German Aerospace Center (DLR)}, \orgaddress{\country{Germany}}}


\corres{*Tarik Sahin, Institute for Mathematics and Computer-Based Simulation (IMCS), 
University of the
Bundeswehr Munich, Werner-Heisenberg-Weg 39, D-85577 Neubiberg, Germany. \email{tarik.sahin@unibw.de}}

\presentaddress{Werner-Heisenberg-Weg 39, D-85577
Neubiberg, Germany.}

\abstract[Abstract]{This paper explores the ability of physics-informed neural networks (PINNs) to solve forward and inverse problems of contact mechanics for small deformation elasticity. 
We deploy PINNs in a mixed-variable formulation enhanced by output transformation  
to enforce Dirichlet and Neumann boundary conditions as hard constraints. Inequality constraints of contact problems, namely
\textit{Karush-Kuhn-Tucker} (KKT) type conditions, are enforced as soft constraints by incorporating them into the loss
function during network training. To formulate the loss function contribution of KKT constraints, 
existing approaches applied to elastoplasticity problems are investigated and 
we explore a nonlinear complementarity problem (NCP) function, namely 
\textit{Fischer-Burmeister}, which possesses advantageous characteristics in terms of optimization. 
Based on the Hertzian contact problem, we show that PINNs can serve as pure partial differential equation (PDE) solver, as data-enhanced forward
model, as inverse solver for parameter identification, and as fast-to-evaluate surrogate model. 
Furthermore, we demonstrate the importance of choosing proper hyperparameters, e.g. loss weights, and a combination of 
\textit{Adam} and \textit{L-BFGS-B} optimizers aiming for better results in terms of accuracy and training time. }


\jnlcitation{\cname{%
\author{Sahin T.},
\author{v. Danwitz Max} and
\author{Popp A.}},
(\cyear{2023}).
\ctitle{Solving Forward and Inverse Problems of Contact Mechanics using Physics Informed Neural Networks},
\cjournal{},
\cvol{}.}

\maketitle

\section{Introduction}\label{intro}
Machine learning approaches usually require a large amount of simulation or experimental
data, which might be challenging to acquire due to the complexity of simulations and the cost of experiments. 
Also, data scarcity can cause data-driven techniques to perform poorly in terms of accuracy. This
is particularly true when using real-world observations that are noisy or datasets that are incorrectly
labeled, as there is no physics-based feedback mechanism to validate the predictions. To tackle this problem,
physics-informed neural networks (PINNs) have been developed.
PINNs integrate boundary or initial
boundary value problems and measurement data into the neural network’s loss function to compensate
for the lack of sufficient data and the black-box behavior of purely data-driven techniques \cite{raissi2019physics}.
In terms of forward problems, PINNs can serve as a partial differential equation (PDE) solver 
even in cases where domains are irregular. 
This is because PINNs utilize automatic differentiation and therefore do not require any connectivity of the sampling points, 
making them a mesh-free method \cite{arend2022truly}. Moreover, PINNs can break the curse of dimensionality 
when approximating functions in higher dimensions \cite{grohs2023proof, poggio2017and}. 
Additionally, PINNs are a good candidate for addressing inverse problems due to the easy integration of measurement data 
\cite{depina2022application}.

To take advantage of these benefits, PINNs have been employed in various fields of engineering and science 
including geosciences \cite{smith2022hyposvi}, 
fluid mechanics \cite{rao2020physics} \cite{eivazi2022physics}, 
optics and electromagnetics \cite{chen2020physics} \cite{beltran2022physics} \cite{khan2022physics}, 
and industrial applications, e.g., fatigue prognosis of a wind turbine main bearing \cite{yucesan2022hybrid}.
Based on sensor data of a physical object,
PINNs can be used in hybrid digital twins of civil engineering structures \cite{von2023hybrid} 
and for critical infrastructure protection \cite{BrucherseiferDLR}.
Particularly in solid mechanics, PINNs have been developed for solving problems of 
linear elasticity, elastodynamics, elastoplasticity \cite{haghighat2023constitutive}, 
and inverse problems for parameter identification \cite{bharadwaja2022physics}. Rao et al.\cite{rao2021physics} propose PINNs in a mixed-variable formulation to solve 
elastodynamic problems inspired by hybrid finite element analysis \cite{zienkiewicz2001displacement}. 
They introduce displacement and stress components as neural network output to enforce boundary conditions as hard constraints
by deploying additional parallel networks. 
Also, it is claimed that a mixed-variable formulation enhances the accuracy and ease of training for the network. 
Samaniego et al. \cite{samaniego2020energy} utilize energy methods to develop PINNs for solving various examples in computational mechanics, i.e. 
elastodynamics, hyperelasticity and phase field modeling of fracture. 
Lu and colleagues develop 
physics-informed neural networks with hard constraints (hPINNs) to perform topology optimization \cite{lu2021physics}.
The authors enhance the loss formulation with the penalty method and the augmented Lagrangian method to 
enforce inequality constraints as hard constraints. Moreover, they deploy output transformation to enforce
equality constraints explicitly for simple domains as introduced in the study of Lagaris et. al \cite{lagaris1998artificial}. 
In another study, Haghighat and colleagues utilize PINNs in the field 
of solid mechanics to tackle inverse problems and construct surrogate models \cite{haghighat2021physics}. Their approach involves 
parallel networks based on the mixed-variable formulation for linear elasticity, and they expand their 
methodology to address nonlinear elastoplasticity problems including classical \textit{Karush-Kuhn-Tucker} (KKT) type inequality constraints.
They enforce KKT constraints as soft constraints via a sign function, which has discontinuous gradients. 
As an extension of their previous work on elastoplasticity, Haghighat et al. \cite{haghighat2023constitutive} deploy
PINNs for constitutive model characterization and discovery through calibration by macroscopic mechanical testing on materials.
As an alternative to the sign function, they adopt the Sigmoid function to enforce KKT constraints,  
since Sigmoid has well-defined gradients, but requires an additional hyperparameter. 

As far as the authors are most aware, no previous work has been conducted on PINNs to solve contact mechanics problems.
Here, we focus on the novel application of PINNs for contact mechanics for small deformation elasticity 
including benchmark examples, e.g. contact between an
elastic block and a rigid flat surface, as well as the Hertzian contact problem. To enforce displacement and traction boundary conditions, 
we deploy PINNs with output transformation in the mixed-variable
formulation inspired by the Hellinger–Reissner principle \cite{reissner1950variational} in which 
displacement and stress fields are defined as network outputs. 
Additionally, contact problems involve a well-known set of \textit{Karush-Kuhn-Tucker} type inequality and equality 
constraints sometimes also referred to as Hertz-Signorini-Moreau conditions in the contact mechanics community. 
We enforce this given set of equations as soft constraints via three different methods: sign-based method, Sigmoid-based method
and a nonlinear complementarity problem (NCP) function, namely the \textit{Fischer-Burmeister} function. 
NCP functions enable reformulating inequalities as a system of equations, 
and have proven particularly robust and efficient for the design of semi-smooth Newton 
methods in contact analysis \cite{Popp2018a} \cite{seitz2015} \cite{Deng2015} \cite{li2012}.  

To validate our PINN formulation for contact mechanics, two examples are investigated. The first example involves 
the contact between an elastic block and a rigid flat surface where all points in the possible contact area will actually be in contact. 
The second example is the famous Hertzian contact problem, where the actual contact area will be determined as part of the solution procedure. 
Furthermore, we illustrate four distinct PINN application cases for the Hertzian contact problem. 
In the first use case, we deploy the PINN as a pure forward solver to validate our approach by comparing results with a finite element simulation. 
PINNs can easily incorporate external data, such as measurements or simulations. 
In the second scenario, we therefore utilize displacement and stress fields obtained through FEM (in the sense of "virtual experiments") to enhance the accuracy of 
our PINN model. The third application is to deploy PINNs to solve inverse problems, 
particularly identifying the prescribed external load in the Hertzian contact problem based on FEM data. 
As a fourth and final example, the load (external pressure) is considered as
another network input to construct a fast-to-evaluate surrogate model, which predicts displacement and stress fields for unseen pressure inputs. 
In the very active research field of physics-informed machine learning further advanced techniques, such as variational 
PINNs (VPINNs) \cite{berrone2022} \cite{kharazmi2021} and integrated finite element neural networks (I-FENNs) \cite{pantidis2023}, 
have been proposed recently. In particular, VPINNs as a Petrov-Galerkin scheme as compared to collocation in standard PINNs might 
be of interest for (non-smooth) contact problems. However, these recent developments are beyond the scope 
of this first study and subject to future work.

The remainder of this article is structured as follows: Section \ref{sec2} summarizes the fundamental equations and 
constraints of contact mechanics with small deformation elasticity. Also, the basics of the so-called mixed-variable formulation based on 
the Hellinger–Reissner principle are given. In Section \ref{sec3}, a generalized formulation of PINNs 
with output transformation is outlined in detail, which in principle allows for the solution of arbitrary partial differential equations (PDEs). 
We narrow the field of interest down to solid and contact mechanics problems based on a mixed-variable formulation, and therefore different methods 
to enforce KKT constraints are explained. Several benchmark examples
are analyzed in Section \ref{sec4}, including the Lam\'e problem of elasticity, 
contact between an elastic block and a rigid domain, and the Hertzian contact problem.  
Section \ref{sec5} concludes the paper by summarizing our key findings and providing an outlook on future research directions.

\section{Problem Formulation}\label{sec2}
\ifhidden
{}
\fi

\subsection{Contact mechanics}\label{sec2_1}
\begin{figure}[thbp]
    \centering
    \def\dist{\hspace*{5.25cm}}
    \stackinset{c}{-0.2in}{b}{-0.25in}{(a) \dist (b) \dist (c)}{
    \includegraphics[scale=0.625]{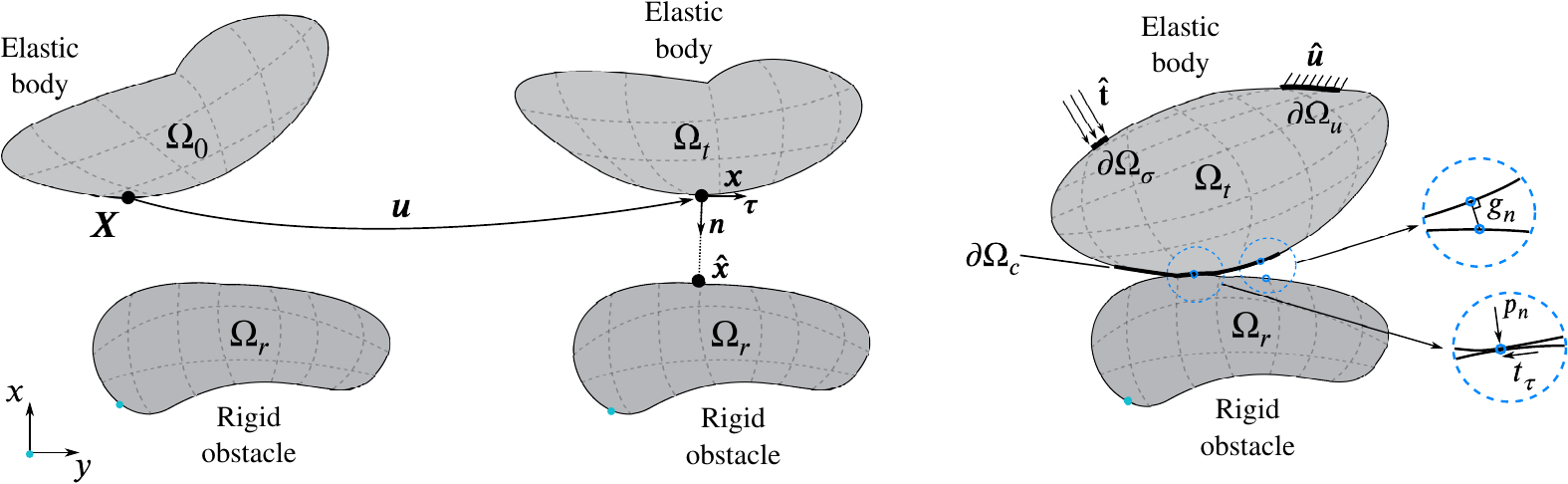} 
    } 
    \caption{Contact problem between an elastic body and a rigid obstacle. 
    (a) Reference configuration, (b) current configuration, 
    (c) accompanying boundary conditions, illustration of the gap $g_n$, tangential traction $t_{\boldsymbol{\tau}}$ and contact pressure $p_n$.}
    \label{fig:contact_illus}
\end{figure}

We consider a 2D contact problem between an elastic body and a fixed rigid obstacle as illustrated in Fig. \ref{fig:contact_illus}.
In the reference configuration, the elastic body is denoted by $\Omega_{0}$,
and in the current configuration, it is represented by $\Omega_{t}$ while the rigid obstacle has the same configuration $\Omega_{r}$.  
Fig. \ref{fig:contact_illus}c shows a configuration for which two bodies come into contact. The surface of the elastic body
can be partitioned into three sections: the Dirichlet boundary $\partial \Omega_{u}$, where displacements are prescribed, 
the Neumann boundary $\partial \Omega_{\sigma}$, where tractions are given, 
and the potential contact boundary $\partial \Omega_{c}$
where contact constraints are imposed. The actual contact surface is a subset of $\partial \Omega_{c}$ and is sought for 
during the solution procedure. 

Let us consider the boundary value problem (BVP) of small deformation
elasticity 
\begin{align}
    \label{eq:lin_equi}
    \boldsymbol{\nabla \cdot } \boldsymbol{\sigma} + \bodyforce & = \mathbf{0} \hspace{0.9cm} \text{in } \Omega, \hspace*{0.91cm} \rightarrow \text{(BE)} \\ 
    \boldsymbol{u} &=\boldsymbol{\hat u}  \hspace{0.9cm} \text{on } \partial \Omega_u, \hspace*{0.5cm} \rightarrow \text{(DBC)} \\ 
    \boldsymbol{\sigma} \cdot \mathbf{n} &=\mathbf{\hat t} \hspace{0.95cm} \text{on } \partial \Omega_{\sigma} \hspace*{0.53cm} \rightarrow \text{(NBC)} 
\end{align}

\afterequation
where $\boldsymbol{\sigma}$ denotes the Chauchy stress tensor, $\boldsymbol{u}$ is the displacement vector
representing the so-called primal variable, $\bodyforce$ denotes the body force vector, and $\boldsymbol{n}$ is the unit outward 
normal vector. Prescribed displacements are represented by 
$\boldsymbol{\hat u}$ on $\partial \Omega_u$, and $\mathbf{\hat t}$ denotes
prescribed tractions on $\partial \Omega_{\sigma}$. Abbreviations BE, DBC and NBC denote the balance equation, 
Dirichlet boundary condition and Neumann boundary condition, respectively. 
The kinematic equation (KE) and constitutive equation (CE) for the deformable body are
expressed as 
\begin{align}
    \label{eq:compatibility}
    \bm{\varepsilon} &= \frac{1}{2}(\boldsymbol{\nabla}\boldsymbol{u}+\boldsymbol{\nabla}\boldsymbol{u}^T), 
    \hspace*{0.7cm} \rightarrow \text{(KE)} \\
    \label{eq:hooke}
    \boldsymbol{\sigma} &= \mathbb{C} : \bm{\varepsilon}. 
    \hspace*{1.8cm} \rightarrow \text{(CE)}
\end{align}

\afterequation
Here, $\bm{\varepsilon}$ is the infinitesimal strain tensor and $\mathbb{C}$ is the fourth-order elasticity tensor.
In the specific case of linear isotropic elasticity, the constitutive equation can be expressed via Hooke's law as
\begin{equation}
    \boldsymbol{\sigma} = \lambda \tr(\bm{\varepsilon})\mathbf{I}+2\mu\bm{\varepsilon},
\end{equation}

\afterequation
where $\lambda$ and $\mu$ are the Lam\'e parameters, $\tr(\cdot)$ is the \textit{trace} operator to sum strain 
components on the main diagonal and $\mathbf{I}$ is the identity tensor. 

The displacement vector $\boldsymbol{u}$ can be obtained for the elastic body by describing the motion from the reference configuration $\boldsymbol{X}$ 
to the current configuration $\boldsymbol{x}$ as follows (see Figs. \ref{fig:contact_illus}a,\ref{fig:contact_illus}b) 
\begin{equation}
    \label{eq:displacementvector}
    \boldsymbol{u} = \boldsymbol{x} - \boldsymbol{X}.
\end{equation}
 
The gap function (GF) $g_n$ is defined as a distance measure between elastic and rigid bodies
in the current configuration as 
\begin{equation}
    g_n = -\boldsymbol{n}\cdot(\boldsymbol{x}-\boldsymbol{\hat x}). \hspace*{0.7cm} \rightarrow \text{(GF)}
\end{equation}

\afterequation
The term $\boldsymbol{\hat x}$ denotes the so-called closest point projection of $\boldsymbol{x}$ onto the surface of $\Omega_r$ 
(see Fig. \ref{fig:contact_illus}b).
Since all contact constraints will be defined in the current configuration, 
$p_n$ and $t_{\tau}$ can be obtained by traction vector decomposition (TVD) of the contact traction vector $\mathbf{t}_c$ as

\begin{equation}
    \label{eq:contact_traction1}
    \mathbf{t}_c = p_n \boldsymbol{n}+t_{\tau} \boldsymbol{\tau}, 
    \quad p_n=\mathbf{t}_c \cdot \boldsymbol{n}, 
    \quad t_{\tau}=\mathbf{t}_c \cdot \boldsymbol{\tau}, \hspace*{0.4cm} \rightarrow \text{(TVD)}
\end{equation}

where 
\begin{equation}
    \label{eq:contact_traction2}
    \mathbf{ t}_c = \boldsymbol{\sigma} \cdot \boldsymbol{n} \quad\quad \text{on } \partial \Omega_c. 
    \hspace*{0.4cm} \rightarrow \text{(CST)}
\end{equation}

\afterequation
Cauchy's stress theorem (CST)
states that the stress tensor $\boldsymbol{\sigma}$ maps the normal vector 
to the traction vector $\mathbf{t}_c$.
Note that the boundary vectors $\boldsymbol{n}$ and $\boldsymbol{\tau}$ can be computed on the reference configuration assuming 
small deformation elasticity (linear elasticity) \cite{Santapuri2015}.

For a frictionless contact problem, we define the
classical set of Karush-Kuhn-Tucker (KKT) conditions, commonly also referred to as Hertz-Signorini-Moreau (HSM) conditions, 
and the frictionless sliding condition (FSC) as
\begin{align}
    \label{eq:gn}
    g_n & \geqslant 0, \\[-0.1em]
    \label{eq:pn}
    p_n & \leqslant 0, \hspace*{0.4cm} \rightarrow \text{(KKT)} \\[-0.1em]
    \label{eq:gn_pn}
    p_n g_n & =0, \\[0.25em]
    \label{eq:tn}
    t_{\tau} & =0. \hspace*{0.4cm} \rightarrow \text{(FSC)}
\end{align}

\afterequation
Eq. \ref{eq:gn} enforces
the kinematic aspect of non-penetration as shown in Fig. \ref{fig:contact_illus}c. 
If two bodies are in contact, the gap vanishes, i.e., $g_n=0$. The term $p_n$ denotes 
the normal component of the contact traction, i.e., the contact pressure. Correspondingly, 
Eq. \ref{eq:pn} guarantees that no adhesive stresses are allowed in the contact zone. 
Furthermore, the complementarity requirement in Eq. \ref{eq:gn_pn} necessitates that the gap should be zero 
when there is a non-zero contact pressure (point in contact),  
and the contact pressure should be zero if there is a positive gap (point not in contact). 
It should be noted that the tangential component of the traction vector vanishes for frictionless contact, resulting in Eq. \ref{eq:tn}. 
For additional information regarding more complex contact constitutive laws including friction, we
refer to \cite{Popp2018a, popp2010, popp2014, wriggers2006computational}. 

\subsection{Mixed-variable formulation: the Hellinger–Reissner principle}\label{sec2_3}

\begin{figure}[thbp]
    \centering
    \includegraphics[scale=0.65]{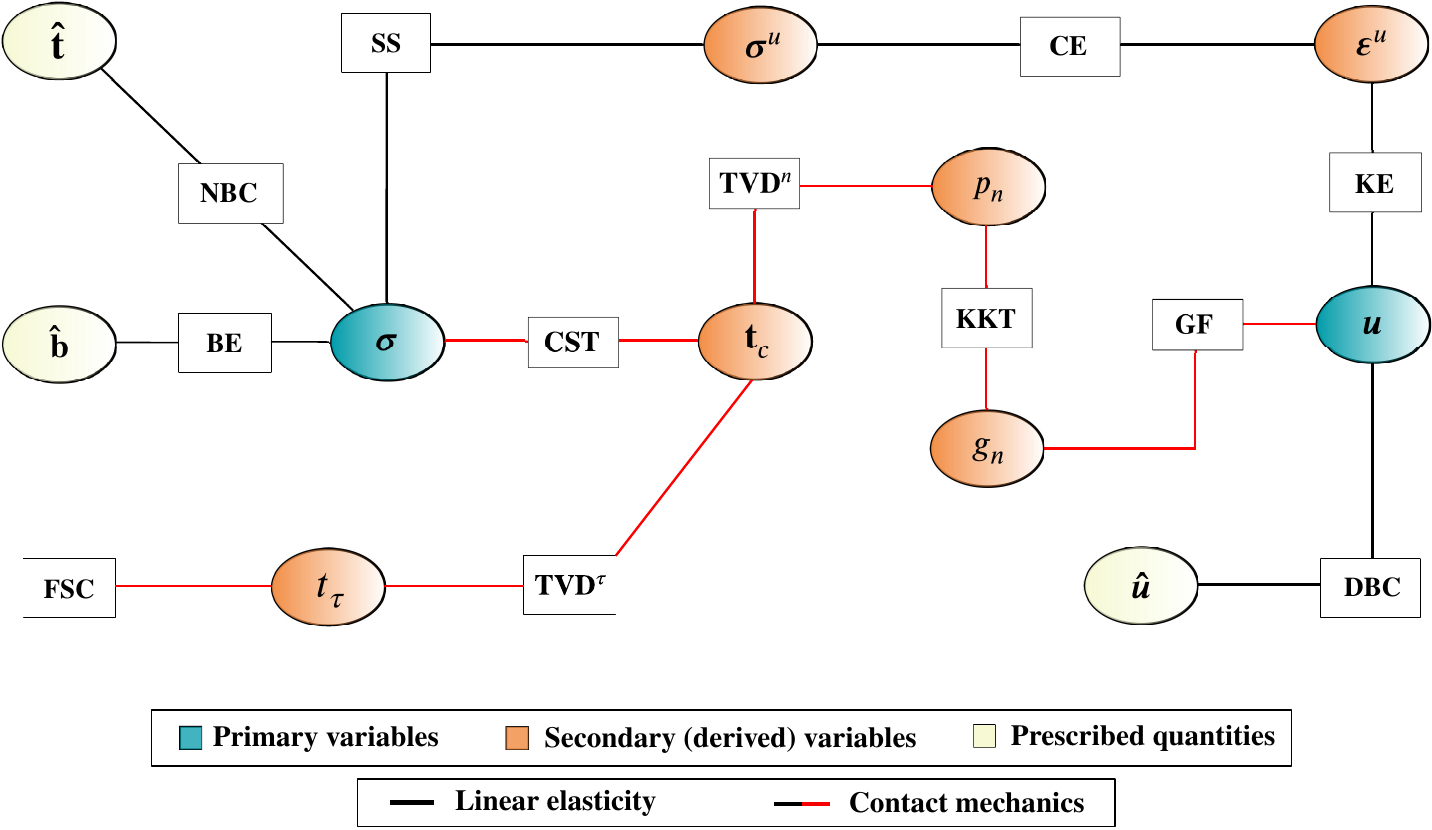}  
    \caption{Tonti's diagram of Hellinger–Reissner (HR) principle for contact problems with small strain theory and 
    frictionless sliding condition.}
    \label{fig:mixed_formulation_base}
\end{figure}

Inspired by the Hellinger–Reissner principle \cite{zienkiewicz2001displacement} \cite{reissner1950variational}, we construct a Tonti's diagram \cite{tonti1976reason}
to solve contact problems based on two primary variables:
displacements $\boldsymbol{u}$ and stresses $\boldsymbol{\sigma}$ (see Fig. 
\ref{fig:mixed_formulation_base}).
The secondary (slave) variables are intermediate and they are derived from the primary variables, i.e., 
$\boldsymbol{\sigma}^u$ and $\boldsymbol{\varepsilon}^u$.

As shown in Fig. \ref{fig:mixed_formulation_base}, the Tonti diagram summarizes the governing equations Eq. 
\ref{eq:lin_equi}-\ref{eq:tn}.
The stress-to-stress coupling (SS) between the primary variable $\boldsymbol{\sigma}$ and
the secondary variable $\boldsymbol{\sigma}^u$, defined as $\boldsymbol{\sigma} = \boldsymbol{\sigma}^u$, 
ensures that the two master fields remain compatible.

\section{Physics-Informed Neural Networks for Solid and Contact Mechanics}\label{sec3}

\subsection{Generic PINNs with output transformation}\label{sec3_1}
\ifhidden
{}
\fi

A generalized formulation for partial differential equations can be expressed in residual form with accompanying boundary conditions as 
\begin{equation}
    \begin{aligned}
    \label{eq:governing}
    \mathcal{R}\bigl[ \latentsol (\pdeInput) \bigr]&=\boldsymbol{0}, \quad \text{on} \ \Omega, \\
    \mathcal{B}\bigl[ \latentsol (\pdeInput) \bigr]-g(\pdeInput)&=\boldsymbol{0}, \quad \text{on} \ \partial \Omega.
    \end{aligned}
\end{equation}

\vspace*{2.5mm}
\noindent Here, $\mathcal{R}[\cdot]$ denotes a differential operator acting on a unknown solution $\latentsol$, 
$\mathcal{B}[\cdot]$ is the boundary operator, $g(\pdeInput)$ represents the prescribed boundary condition,
$\pdeInput$ are the spatial coordinates that span the domain $\Omega$ and
the boundary $\partial \Omega$. 

Consider a fully-connected $L$-layers neural network to construct an approximated solution $\pinnsol$ to the BVP as follows \cite{lawrence1993} \cite{kollmannsberger2021}
\begin{equation}
    \label{eq:approx}
    \latentsol \approx \pinnsol := (\netOutput)^{'}, \quad  \netOutput : \mathbb{R}^{d} \rightarrow \mathbb{R}^{n}.
\end{equation}

\vspace*{2.5mm}
\noindent $\netOutput$ represents the network output in the output layer $L$, trainable network parameters, namely weights and biases, are denoted as 
$\boldsymbol{\theta}$, $(.)'$ represents a user-defined output transformation \cite{lu2021physics}, and
$\netInput$ is the network input such that $\pdeInput \subset \netInput$. 
Note that $\netInput$ can consist of spatial coordinates, time, and other additional input parameters.   
The network output is calculated using recursive $L-1$ element-wise operations between the input layer and hidden layers as 

\begin{equation}
    \begin{aligned}
        \text{input layer}\ \  & \ \rightarrow \ \netOutputLayer{1} = \netInput, \\
        \text{hidden layers} & \ \rightarrow \ \netOutputLayer{l} =  \psi(\boldsymbol{\theta}^{l} \cdot \netOutputLayer{l-1}),\ \text{for }2 \leqslant l \leqslant L-1, \\
        \text{output layer} \ & \ \rightarrow \ \netOutputLayer{L} =  \boldsymbol{\theta}^{L} \cdot \netOutputLayer{L-1},
    \end{aligned}   
\end{equation}

\vspace*{2.5mm}
\noindent where $\psi$ denotes the activation function that adds non-linearity to the layer output. 

Output transformation enables the neural network to enforce boundary conditions explicitly.
A user-defined output transformation can be obtained with suitable helper functions as 
\begin{equation}
    \pinnsolComponent{i} = g(\netInput)+ s \cdot h(\netInput) \cdot \netOutputUComponent{i},
\end{equation}

\vspace*{2.5mm}
\noindent where $g$ is the prescribed boundary condition, $s$ is a scaling parameter and 
$h$ is a distance-to-boundary function fulfilling the following conditions:
\begin{equation}
    \label{eq:h_func}
    \centering
    \begin{aligned}
        h(\netInput)&=0, \quad \text{on} \ \partial \Omega, \\ 
        h(\netInput)&>0, \quad \text{in} \ \Omega \setminus \partial \Omega.
    \end{aligned}
\end{equation}

\vspace*{2.5mm}
\noindent 
For simple boundaries $\partial \Omega$, it is relatively easy to define an appropriate distance-to-boundary function. 
On the other hand, it can become a quite challenging task in the case of arbitrary boundaries. 
One method to find a generalized distance-to-boundary function is to use NURBS parametrizations \cite{saidaoui2022deep}.
Moreover, the scaling parameter can help the optimizer avoid getting stuck in
local minima by balancing target governing equations (see Sec. \ref{sec3_4}).

\begin{figure}[thbp]
    \centering
    \includegraphics[scale=0.55]{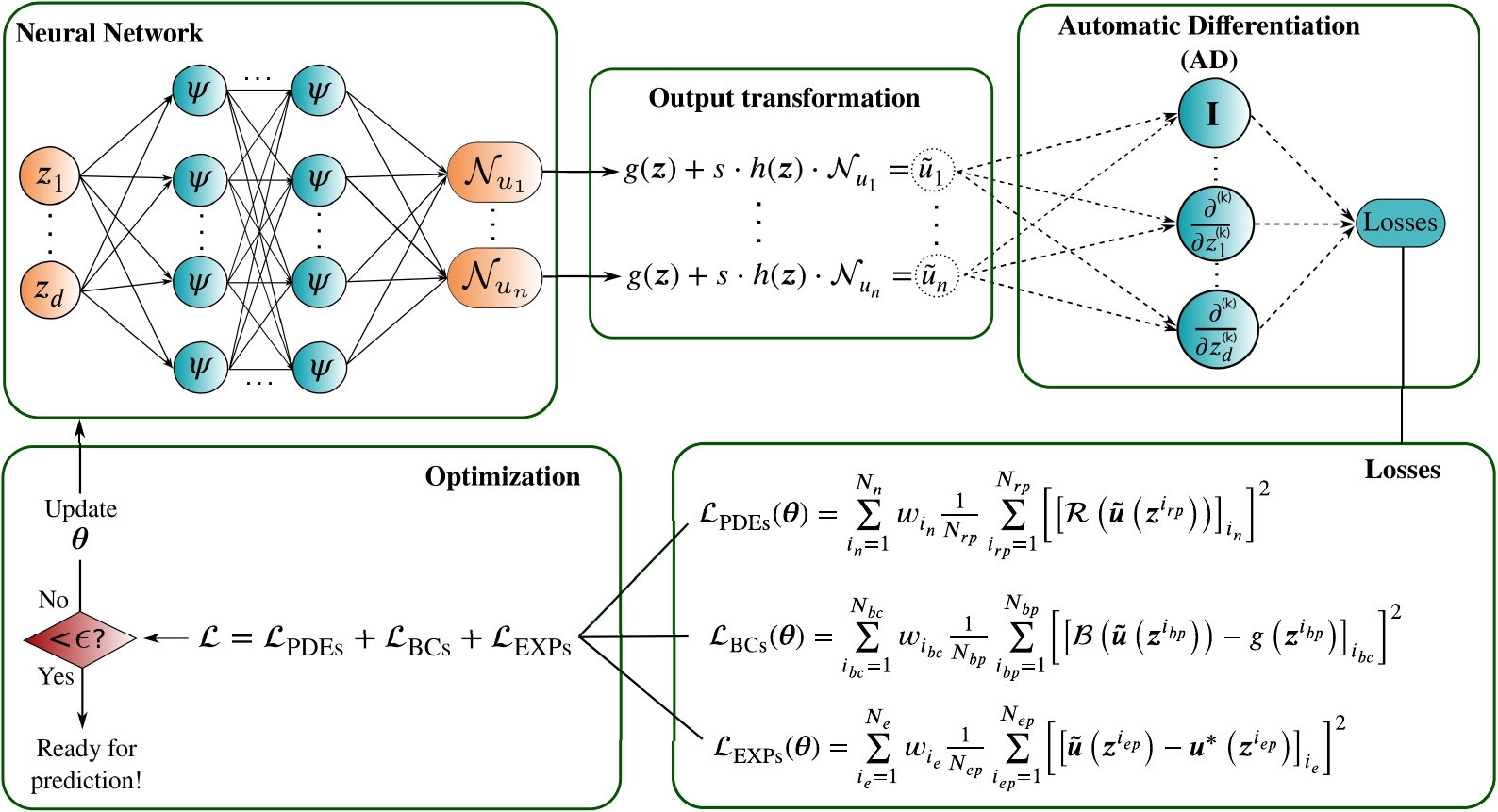}  
    \caption{The general representation of a physics-informed neural network for a BVP. }
    \label{fig:pinn_arch}
\end{figure}
To ensure that $\pinnsol$ is a reasonable approximation of $\latentsol$, the network parameters must be determined accordingly to the BVP. 
As shown in Fig. \ref{fig:pinn_arch}, the overall loss $\pazocal{L}(\boldsymbol{\theta})$ consists of PDE losses $\pazocal{L}_{\mathrm{PDEs}}$, boundary condition losses 
$\pazocal{L}_{\mathrm{BCs}}$ and
experimental data losses $\pazocal{L}_{\mathrm{EXPs}}$. 
Note that PDE derivatives are calculated using automatic differentiation (AD) \cite{baydinAD}. 
To minimize the overall loss, 
the optimization process continues until a prescribed tolerance $\epsilon$ is reached so that optimal network parameters $\boldsymbol{\theta}^*$ are
calculated
\begin{equation}
    \boldsymbol{\theta}^*=\underset{\boldsymbol{\theta}}{\arg \min } \ \mathcal{L}(\boldsymbol{\theta}).
\end{equation}

\vspace*{2.5mm}
For each loss term, the inner loop sums up the \textit{mean squared error} contributions of data points collected inside the domain or on the boundary.
Specifically, $\{\netInput^{i_{rp}}\}_{i_{rp}=1}^{N_{rp}}$ denotes the collocation points in the domain,
$\{\netInput^{i_{bp}}\}_{i_{bp}=1}^{N_{bp}}$ are the boundary points corresponding to the prescribed boundary conditions,
and $\{\netInput^{i_{ep}}\}_{i_{ep}=1}^{N_{ep}}$ represents the points on which measurement data $\boldsymbol{u^*}$ is available.
As the PDE residual $\mathcal{R}$ might have multiple terms and usually more than one boundary condition is defined, 
we use a lower index to explicitly point out that the 
inner summation is done for a given specific component, i.e. $\mathcal{R}_{i_n=1}$.
Consequently, the outer loop sums up the weighted contributions coming from individual components of the loss function. The terms $\{w_{i_n}\}_{i_n=1}^{N_n}$, $\{w_{i_{bc}}\}_{i_{bc}=1}^{N_{bc}}$
and $\{w_{i_e}\}_{i_e=1}^{N_e}$ denote the loss weights for the individual components of $\pazocal{L}_{\mathrm{PDEs}}$,
$\pazocal{L}_{\mathrm{BCs}}$, $\pazocal{L}_{\mathrm{EXPs}}$, respectively. We observe that the weighting of loss terms 
can be quite crucial for the convergence of the overall loss, since it avoids the optimizer expanding greater efforts on loss components that 
have a larger order of magnitude compared to others.
It should be noted that the identical collocation points $\{\netInput^{i_{rp}}\}_{i_{rp}=1}^{N_{rp}}$ are employed in the calculation of every component of $\mathcal{L}_{\mathrm{PDEs}}$ . 
However, the boundary and experimental points may vary in each $\mathcal{L}_{\mathrm{BCs}}$ and $\mathcal{L}_{\mathrm{EXPs}}$ contribution.

\subsection{Application of PINNs to solid and contact mechanics}\label{sec3_2}
\begin{figure}[thbp]
    \centering
    \includegraphics[scale=0.75]{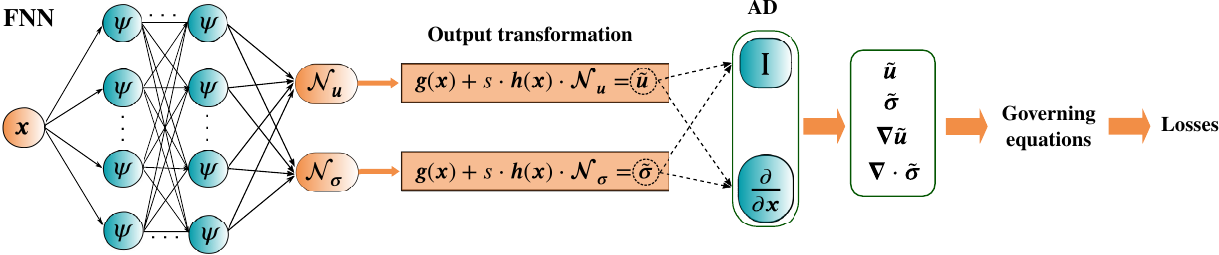}  
    \caption{Physics-informed neural networks in the mixed-variable form to solve quasi-static solid and contact mechanics problems
    without additional network parameters.}
    \label{fig:mixed_formulation}
\end{figure}

\ifhidden
{}
\fi
In the context of solid and contact mechanics problems, 
we use PINNs with output transformation in a mixed-variable formulation. 
In the mixed-variable formulation for quasi-static problems without additional network input parameters (i.e. $\netInput = \boldsymbol{x}$), 
a fully-connected neural network (FNN) maps the given spatial
coordinates $\boldsymbol{x}$ to the displacement vector $\boldsymbol{u}$ and stress tensor $\boldsymbol{\sigma}$.
In other words, the displacement and stress fields are chosen as the quantities of interest that the FNN approximates as
(see Fig. \ref{fig:mixed_formulation})
\begin{equation}
        \label{eq:approx_mixed}
        \pinnUVector := (\mathcal{N}_{\boldsymbol{u}}(\boldsymbol{x};\boldsymbol{\theta}))^{'} 
        \hspace*{0.3cm} \text{and} \hspace*{0.3cm}
        \pinnSVector := (\mathcal{N}_{\boldsymbol{\sigma}}(\boldsymbol{x};\boldsymbol{\theta}))^{'}. \\
\end{equation}

\afterequation
Combining the information provided in Fig. \ref{fig:pinn_arch} for losses of general PINNs with the governing equations of solid mechanics,
we obtain the total loss $\pazocal{L}_E $ for linear elasticity (without contact) in the mixed-variable formulation with additional experimental data as 
\begin{equation}
    \centering
    \label{eq:pinn_elasticity}
        \pazocal{L}_E = \pazocal{L}_{\mathrm{PDEs}} + \pazocal{L}_{\mathrm{DBCs}} + \pazocal{L}_{\mathrm{NBCs}} + \pazocal{L}_{\mathrm{EXPs}},
\end{equation}

\afterequation
where 
\begin{equation}
    \centering
    \label{eq:pinn_elasticity_detail}
    \begin{aligned}
        \mathcal{L}_{\mathrm{PDEs}}&= \sum_{i_{n}=1}^{N_{m}} {w}_{i_{n}} \frac{1}{N_{rp}} \sum_{i_{rp}=1}^{N_{rp}}\Bigl[ \bigl[\boldsymbol{\nabla}\cdot \pinnSVector \left(\boldsymbol{x}^{i_{rp}} \right) + \bodyforce\left(\boldsymbol{x}^{i_{rp}} \right) \bigr]_{i_{n}}\Bigr]^2+
        \sum_{i_{n}=N_{m}+1}^{N_{n}} {w}_{i_{n}} \frac{1}{N_{rp}} \sum_{i_{rp}=1}^{N_{rp}}\Bigl[ \bigl[\pinnSVector\left(\boldsymbol{x}^{i_{rp}} \right) -\mathbb{C} : \tilde{\bm{\varepsilon}}\left(\boldsymbol{x}^{i_{rp}} \right) \bigr]_{i_{n}}\Bigr]^2, \\
        \mathcal{L}_{\mathrm{DBCs}}&=\sum_{i_{bc,D}=1}^{N_{bc,D}} {w}_{i_{bc,D}} \frac{1}{N_{bp,D}} \sum_{i_{bp,D}=1}^{N_{bp,D}}\Bigl[ \bigl[ \pinnUVector\left(\boldsymbol{x}^{i_{bp,D}}\right) - \boldsymbol{\hat u}\left(\boldsymbol{x}^{i_{bp,D}}\right) \bigr]_{i_{bc,D}} \Bigr]^2, \\
        \mathcal{L}_{\mathrm{NBCs}}&=\sum_{i_{bc,N}=1}^{N_{bc,N}} {w}_{i_{bc,N}} \frac{1}{N_{bp,N}} \sum_{i_{bp,N}=1}^{N_{bp,N}}\Bigl[ \bigl[ \pinnSVector\left(\boldsymbol{x}^{i_{bp,N}}\right)\cdot \boldsymbol{n} - \mathbf{\hat t}\left(\boldsymbol{x}^{i_{bp,N}}\right) \bigr]_{i_{bc,N}} \Bigr]^2 \\
        \mathcal{L}_{\mathrm{EXPs}}&=\sum_{i_{e,\boldsymbol{u}}=1}^{N_{e,\boldsymbol{u}}} {w}_{i_{e,\boldsymbol{u}}} \frac{1}{N_{ep,\boldsymbol{u}}} \sum_{i_{ep,\boldsymbol{u}}=1}^{N_{ep,\boldsymbol{u}}}\Bigl[ \bigl[ \pinnUVector\left(\boldsymbol{x}^{i_{ep,\boldsymbol{u}}}\right) - \boldsymbol{u^*}\left(\boldsymbol{x}^{i_{ep,\boldsymbol{u}}}\right) \bigr]_{i_{e,\boldsymbol{u}}} \Bigr]^2
        + \sum_{i_{e,\boldsymbol{\sigma}}=1}^{N_{e,\boldsymbol{\sigma}}} {w}_{i_{e,\boldsymbol{\sigma}}} \frac{1}{N_{ep,\boldsymbol{\sigma}}} \sum_{i_{ep,\boldsymbol{\sigma}}=1}^{N_{ep,\boldsymbol{\sigma}}}\Bigl[ \bigl[ \pinnSVector \left(\boldsymbol{x}^{i_{ep,\boldsymbol{\sigma}}}\right) - \boldsymbol{\sigma^*}\left(\boldsymbol{x}^{i_{ep,\boldsymbol{\sigma}}}\right) \bigr]_{i_{e,\boldsymbol{\sigma}}} \Bigr]^2.
    \end{aligned}
\end{equation}

\afterequation
Here, terms $\mathcal{L}_{\mathrm{DBCs}}$ and $\mathcal{L}_{\mathrm{NBCs}}$ denote losses for Dirichlet and Neumann BCs, respectively, and
the $\mathcal{L}_{\mathrm{EXPs}}$ term represents losses due to additional experimental data.
In the mixed-variable formulation, the $\mathcal{L}_{\mathrm{PDEs}}$ term is constructed in a composite form to 
fulfill both the balance equation (BE) and the stress-to-stress coupling (SS) as depicted in Fig. \ref{fig:mixed_formulation_base}. 
The index $N_{m}$ is used to distinguish the loss weights related to BE and SS. 
Since stress components are directly defined as network outputs, 
traction or Neumann BCs can be imposed as hard constraints using output transformation. 
Moreover, it is sufficient to calculate first-order derivatives of the neural network outputs with respect to the inputs,
since the governing equations in the mixed-variable formulation contain only first-order derivatives. 
An alternative to the mixed-variable formulation is the classical displacement-based formulation in which only $\boldsymbol{u}$
is considered as the network output. However, such an approach requires second-order derivatives for evaluating the balance equation, and 
traction BCs can not be enforced as hard constraints \cite{rao2021physics} \cite{sun2020surrogate}.  

Next, we construct the composite (total)
loss function $\pazocal{L}_{\mathrm{C}}$  for linear elasticity including contact as follows:

\begin{equation}
    \label{eq:loss_contact}
    \begin{aligned}
        \pazocal{L}_{\mathrm{C}} = &\ \pazocal{L}_{\mathrm{E}} + \pazocal{L}_{\mathrm{FS}} + \pazocal{L}_{\mathrm{KKT}}. \\
   \end{aligned}
\end{equation}

\afterequation
The first additional loss term
\begin{equation}
        \pazocal{L}_{\mathrm{FS}} = w^{(\mathrm{fs})} | t_{\tau}|_{\partial \Omega_c} 
\end{equation}

\afterequation
enforces the frictionless sliding condition in the contact zone ${\partial \Omega_c}$ (see Eqs. \ref{eq:contact_traction1} and 
\ref{eq:contact_traction2}). 
For simplicity, we denote the \textit{mean squared error} (MSE) as $|\cdot|$, which can be calculated as 
$\frac{1}{n}\sum_{i=1}^{n}(\cdot)^2$.
The second additional term $\pazocal{L}_{\mathrm{KKT}}$ will be elaborated in the next section. 


\afterequation
While the various ways of evaluating the normal gap $g_n$ are a matter of intense discussions, especially 
within discretization schemes such as the finite element method \cite{wriggers2006computational} \cite{yastrebov2013numerical}, a very simple gap calculation is sufficient here 
due to the fact that we only consider contact problems between an elastic body and a rigid flat surface. 
The normal gap is consistently expressed by evaluating the orthogonal projection of the elastic body onto the rigid flat surface.



\subsection{Enforcing the Karush-Kuhn-Tucker inequality constraints}\label{sec3_3}

There are several methods available to enforce inequality conditions in general. The direct approach is to 
formulate loss functions of inequalities and impose them as soft constraints 
with fixed loss weights \cite{haghighat2021physics} \cite{haghighat2023constitutive}. However, setting large loss weights can cause an ill-conditioned problem \cite{lu2021physics}. 
On the other hand, when small loss weights are chosen, the estimated solution may violate the inequalities. 
To tackle this problem, authors in \cite{lu2021physics} suggest penalty and augmented Lagrangian methods, 
well-known from constrained optimization, which construct 
loss formulations with adaptive loss weights. In the following, we investigate three methods to enforce KKT conditions 
of normal contact problems based on soft constraints.
 
\subsubsection{Sign-based method}\label{ssec_adopted_sign}
\ifhidden
{}
\fi

One possible way to enforce KKT conditions is to use the sign function\cite{haghighat2021physics},  
which leads to 

\begin{equation}
    \label{eq:adopted_sign_loss}
    \begin{aligned}
        \pazocal{L}_{\mathrm{KKT}} = &\ \pazocal{L}_{\tilde{g}_n  \geqslant 0} + \pazocal{L}_{\tilde{p}_n  \leqslant 0} + \pazocal{L}_{\tilde{p}_n \tilde{g}_n  =0} \\
        = &\ w^{(\mathrm{KKT})}_1 \bigl|\frac{1}{2}\bigl(1-\sgn(\tilde{g}_n)\bigr) \ \tilde{g}_n\bigr|_{\partial \Omega_c} + 
        w^{(\mathrm{KKT})}_2 \bigl|\frac{1}{2}\bigl(1+\sgn(\tilde{p}_n)\bigr) \ \tilde{p}_n\bigr|_{\partial \Omega_c} + 
        w^{(\mathrm{KKT})}_3 \bigl|\tilde{p}_n \ \tilde{g}_n\bigr|_{\partial \Omega_c}.
    \end{aligned}
\end{equation}

\begin{figure}[thbp]
    \centering
    \includegraphics[trim={0 0.4cm 0 0.5cm},clip,width=0.6\linewidth]{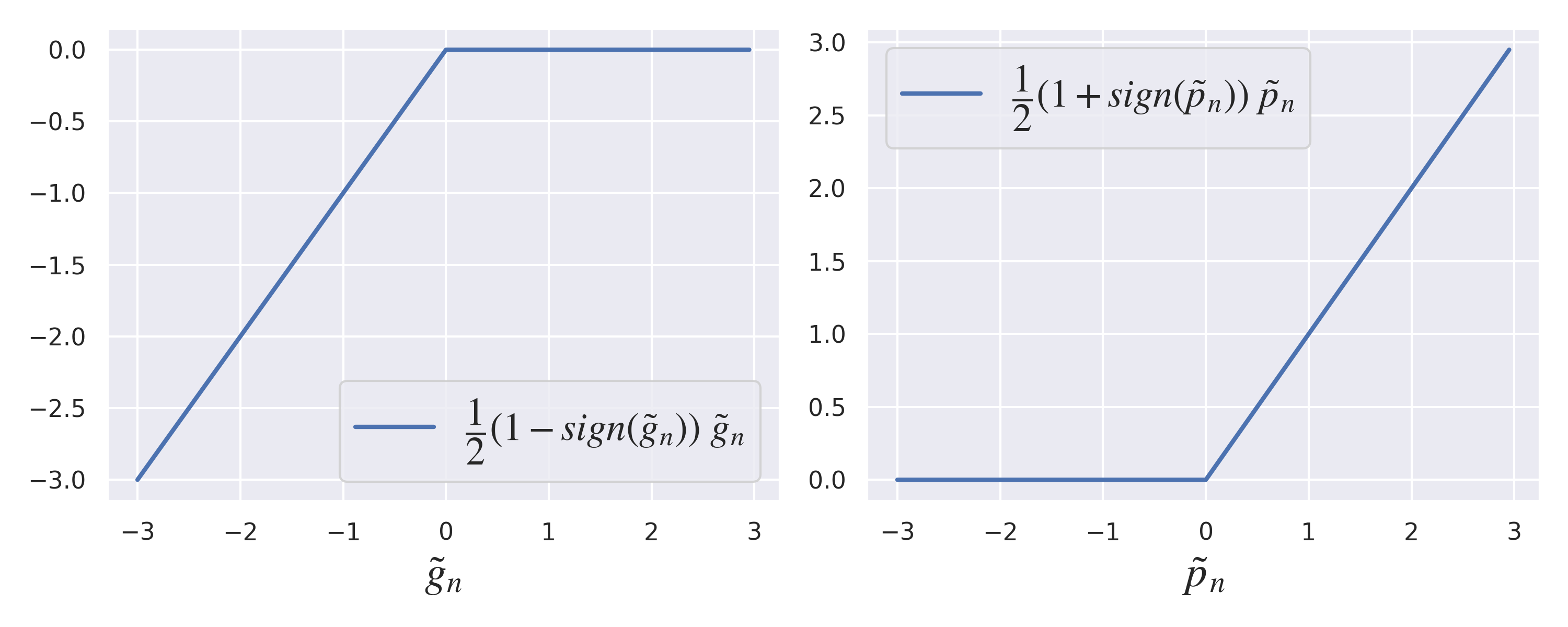}  
    \caption{An illustration of the sign-based function depending on gap $g_n$ and contact pressure $p_n$}
    \label{fig:adopted_sign}
\end{figure}

\afterequation
Here, $\{w^{(\mathrm{KKT})}_i\}_{i=1}^{3}$ represent the loss weights on the corresponding KKT condition. 
Figure \ref{fig:adopted_sign} illustrates that $\frac{1}{2}\bigl(1-\sgn(\tilde{g}_n)\bigr)\tilde{g}_n$ contributes to the loss component
$\pazocal{L}_{g_n \geqslant 0}$ when the gap $g_n$ is less than zero.
On the other hand, $\frac{1}{2}\bigl(1+\sgn(\tilde{p}_n)\bigr)\tilde{p}_n$ contributes to 
$\pazocal{L}_{p_n \leqslant 0}$ if the contact pressure $p_n$ is positive.
Note that the $\sgn$ function has gradient jumps, which is typically not a desired feature
in the context of optimization.

\subsubsection{Sigmoid-based method}\label{ssec_sigmoid}
\ifhidden
{}
\fi

An alternative approach to circumvent discontinuous gradients is to use
the Sigmoid function \cite{haghighat2023constitutive} to obtain $\pazocal{L}_{\mathrm{KKT}}$ as

\begin{equation}
    \label{eq:adopted_loss}
    \begin{aligned}
        \pazocal{L}_{\mathrm{KKT}} = &\ \pazocal{L}_{\tilde{g}_n   \geqslant 0} + \pazocal{L}_{\tilde{p}_n   \leqslant 0} + \pazocal{L}_{\tilde{p}_n  \tilde{g}_n   =0} \\
        = &\ w^{(\mathrm{KKT})}_1 \bigl|\dfrac{1}{1+e^{\delta \tilde{g}_n}} \tilde{g}_n\bigr|_{\partial \Omega_c} + 
        w^{(\mathrm{KKT})}_2 \bigl|\dfrac{1}{1+e^{-\delta \tilde{p}_n}} \tilde{p}_n\bigr|_{\partial \Omega_c} + 
        w^{(\mathrm{KKT})}_3 \bigl|\tilde{p}_n \ \tilde{g}_n\bigr|_{\partial \Omega_c},
    \end{aligned}
\end{equation}

\begin{figure}[thbp]
    \centering
    \includegraphics[scale=0.4]{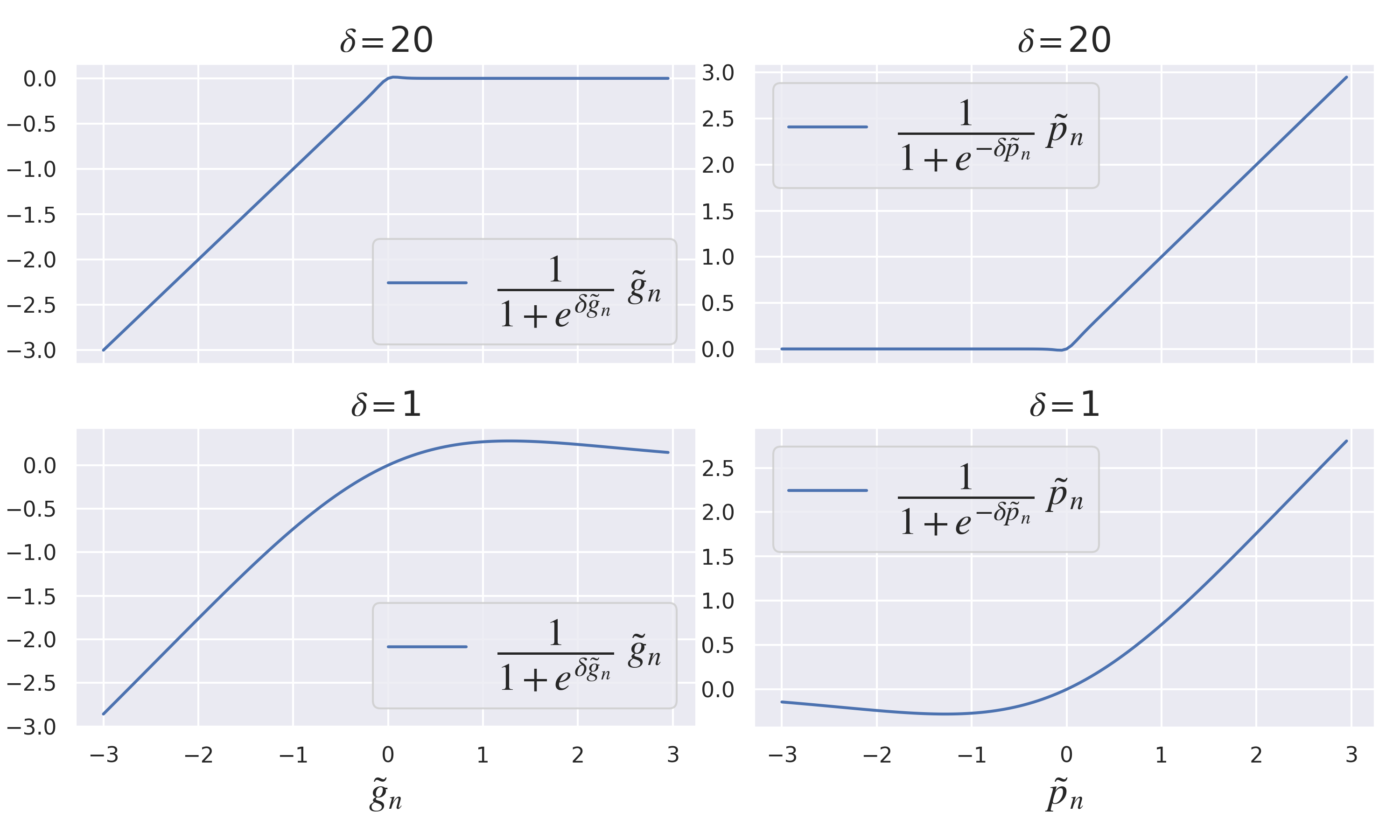}  
    \caption{An illustration of the Sigmoid-based function depending on gap $g_n$ and contact 
    pressure $p_n$ for different $\delta$ values.}
    \label{fig:adopted_sigmoid}
\end{figure}

\afterequation
\afterequation
where $\delta$ is the \textit{steepness} parameter that controls the transition
between zero and non-zero loss contributions. As depicted in Fig. \ref{fig:adopted_sigmoid}, 
the Sigmoid function avoids gradient jumps through exponential regularization term. 
However, when $\delta$ is chosen too small, e.g. $\delta=1$, then significant unphysical loss function values
are obtained for $\tilde{g}_n > 0$ and $\tilde{p}_n<0$. On the other hand, setting $\delta$ too large recovers the sign-based implementation.
Therefore, a parameter study must be conducted to find the optimal parameter $\delta_{opt}$.

\subsubsection{A nonlinear complementarity problem function: \textit{Fischer-Burmeister}}\label{ssec_fb}
\ifhidden
{}
\fi
Nonlinear complementarity problem (NCP) functions are developed based on reformulating inequalities
as equalities \cite{li2020optimization}.  One popular choice of NCP function is the 
Fischer-Burmeister function \cite{fischer1992} expressed as 

\begin{equation}
    \label{eq:fischer_burmeister}
    \phi_{\mathrm{FB}}(a,b):=a+b-\sqrt{a^2+b^2}=0 \quad \Longleftrightarrow \quad a \geqslant 0, b \geqslant 0, ab=0.
\end{equation}

\afterequation
By setting $a=\tilde{g}_n $ and $b=-\tilde{p}_n $ in the Fischer-Burmeister function, we obtain $\pazocal{L}_{\mathrm{KKT}}$ as follows

\begin{equation}
    \label{eq:fischer_burmeister_loss}
    \begin{aligned}
        \pazocal{L}_{\mathrm{KKT}} = &\ w^{(\mathrm{KKT})} \bigl|\phi_{\mathrm{FB}}(\tilde{g}_n ,-\tilde{p}_n )\bigr|_{\partial \Omega_c}
                                   = w^{(\mathrm{KKT})} \Bigl|\tilde{g}_n  - \tilde{p}_n  - \sqrt{{\tilde{g}_n}^2+\tilde{p}_n ^2}\Bigr|_{\partial \Omega_c}.
    \end{aligned}
\end{equation}

\begin{figure}[thbp]
    \centering
        \begin{subfigure}[c]{0.45\linewidth}
            \centering
            \stackinset{c}{0.2in}{t}{0.1in}{(a) \hspace*{4.7cm} (b)}{
            \includegraphics[width=0.6\linewidth]{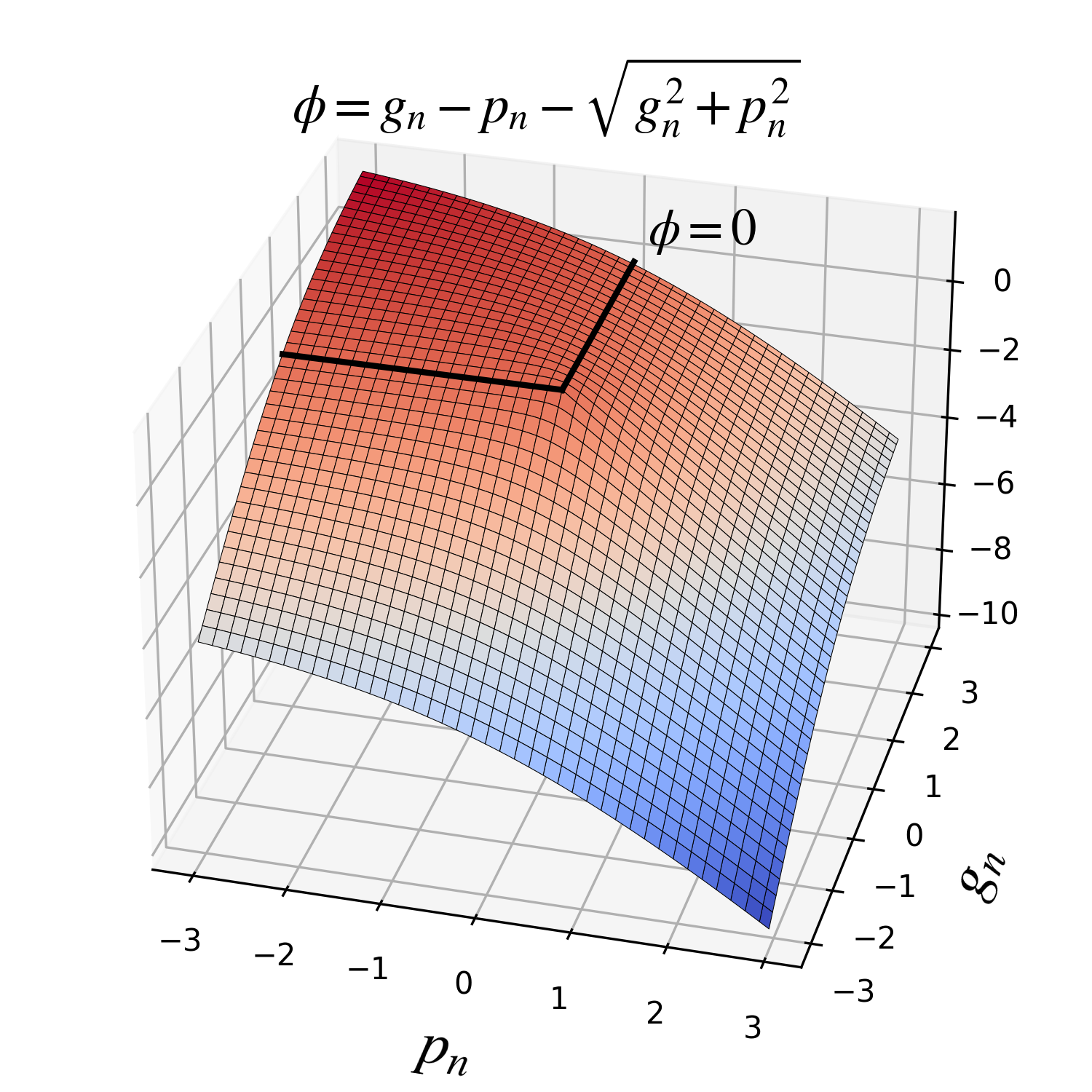}}
            \label{fig:dragratio} 
        \end{subfigure}
        \hspace{-3.0cm}
        \begin{subfigure}[c]{0.45\linewidth}
            \centering
            \vspace{0.5cm}
            \includegraphics[width=0.6\linewidth]{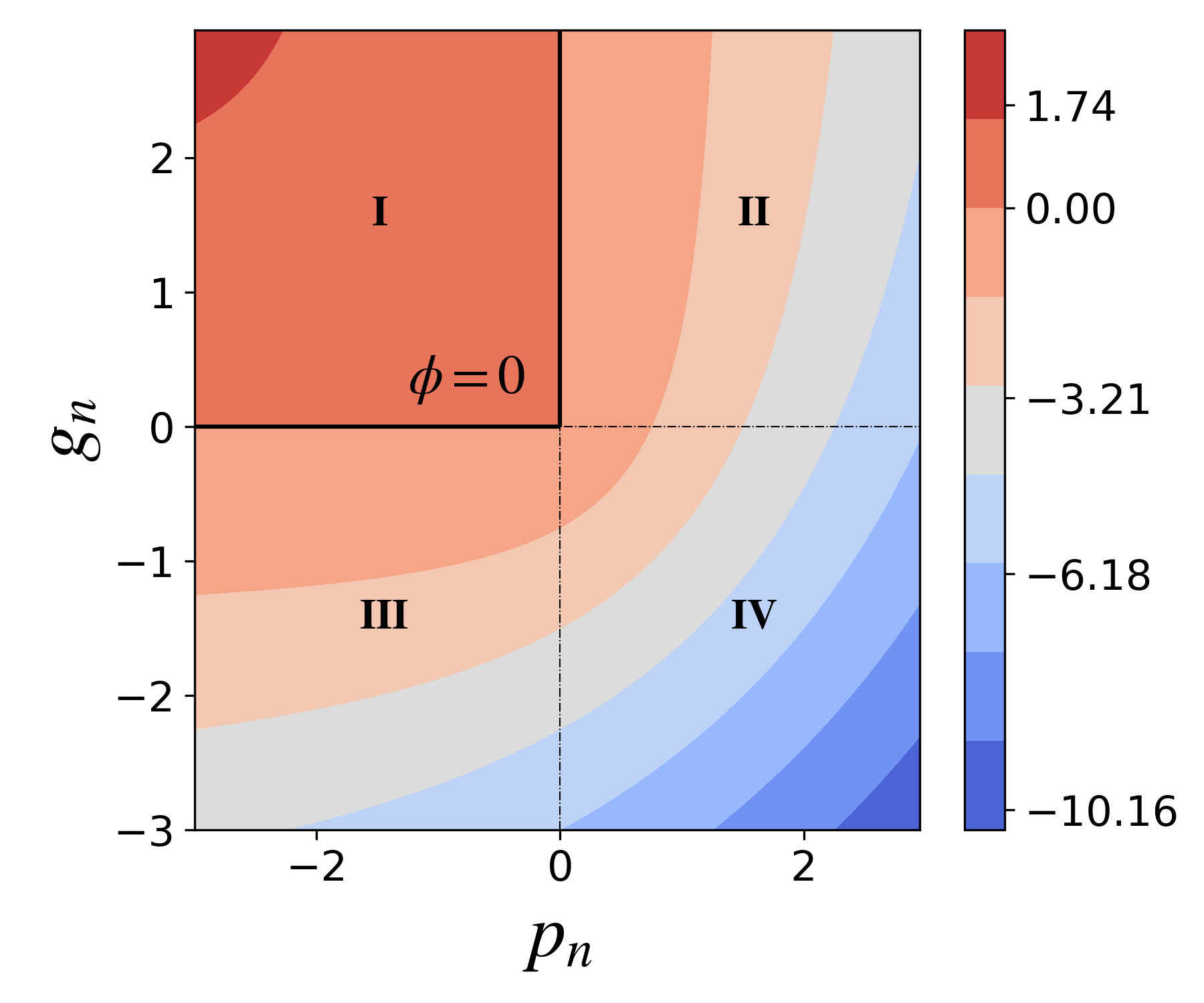}
            \label{fig:dragratio2}
        \end{subfigure}
    \centering
    \caption{The Fischer-Burmeister NCP function depending on gap $g_n$ and contact pressure $p_n$ as a 3D plot (a) and as a 2D contour plot (b).}
    \label{fig:fischer_burmeister}
\end{figure}

\afterequation
The Fischer-Burmeister function is a particularly suitable choice for typical loss calculations
based on the mean squared error (MSE), 
since $(\phi_{\mathrm{FB}})^2$ is continuously differentiable also at $a=b=0$ 
as reported \cite{sun1999ncp}. As shown in Fig. \ref{fig:fischer_burmeister}, the largest loss 
contribution comes from section IV in which both excessive penetrations, i.e. $g_n \leqslant 0$, and large adhesive stresses,
i.e. $p_n \geqslant 0$, are present. We refer to \cite{fischer1992} \cite{sun1999ncp} \cite{li2020optimization} \cite{bartel2019investigations}
for more details about the Fischer-Burmeister function.
Note also that having fewer loss terms generally
eases the optimization process as well as parameter tuning, and in contrast to the previous variants in 
Section \ref{ssec_adopted_sign} and \ref{ssec_sigmoid} only one single loss weight is required in the case of the Fischer-Burmeister NCP function.  

\subsection{Algorithmic challenges}


\subsubsection{Domination of $p_n$ over $g_n$}\label{sec3_4}

\begin{figure}[thbp]
    \centering
    \includegraphics[width=0.8\linewidth]{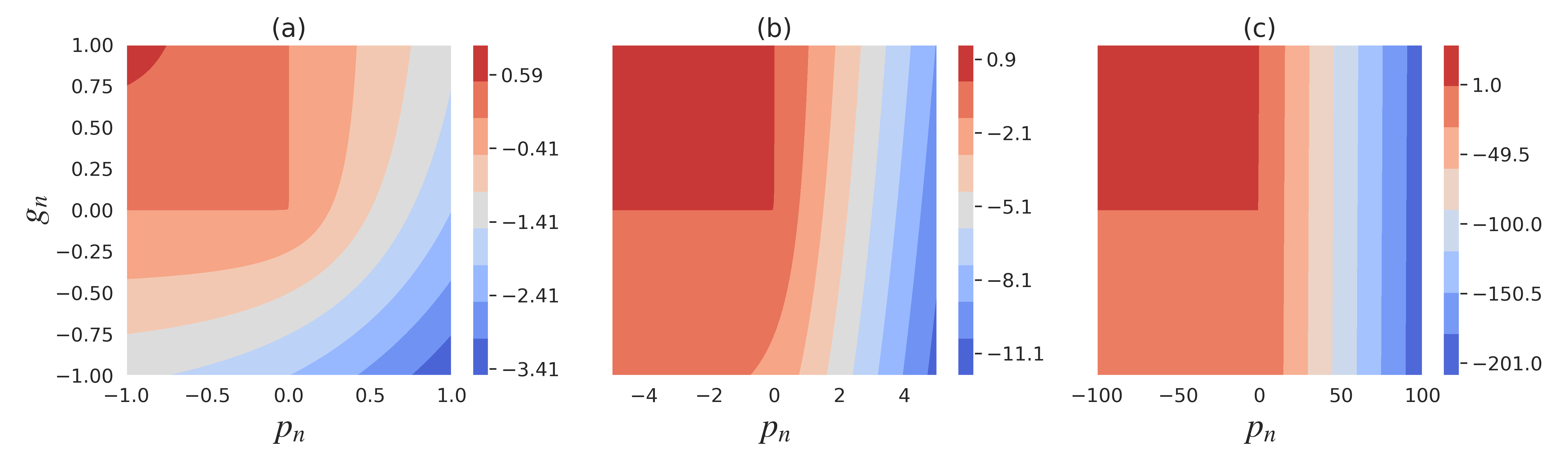}  
    \caption{Domination of $p_n$ over $g_n$: (a) no domination, (b) intermediate domination and (c) strong domination.}
    \label{fig:domination}
\end{figure}

The Fischer-Burmeister NCP function reduces a set of inequality constraints into a single equation. 
However, it might cause domination of one term over another.
As explained in the previous sections, $g_n$ is derived from the displacement field $\boldsymbol{u}$ and $p_n$ is derived from
the Cauchy stress field $\boldsymbol{\sigma}$. Depending on the chosen problem parameters, e.g., stiffness, the estimated quantities
$\boldsymbol{u}$ and $\boldsymbol{\sigma}$
might have different scales, and then so do $g_n$ and $p_n$.
As illustrated in Fig. \ref{fig:domination}, when $g_n$ and $p_n$ have similar scales, there is no domination. 
But increasing $p_n$ causes domination of $p_n$ over $g_n$. In terms of optimization, the neural network will then expand
more effort into minimizing the large-scale quantity $p_n$, which might cause unacceptable violations of constraints related to $g_n$.
To tackle this problem, non-dimensionalization techniques can be used to ensure that $p_n$ and $g_n$ have similar scales \cite{xu2022transfer}. 

\subsubsection{Importance of output scaling}\label{outputscale}
Output scaling is a functionality of output transformation that can prevent the optimizer to get stuck in local minima. 
In the context of solid and contact mechanics, the PINNs estimate the displacement field $\boldsymbol{u}$ and the stress field $\boldsymbol{\sigma}$ as neural network outputs. 
Assuming that no transfer learning is used, the first estimations of outputs are done randomly due to random
initialization of the neural network parameters. Also, 
there are well-known and popular methods that generate initial network estimations following a normal distribution with a mean of zero and a standard
deviation of one \cite{geron2022hands}, e.g., \textit{Glorot} initialization. 
However, using such initialization methods could cause convergence issues because of the output quantities having similar magnitudes. 
This issue can be explained through the example of the SS loss term (see Eq. \ref{eq:pinn_elasticity_detail})

\begin{equation}
    \mathcal{L}_{\mathrm{SS}}= \sum_{i_{n}=N_m+1}^{N_{n}} {w}_{i_{n}} \frac{1}{N_{rp}} \sum_{i_{rp}=1}^{N_{rp}}\Bigl[ \bigl[\boldsymbol{\sigma}\left(\boldsymbol{x}^{i_{rp}} \right) -\mathbb{C} : \bm{\varepsilon}\bigr]_{i_{n}}\Bigr]^2.
\end{equation}

\afterequation
Minimization of $\mathcal{L}_{\mathrm{SS}}$ requires the condition $\boldsymbol{\sigma}=\mathbb{C} : \bm{\varepsilon}$ to hold. 
In case very large values of the material properties are chosen, e.g. Young's modulus, the term $\mathbb{C} : \bm{\varepsilon}$ will initially dominate $\boldsymbol{\sigma}$,
which means a large loss contribution has to be handled by the optimizer. Therefore, 
to minimize the loss, either a significant increase in $\boldsymbol{\sigma}$ or a significant decrease in $\boldsymbol{u}$ is required.
Such large increments in the optimization procedure are troublesome, as the gradient of the employed $\tanh()$ activation function tends to zero for large function arguments, 
which is also referred to as the vanishing gradient problem \cite{Zheng2021}. 
To ease optimization, the network output $\boldsymbol{u}$ 
can be scaled by the inverse of the Young's modulus, 
i.e. by $\frac{1}{E}$. This ensures that the initial magnitudes of both terms in the SS condition are comparable, 
which summarizes the benefits of output scaling in a nutshell.

\section{Numerical examples}\label{sec4}

In the following, we investigate three numerical examples. 
The first example is the well-known Lam\'e problem of elasticity, which is 
considered as a preliminary test without contact to
verify that our PINN framework works as expected, including in particular the 
hard enforcement of DBC and NBC with output transformation.
Afterward, our investigation focuses on examining two contact examples: 
a contact problem between a simple square block and a rigid flat surface, 
and the Hertzian contact problem. The main difference between the two contact examples is the fact
that the actual contact area has to be identified by the PINN in the Hertzian example, 
while the potential and actual contact areas are the same in the case of 
the square block and the rigid flat surface. Note that
2D plane strain conditions are considered throughout the entire section and body forces are neglected.  

All numerical examples have the following common settings.
The PINN maps spatial coordinates $\boldsymbol{x}=(x,y)$ as inputs to transformed mixed-form outputs
$(\pinnUVector,\pinnSVector)=(\pinnUComponent{x},\pinnUComponent{y},\pinnSComponent{xx},\pinnSComponent{yy},\pinnSComponent{xy})$. 
Networks are initialized using the \textit{Glorot uniform} initializer, and the \textit{tanh} function is chosen as activation
function. Models are first trained using the stochastic gradient descent optimizer \textit{Adam}\cite{Kingma2014AdamAM} with a learning rate, 
$lr=0.001$, for 2000 epochs, and then we switch to the limited memory BFGS algorithm including box constraints (\textit{L-BFGS-B})\cite{ciyoulbfgs} until
one of the stopping criteria is met \cite{byrd1995limited} \cite{markidisoptimizer}.  
Our workflow is developed based on the \textit{DeepXDE} package \cite{lu2021deepxde}
and we refer to the \textit{DeepXDE} documentation for default \textit{L-BFGS-B} options.
Note that training points are generated by \textit{GMSH}, since it has strong capabilities
in mesh generation and visualization, and provides boundary normals at arbitrary query points \cite{gmsh}.


As a common error metric, we report the vector-based relative $L_2$ errors for displacement $\relErrorU$ 
and stress fields $\relErrorS$ as follows

\begin{equation}
    \begin{aligned}
    \relErrorU &= \frac{\sqrt{\sum_{i}\sum_{j=1}^{N_{\text{test}}}
    \left(\pinnUComponent{i}(\boldsymbol{x}^j)-\genericSol_{i}(\boldsymbol{x}^j)\right)^2}}
    {\sqrt{\sum_{i}\sum_{j=1}^{N_{\text{test}}}
    \left(\genericSol_{i}(\boldsymbol{x}^j)\right)^2}}
    \quad \text{for } i=(x,y),\\[5pt]
    \relErrorS &= \frac{\sqrt{\sum_{i}\sum_{j=1}^{N_{\text{test}}}
    \left(\pinnSComponent{i}(\boldsymbol{x}^j)-\sigma_{i}(\boldsymbol{x}^j)\right)^2}}
    {\sqrt{\sum_{i}\sum_{j=1}^{N_{\text{test}}}
    \left(\sigma_{i}(\boldsymbol{x}^j)\right)^2}} 
    \quad \text{for } i=(xx,yy,xy).
    \end{aligned}
\end{equation}

\afterequation
Here, $\pinnUComponent{i}$ and $\pinnSComponent{i}$ denote PINN solutions and $u_i$ and $\sigma_i$ denote reference solutions 
that are obtained analytically or numerically. Also, the index  
$j=(1,\ldots, N_{\text{test}})$ runs over the test points that are generated using structured meshes. 
We refer to Appendix \ref{append_error_comparison} for an error comparison between vector-based and integral-based error measurements. 

\subsection{Lam\'e problem of elasticity}\label{sec_lame}
\ifhidden
{}
\fi

\begin{figure}[thbp]
    \centering
    \includegraphics[width=0.7\linewidth]{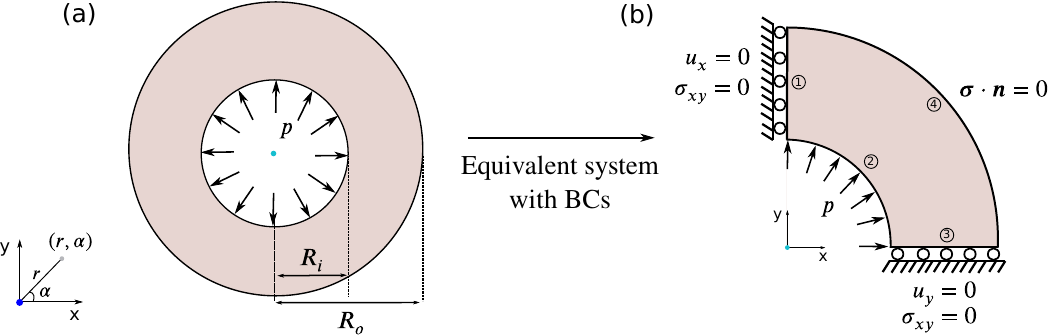}  
    \caption{The Lam\'e problem of elasticity.  (a) The full problem under internal pressure, and (b) the equivalent problem due to the axisymmetrical nature of both geometry and loading, and accompanying boundary conditions.}
    \label{fig:lame_geometry}
\end{figure}

In the first example, we study a benchmark example without contact, namely the well-known Lam\'e problem of a cylinder, subjected to 
an internal pressure $p$ (see Fig. \ref{fig:lame_geometry}a). 
Since the problem is geometrically axisymmetric and the internal pressure is applied to the entire inner boundary, 
only a quarter of the annulus is considered. 
The analytical solution for the stress and displacement field can be derived in polar coordinates \{$r,\alpha$\} as 
\cite{timoshenko1951theory}

\begin{equation}
    \begin{aligned}
        \sigma_{r r} &=\frac{R_{i}^{2} p}{R_{o}^{2}-R_{i}^{2}}\left(1-\frac{R_{o}^{2}}{r^{2}}\right), \\
        \sigma_{\alpha \alpha} &=\frac{R_{i}^{2} p}{R_{o}^{2}-R_{i}^{2}}\left(1+\frac{R_{o}^{2}}{r^{2}}\right), \\
        \sigma_{r \alpha} &=0, \\
        u_{r} &=\frac{(1+\nu) p R_i^2 R_o^2}{E(R_o^2-R_i^2)}\left(1+\frac{r}{R_o^2}\right).
    \end{aligned}
\end{equation}


\afterequation
Here, $R_i$ and $R_o$ are the inner and outer radius of the annulus, respectively, $p$ 
represents the internal pressure applied on the inner radius, $E$ is the Young's modulus and $\nu$ is the Poisson's 
ratio. For our specific setup, we set $R_i=1$, $R_o=2$, $p=1$, $E=2000$, and $\nu=0.3$.
Note that our formulation is based on Cartesian coordinates. Thus, polar transformation is performed to compare results. 





To solve the Lam\'e problem, we deploy a PINN in the mixed variable formulation (see Eq. \ref{eq:loss_elasticity_2d}). 
The following output transformation
is applied to enforce displacement and traction BCs as hard constraints on the edges numbered 1 and 3:
\begin{equation}
    \pinnUComponent{x} = \frac{x}{E} \pazocal{N}_{u_x}, \quad
    \pinnUComponent{y} = \frac{y}{E} \pazocal{N}_{u_y}, \quad
    \pinnSComponent{xy} = xy \pazocal{N}_{\sigma_{xy}}.
\end{equation}

\afterequation
As can be seen in Fig. \ref{fig:lame_geometry}, the constructed output transformation fulfills the displacement boundary conditions at $x=0$ and $y=0$. Also, 
the displacement outputs are scaled by $1/E$ to ease the optimization process (see Section \ref{outputscale}).
Traction boundary conditions are enforced as hard constraints for edges numbered 1 and 3, which lets us represent zero shear stresses there. 
Since traction boundary conditions on edges 2 and 4 contain a coupling of normal and shear stresses, we 
enforce them as soft constraints. 

The employed PINN is a fully connected neural network consisting of 3 hidden layers of 50 neurons each as indicated in Table \ref{tab:lame_table}.
We train the network using 330 training points, of which 262 are located within the domain 
and the remaining 68 points are on the boundary. We refer to  the introductory paragraph of Section \ref{sec4} 
for further settings of both network and optimizer. 

\begin{table}[thbp]
    \centering
    \begin{tabular}{@{}lccccccc@{}}
        & \multicolumn{1}{c}{\begin{tabular}[c]{@{}c@{}}hidden \\ layers\end{tabular}} 
        & \multicolumn{1}{c}{\begin{tabular}[c]{@{}c@{}}no. of training \\ points\end{tabular}} 
        & \multicolumn{1}{c}{\begin{tabular}[c]{@{}c@{}}no. of test \\ points\end{tabular}}
        & \multicolumn{1}{c}{\begin{tabular}[c]{@{}c@{}}training \\ time (s)\end{tabular}}
        & \multicolumn{1}{c}{\begin{tabular}[c]{@{}c@{}}prediction \\ time (s)\end{tabular}}
        & $\relErrorU$ (\%)
        & $\relErrorS$ (\%)
        \\ \midrule
        \begin{tabular}[c]{@{}l@{}}Lam\'e problem \\ of elasticity\end{tabular} 
        & 3x50  
        & 330  
        & 7104
        & 27.15
        & 0.001 
        & 0.017
        & 0.043 
    \end{tabular}
    \caption{The structure of hidden layers, number of training and test points, performance measurements, and errors 
    for the Lam\'e problem of linear elasticity.}
    \label{tab:lame_table}
\end{table}

\begin{figure}[thbp]
    \centering
    \stackinset{c}{-0.9in}{t}{}{(a) \hspace*{4.7cm} (b) \hspace*{4.8cm} (c)}{%
            \includegraphics[width=0.9\linewidth]{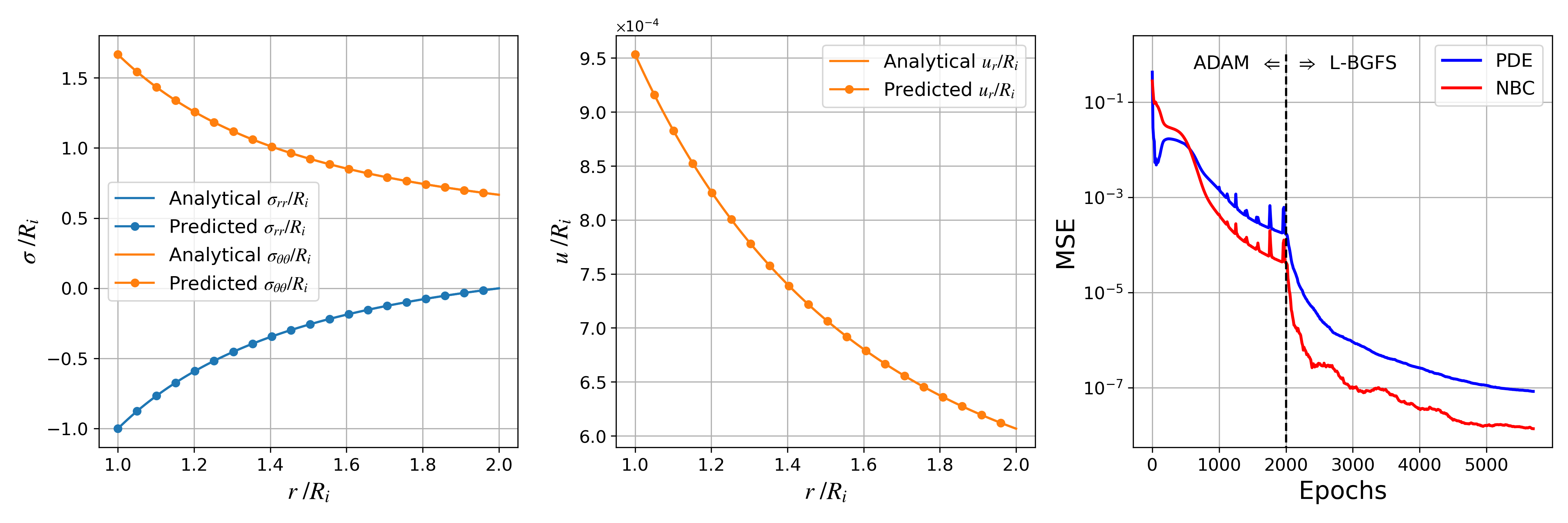}}
    \caption{Results for the Lam\'e problem of linear elasticity. (a) Comparison of normalized stresses obtained from the analytical solution and the predicted values, 
    (b) comparison of normalized displacements, 
    and (c) evolution of the MSE for the summed PDE loss and NBC loss.}
    \label{fig:lame_results}
\end{figure}

Fig. \ref{fig:lame_results}(a,b) show a comparison of the normalized stress and displacement solutions in radial direction.
Relative $L_2$ errors for displacements and stresses are calculated on test points as 0.017\% and 0.043\%, respectively. 
While the network is trained with \textit{Adam}, the convergence rate decreases along with epochs as shown in Fig. \ref{fig:lame_results}c.
Applying \textit{L-BFGS-B} just after \textit{Adam} increases the convergence rates and leads to a further significant reduction of PDE and NBC losses. 
The average MSE for the PDE loss reaches approximately 1.06e-7, while the average MSE for the NBC loss is approximately 1.48e-8 
when all stopping criteria are met. We observe that deploying \textit{Adam} and \textit{L-BFGS-B} optimizers in a sequential order is one of the key points to 
obtain a good accuracy since \textit{Adam} avoids rapid convergence to a local minimum, which has also been mentioned in \cite{markidisoptimizer}.
Using a standard multi-core workstation as hardware, training takes 27.15 \textit{s} and prediction takes 0.001 \textit{s}.

\subsection{Contact between an elastic block and a rigid surface}\label{sec_patch}
\begin{figure}[thbp]
    \centering
    \includegraphics[width=0.65\linewidth]{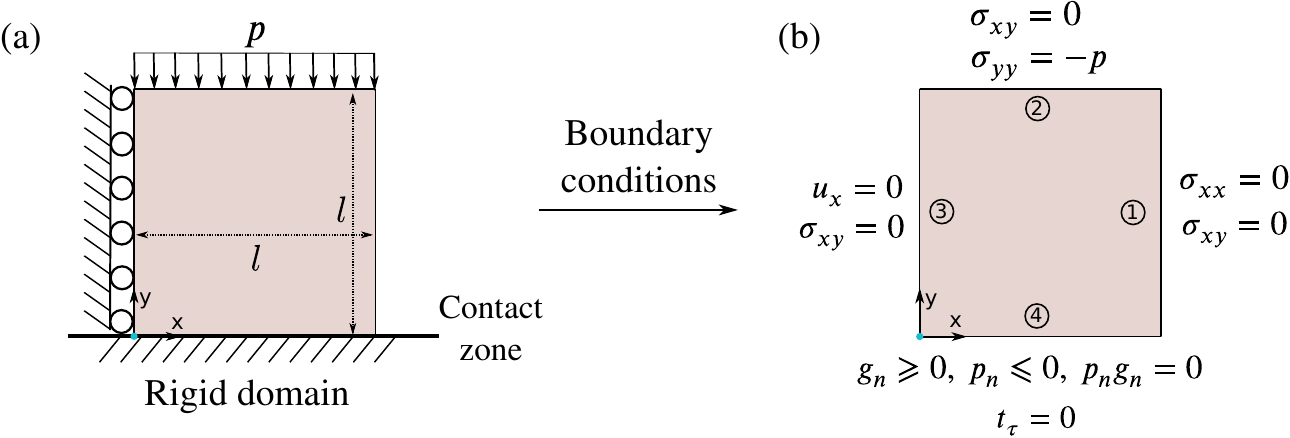}
    \caption{Contact problem between an elastic block and a rigid flat surface,
    (a) geometry, supports and loading,   
    and (b) an equivalent system including all relevant boundary conditions.}
    \label{fig:single_block_geometry}
\end{figure}

\ifhidden
{}
\fi
The second example is a contact problem between a linear-elastic block and a rigid flat surface as depicted in Fig. \ref{fig:single_block_geometry}a.
The elastic block is subjected to an external pressure on its top surface and constrained in the horizontal direction on its left surface.   
The analytical solution \cite{Atluri2006} can easily be derived as follows 
\begin{align*}
    u_x = \frac{-p}{E}\nu (1+\nu) x, & \quad u_y = \frac{p}{E}(1-\nu^2) y, \\
    & \text{and} \\
    \sigma_{yy} = -p, & \quad \sigma_{xy}=0.
\end{align*}

\afterequation
For our specific setup, we set the material parameters as $E=1.33$, $\nu=0.33$, the edge length of the square block as $l$=1 and the pressure as $p$=0.1.

Similar as before, we apply the following output transformation to enforce displacement and traction boundary
conditions on the edges numbered 1, 2 and 3 as
\begin{equation}
    \pinnUComponent{x} = x \pazocal{N}_{u_x}, \quad
    \pinnSComponent{xx} = (l-x) \pazocal{N}_{\sigma_{xx}}, \quad
    \pinnSComponent{yy} = -p + (l-y)\pazocal{N}_{\sigma_{xy}}, \quad
    \pinnSComponent{xy} = x(l-y)(l-x) \pazocal{N}_{\sigma_{xy}}.
\end{equation}

\afterequation
These output transformations can easily be derived based on Fig. \ref{fig:single_block_geometry}b. For instance, 
the normal stress $\sigma_{yy}$ in the loading direction is equal to $-p$ at $y=l$. Thus, we choose $g(\boldsymbol{x})=-p$ 
, and $h(\boldsymbol{x})=(l-y)$ so that the requirements given in Eq. \ref{eq:h_func} are fulfilled.
On the other hand, the contact constraints at the bottom edge are enforced as soft constraints. 
Contact constraints are enforced using the three different 
methods that were explained earlier in Section \ref{sec3_3}: the sign-based method, the Sigmoid-based method and the
\textit{Fischer-Burmeister} NCP function. For the Sigmoid-based method, we choose $\delta_{g_n}=10$ and $\delta_{p_n}=100$. 
For training, 514 points are used
(434 points lie within the domain and 80 points lie on the boundary), and 11827 points are used for testing. 

\begin{figure}[thbp]
    \def\figsize{0.2}
    \newcommand{\trimsize}{18}
    \def\hskiplocal{\hskip -0.75cm}
    \def\hskiplocalerr{\hskip -0.9cm}
    \centering
    \def\arraystretch{0.8}%
    \begin{tabular}{c@{\hskip -0.225cm} c@{\hskiplocalerr}c @{\hskiplocal} c@{\hskiplocalerr}c @{\hskiplocal} c@{\hskiplocalerr}c}
    & $u_{y}$ & $E_{abs}^{u_{y}}$ 
    & $\sigma_{yy}$ & $E_{abs}^{\sigma_{yy}}$ 
    & $\sigma_{xy}$ & $E_{abs}^{\sigma_{xy}}$ \\
    \rownameform{\hspace*{-1.5cm}\textbf{sign-based}}&
    \includegraphics[trim={0cm \trimsize cm 0cm 0cm},clip, width=\figsize\linewidth]{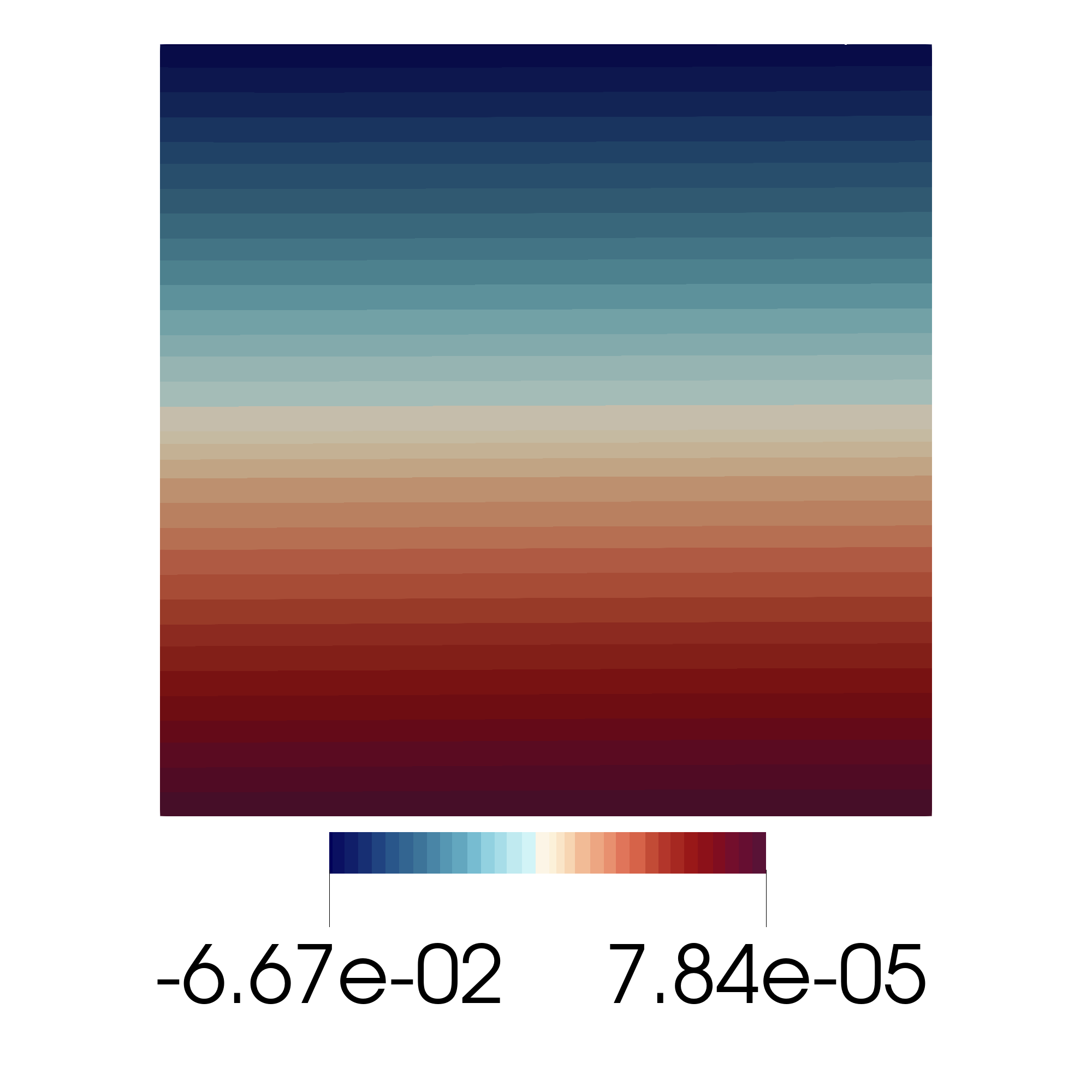}&
    \includegraphics[trim={0cm \trimsize cm 0cm 0cm},clip, width=\figsize\linewidth]{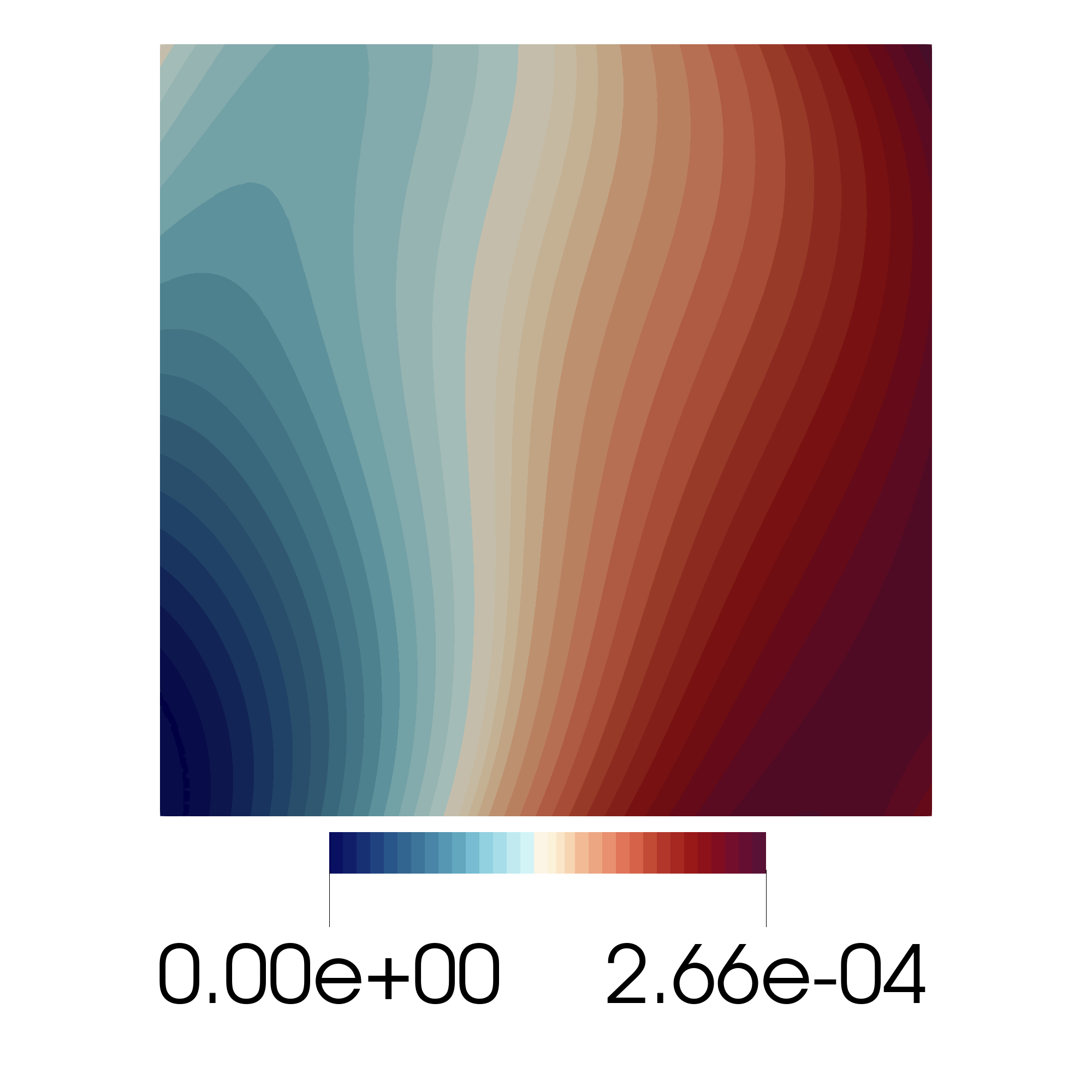}&
    \includegraphics[trim={0cm \trimsize cm 0cm 0cm},clip, width=\figsize\linewidth]{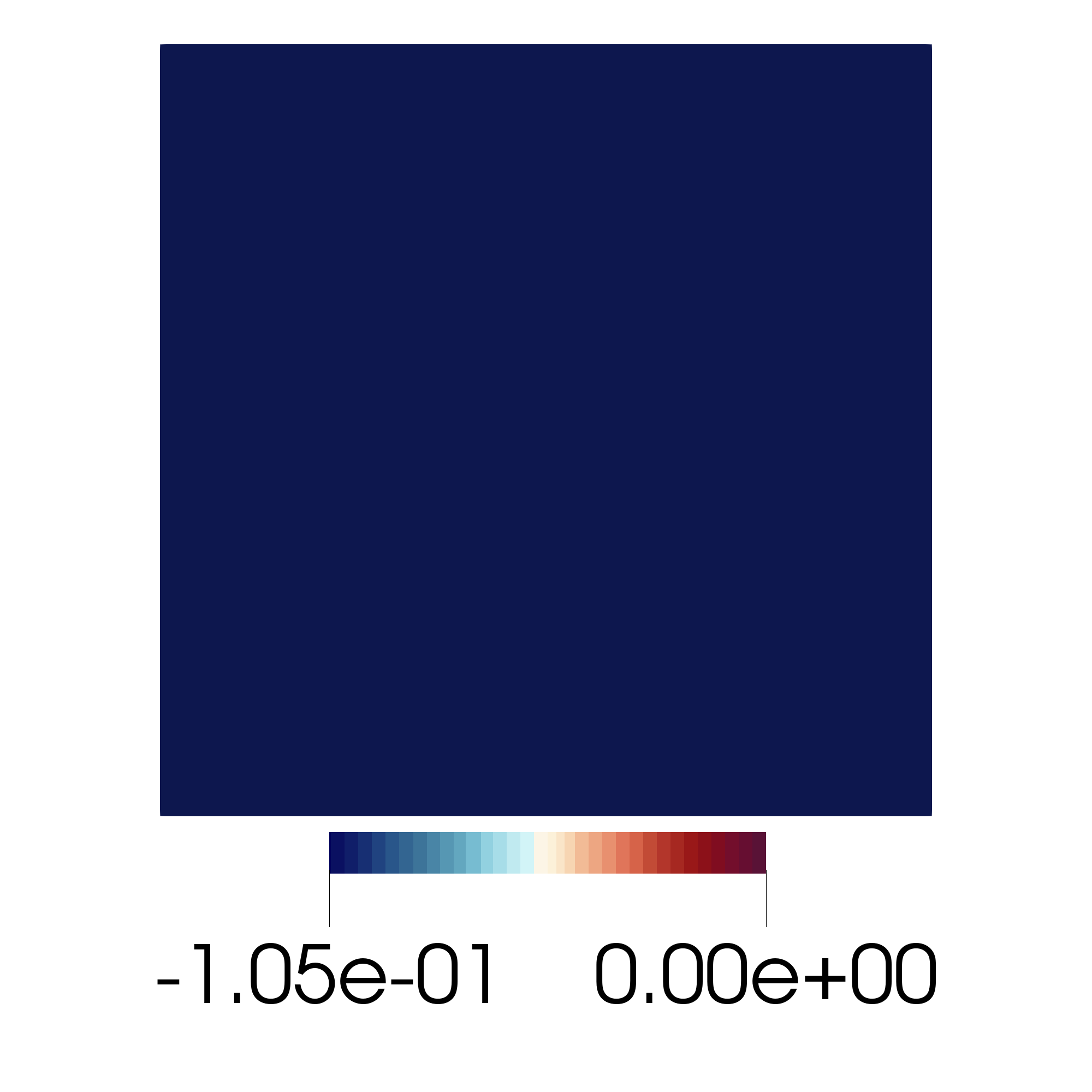}&
    \includegraphics[trim={0cm \trimsize cm 0cm 0cm},clip, width=\figsize\linewidth]{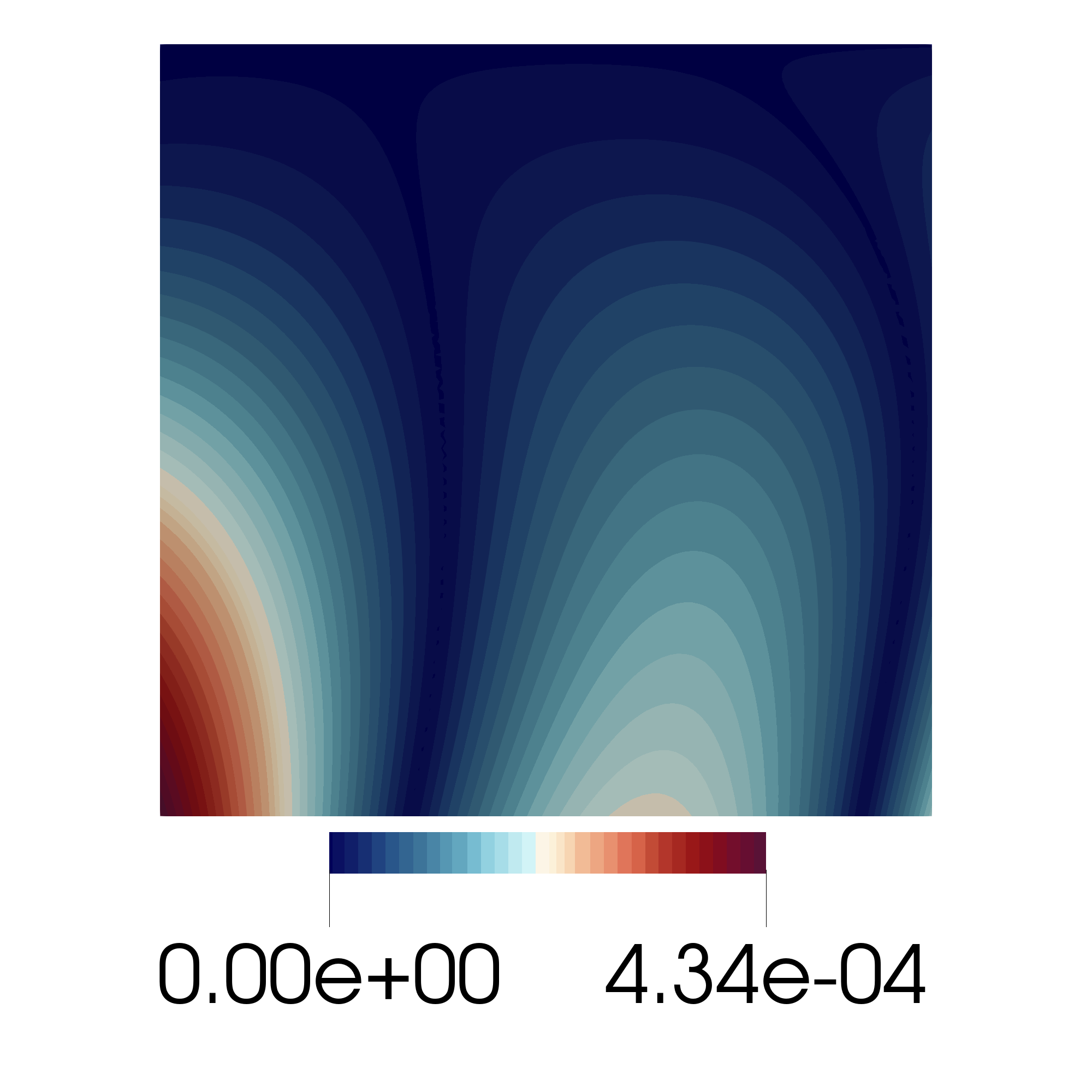}&
    \includegraphics[trim={0cm \trimsize cm 0cm 0cm},clip, width=\figsize\linewidth]{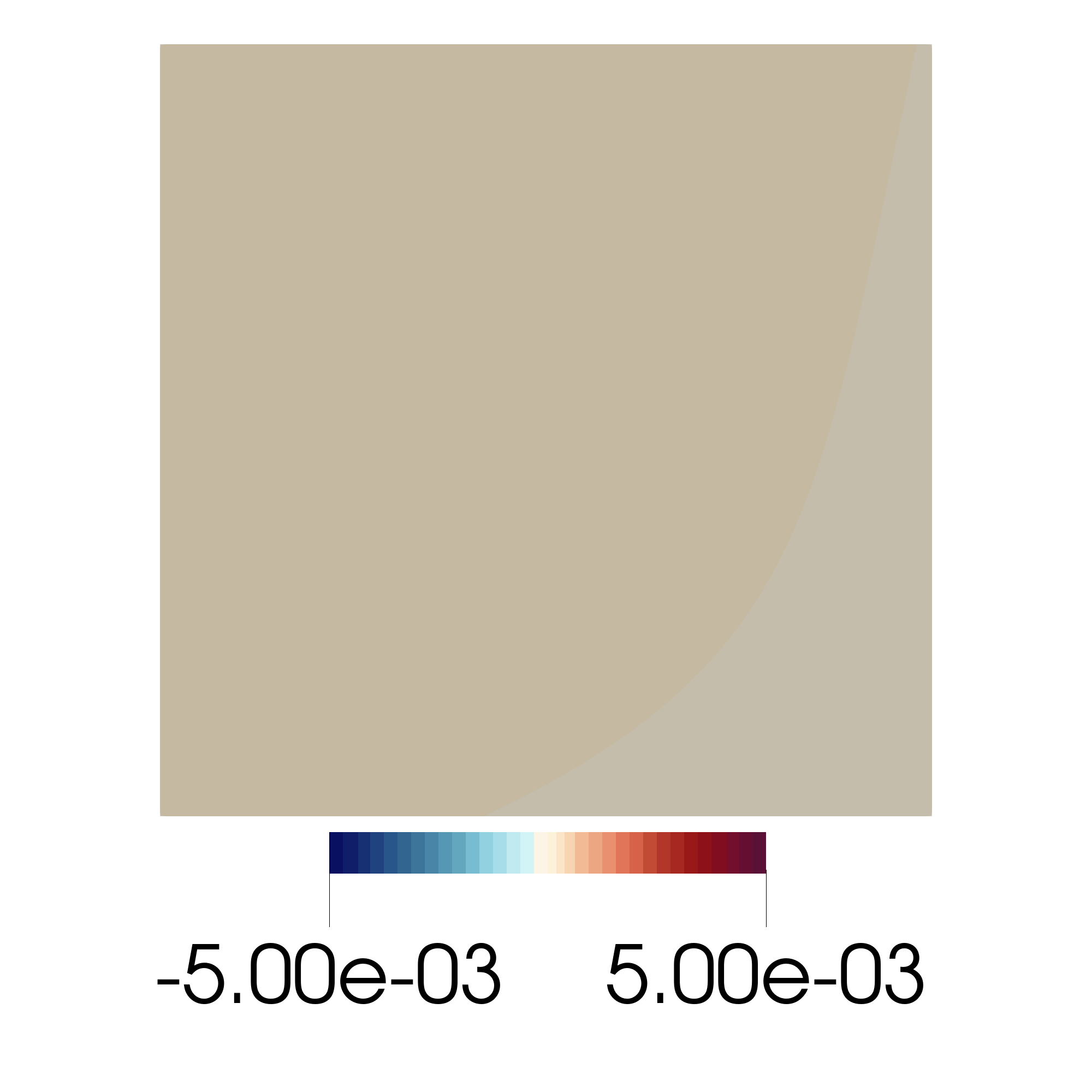}&
    \includegraphics[trim={0cm \trimsize cm 0cm 0cm},clip, width=\figsize\linewidth]{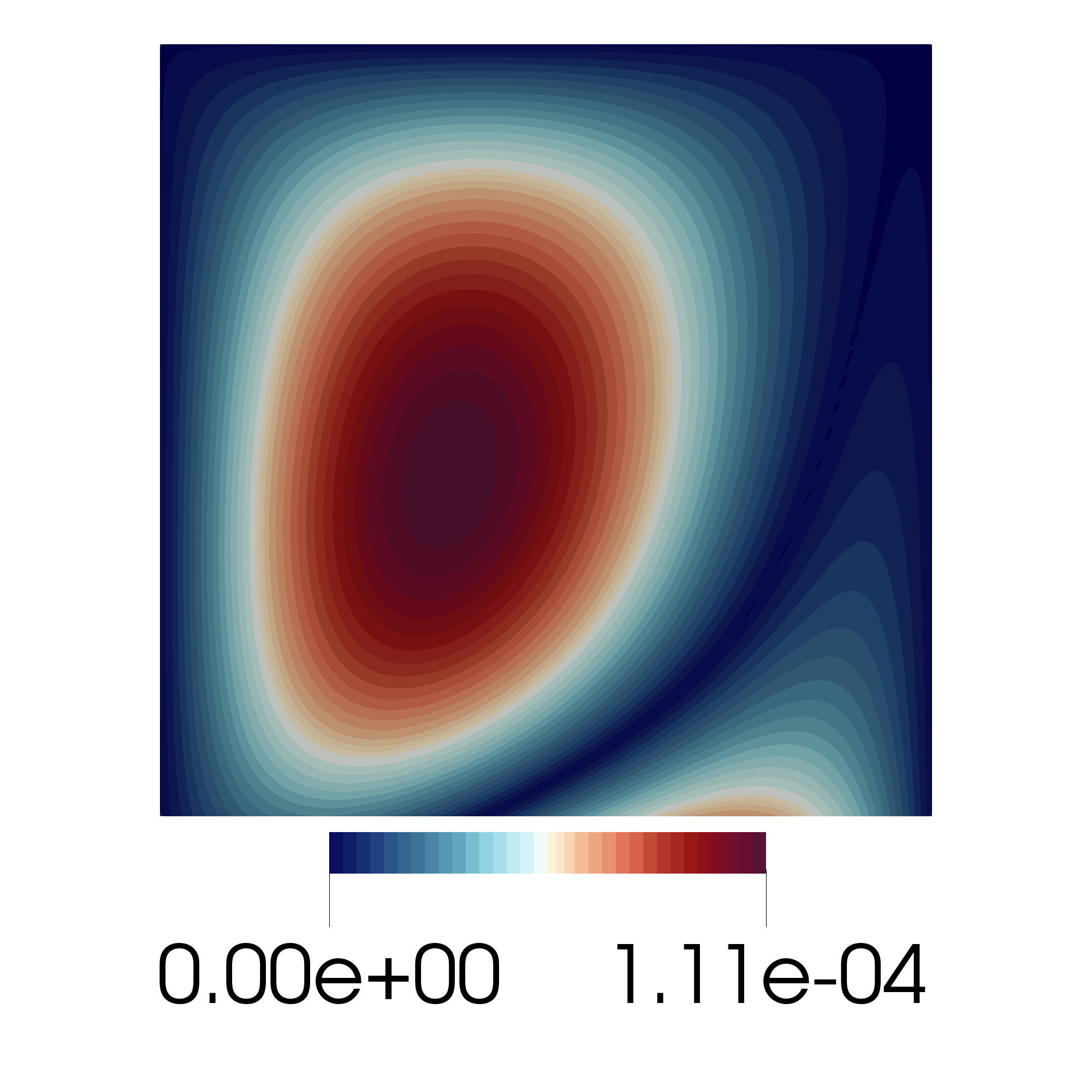}\\[-1ex]
    \rownameform{\hspace*{-1.6cm}\textbf{{Sigmoid-based}}}&
    \includegraphics[trim={0cm \trimsize cm 0cm 0cm},clip, width=\figsize\linewidth]{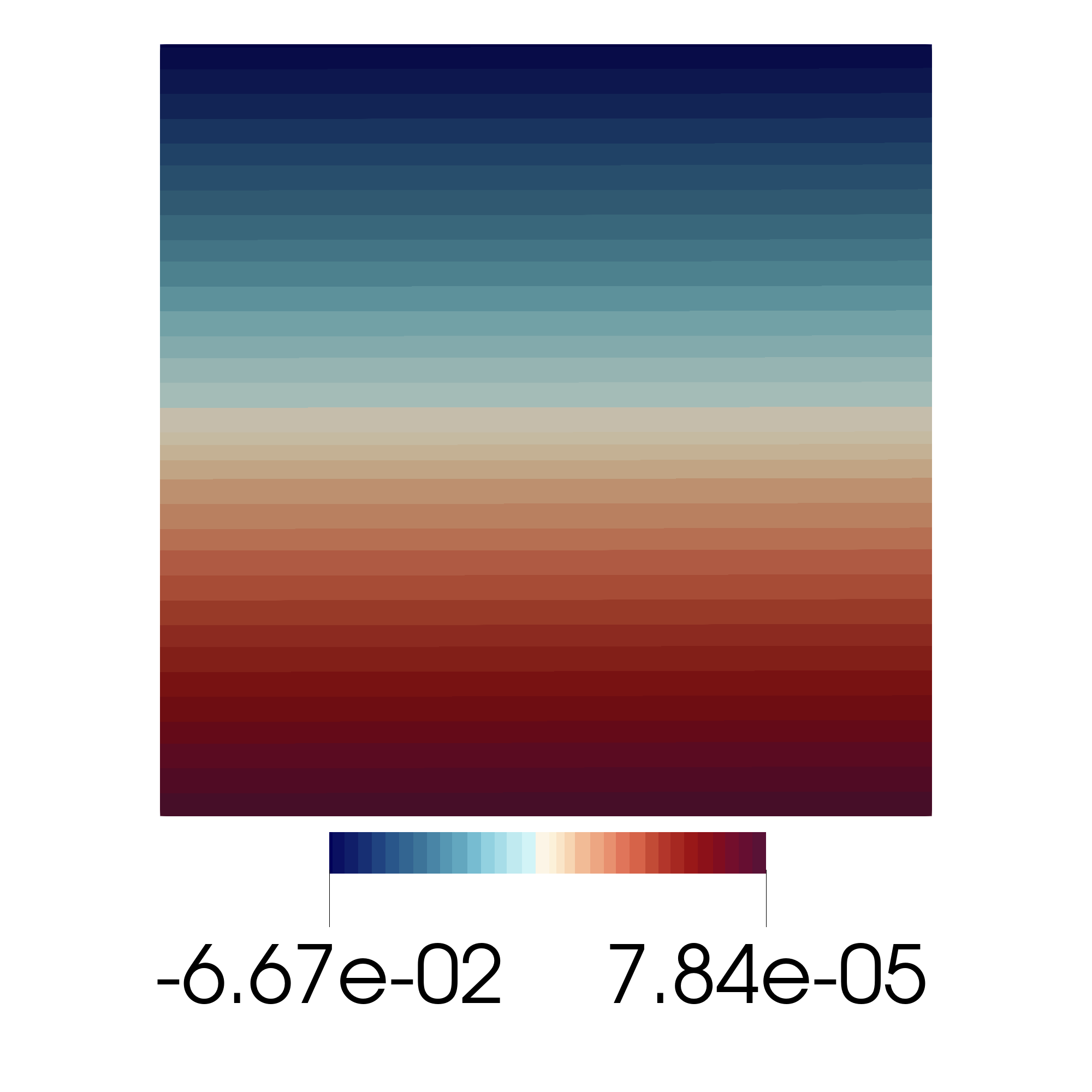}&
    \includegraphics[trim={0cm \trimsize cm 0cm 0cm},clip, width=\figsize\linewidth]{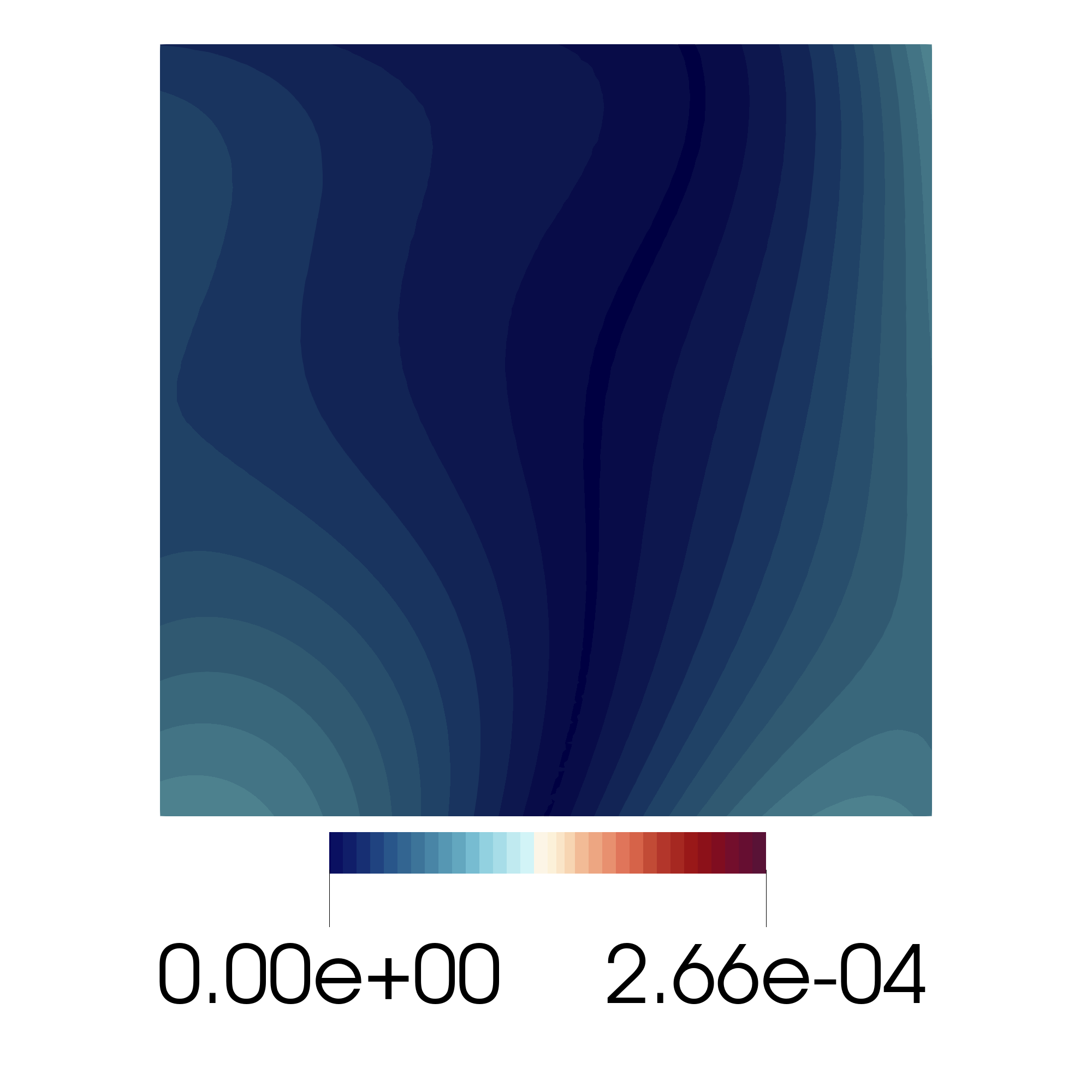}&
    \includegraphics[trim={0cm \trimsize cm 0cm 0cm},clip, width=\figsize\linewidth]{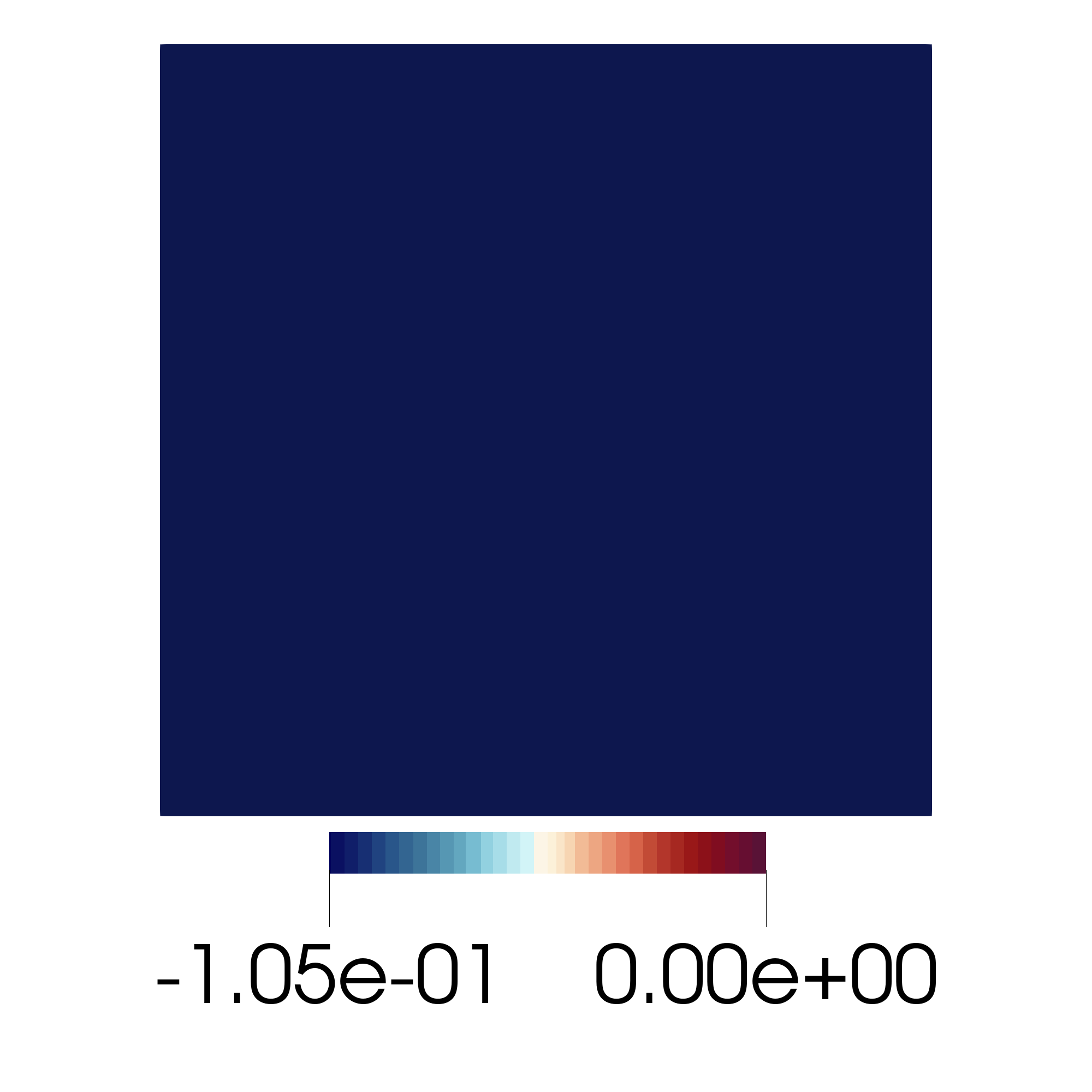}&
    \includegraphics[trim={0cm \trimsize cm 0cm 0cm},clip, width=\figsize\linewidth]{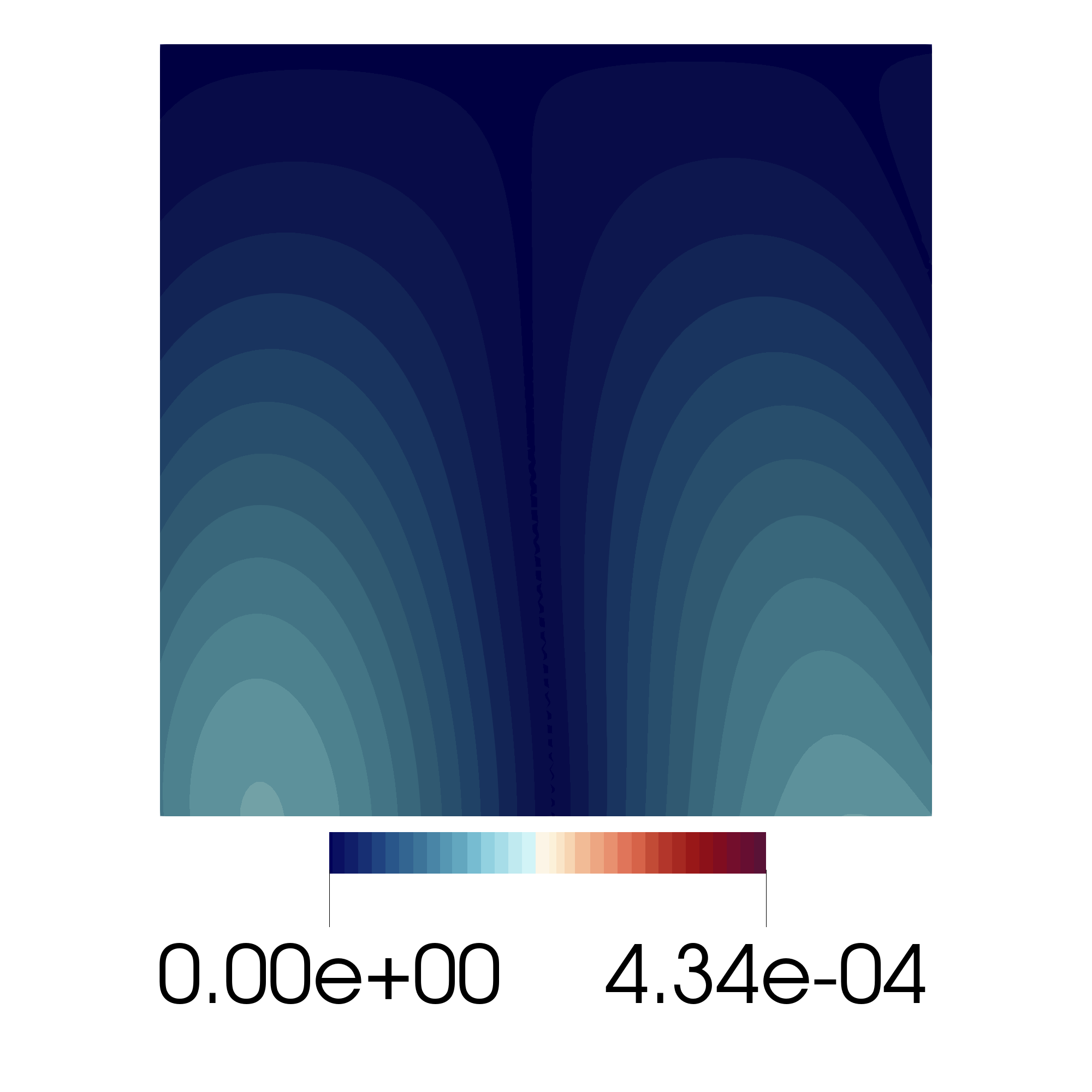}&
    \includegraphics[trim={0cm \trimsize cm 0cm 0cm},clip, width=\figsize\linewidth]{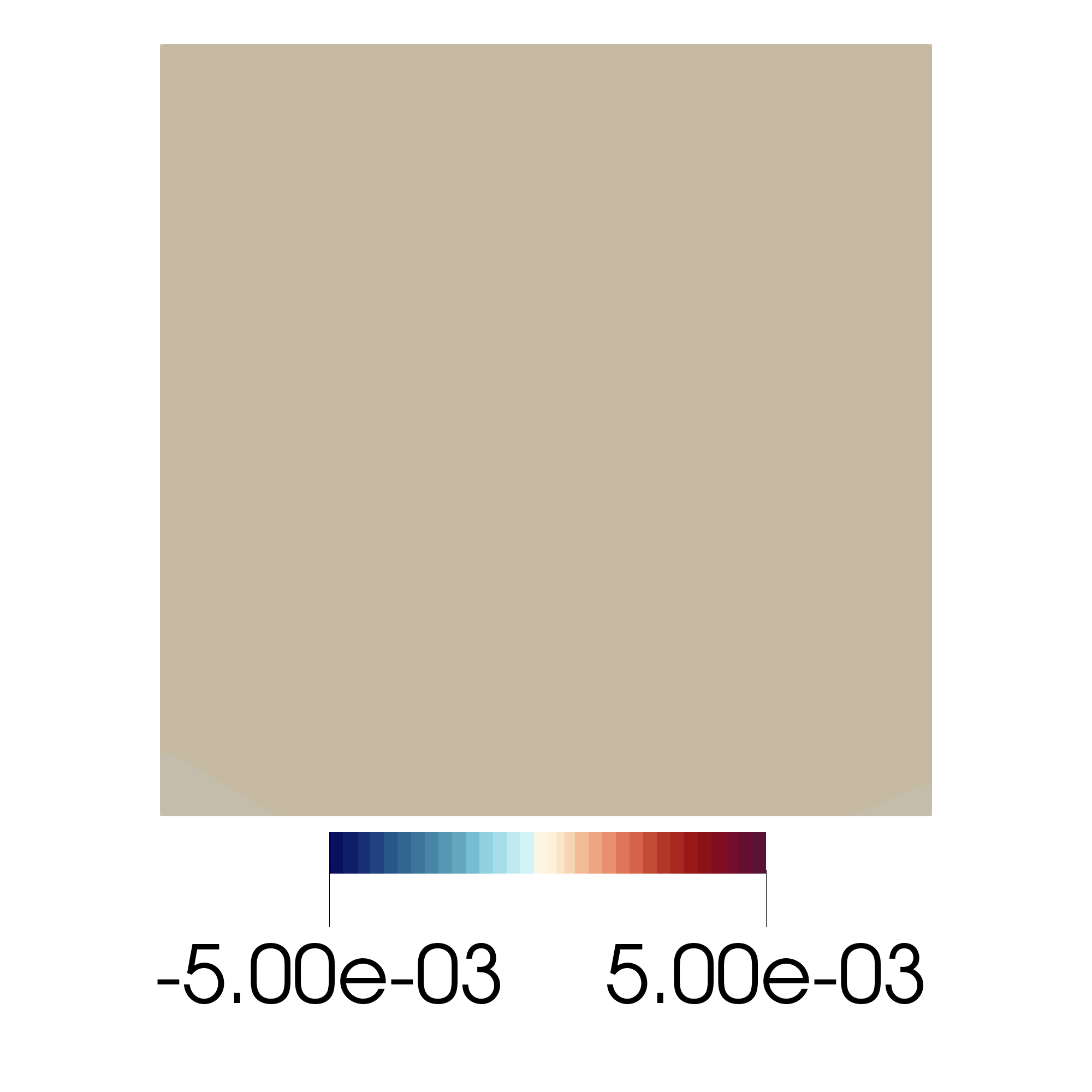}&
    \includegraphics[trim={0cm \trimsize cm 0cm 0cm},clip, width=\figsize\linewidth]{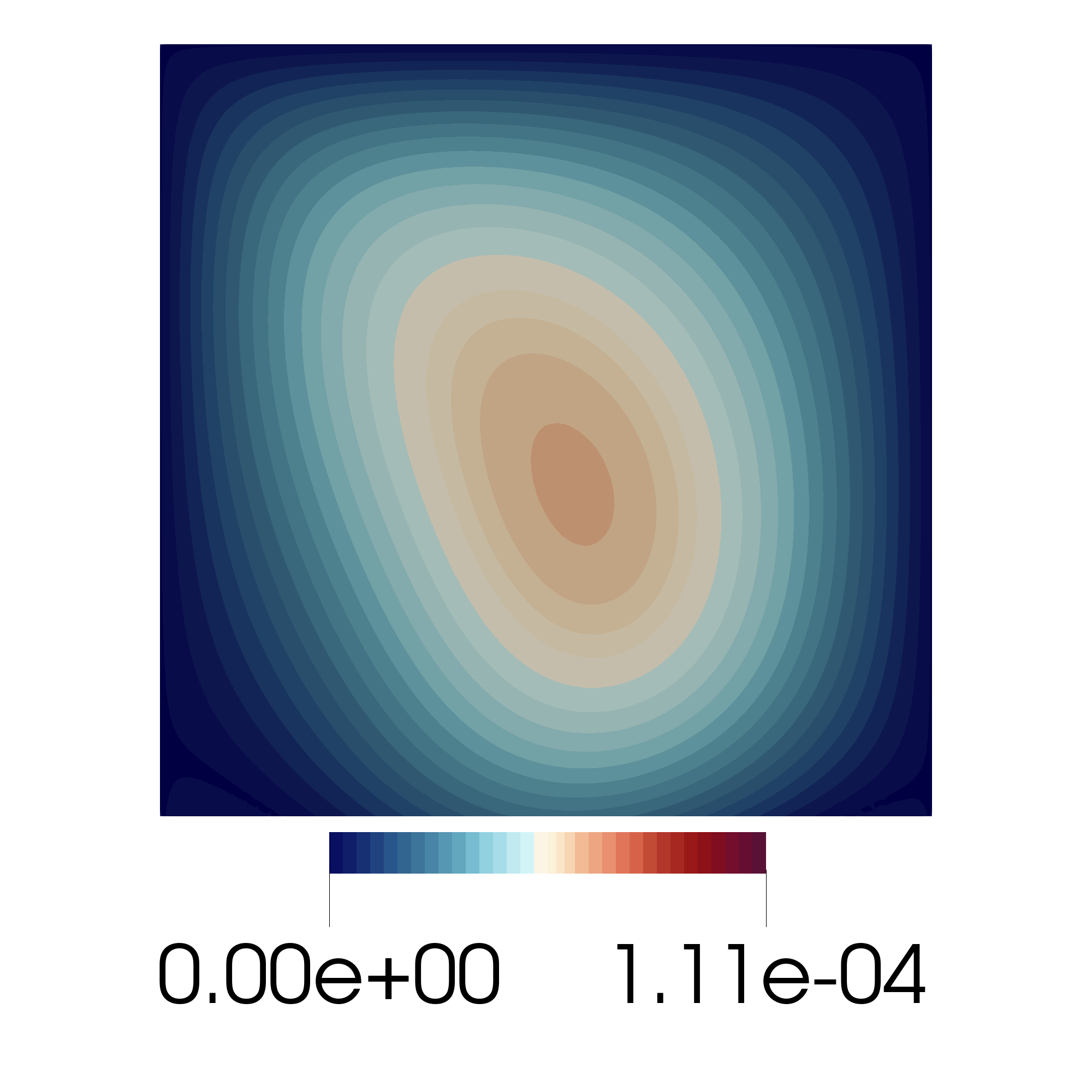}\\[-1ex]
    \rownameform{\hspace*{-0.5cm}\textbf{\textit{Fischer-Burmeister}}}&
    \includegraphics[width=\figsize\linewidth]{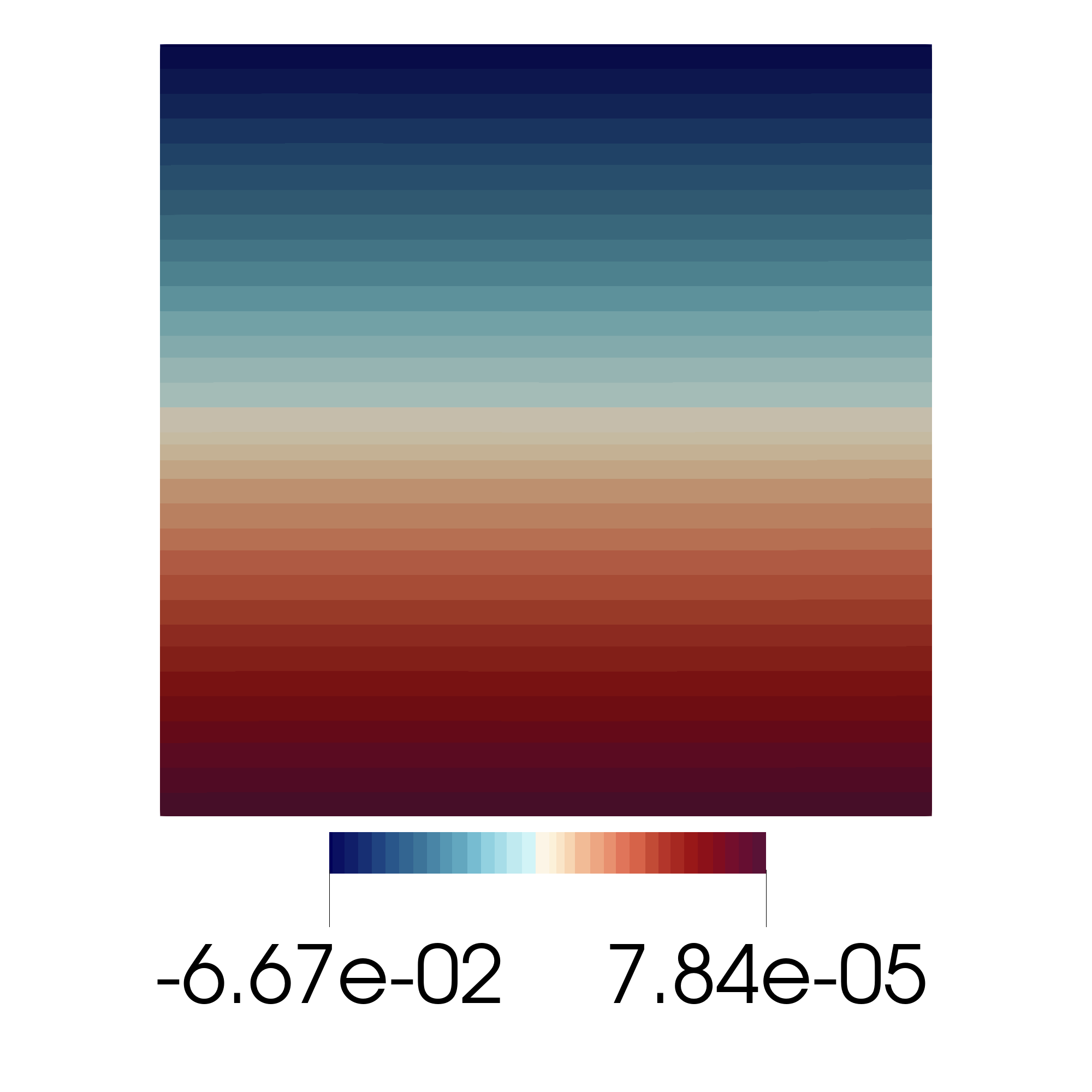}&
    \includegraphics[width=\figsize\linewidth]{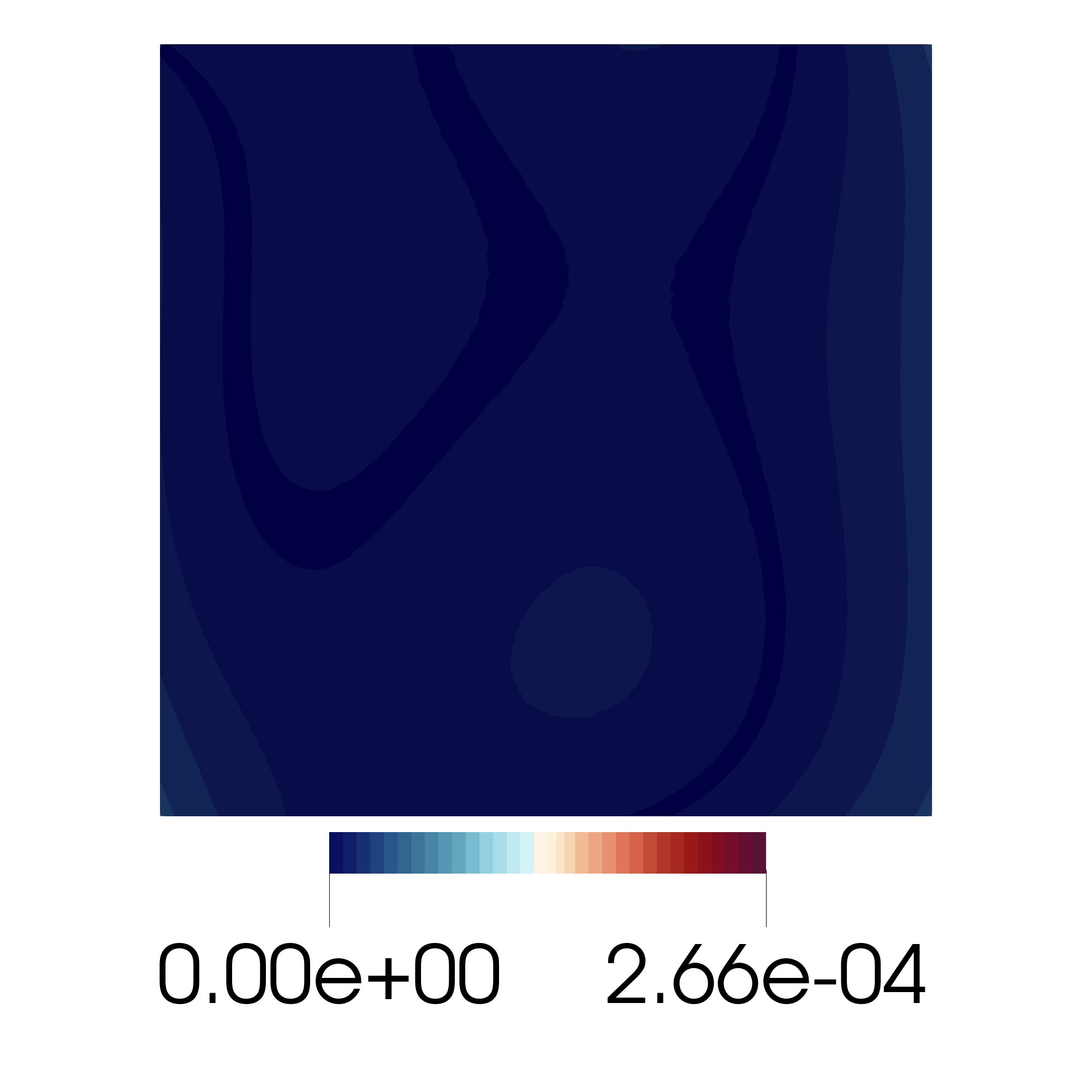}&
    \includegraphics[width=\figsize\linewidth]{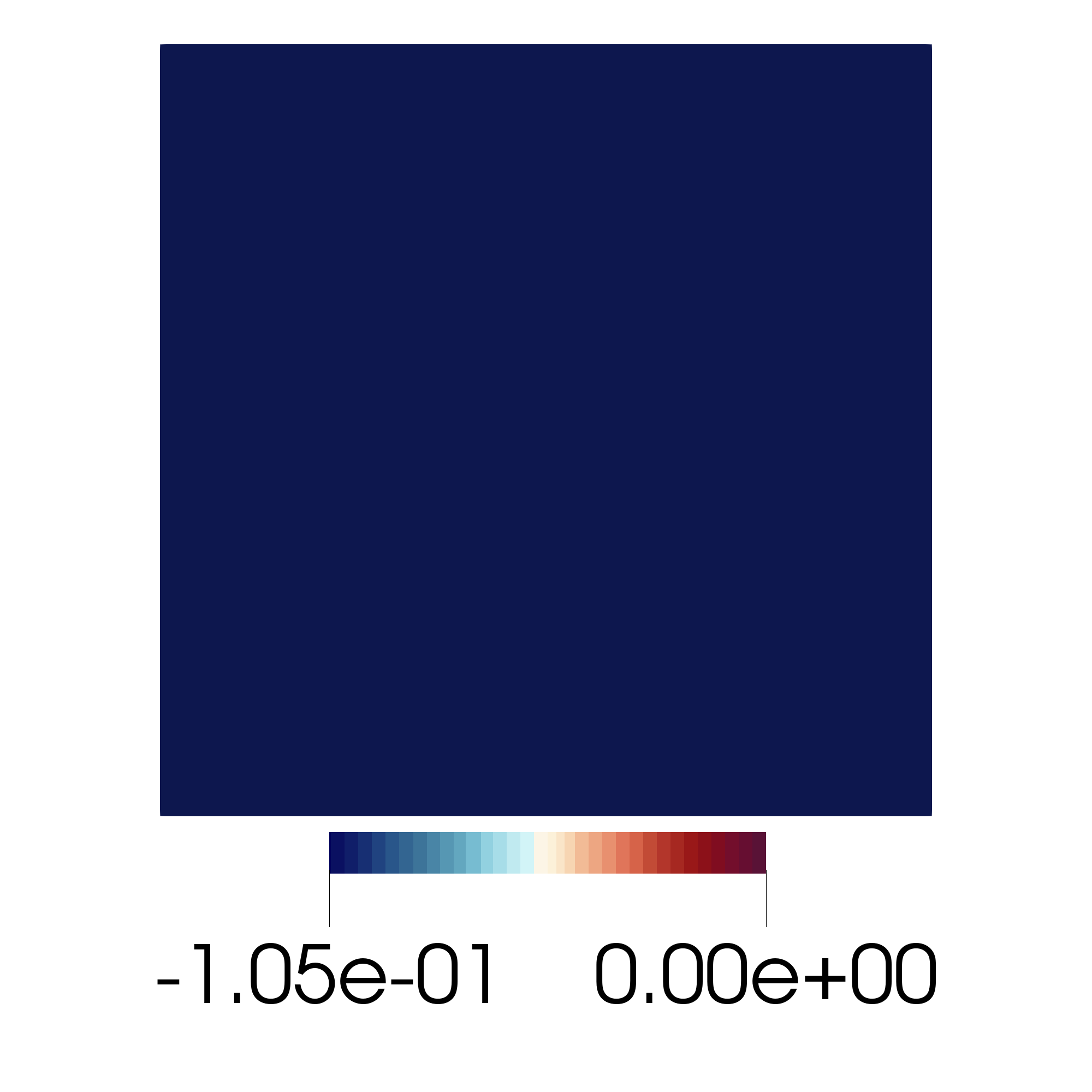}&
    \includegraphics[width=\figsize\linewidth]{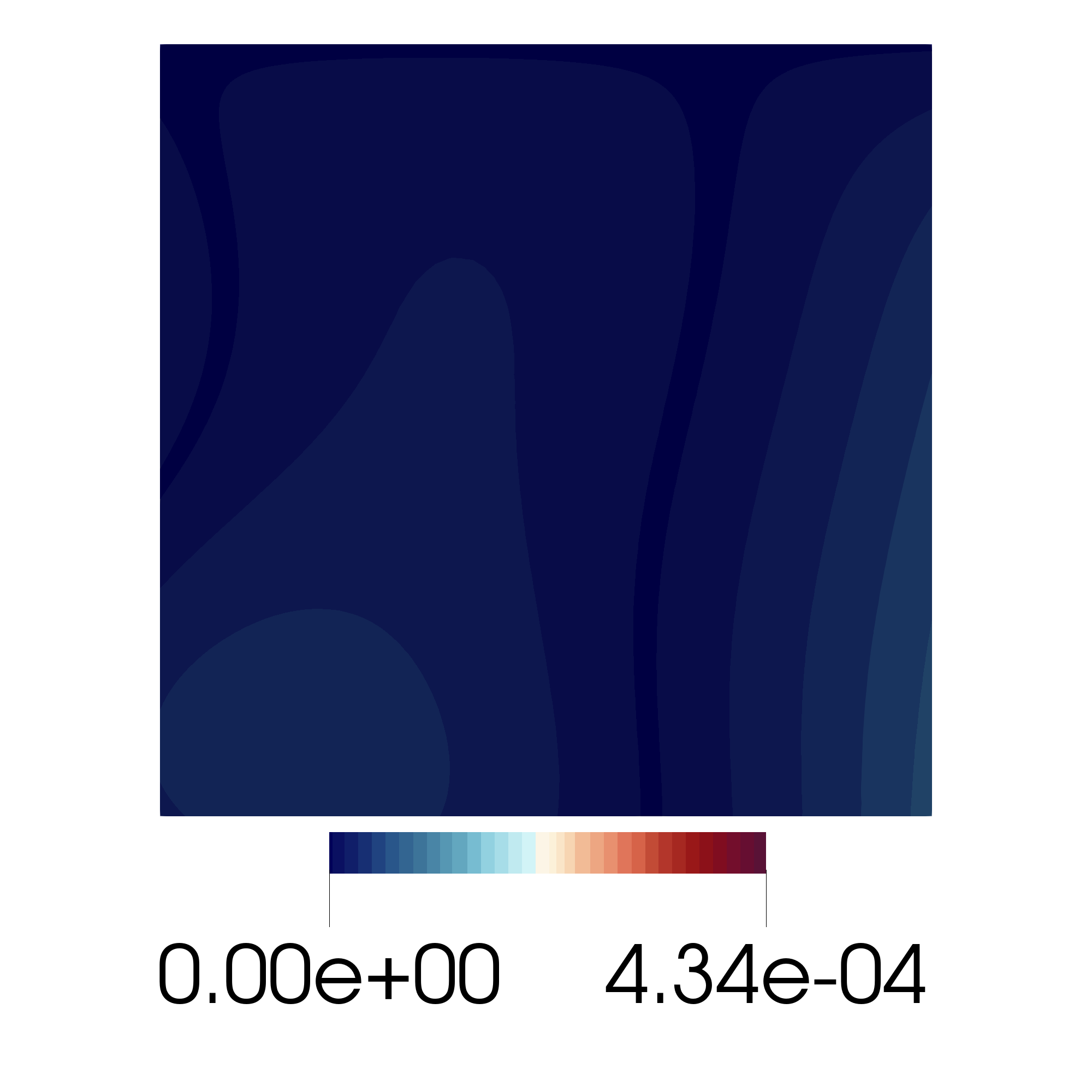}&
    \includegraphics[width=\figsize\linewidth]{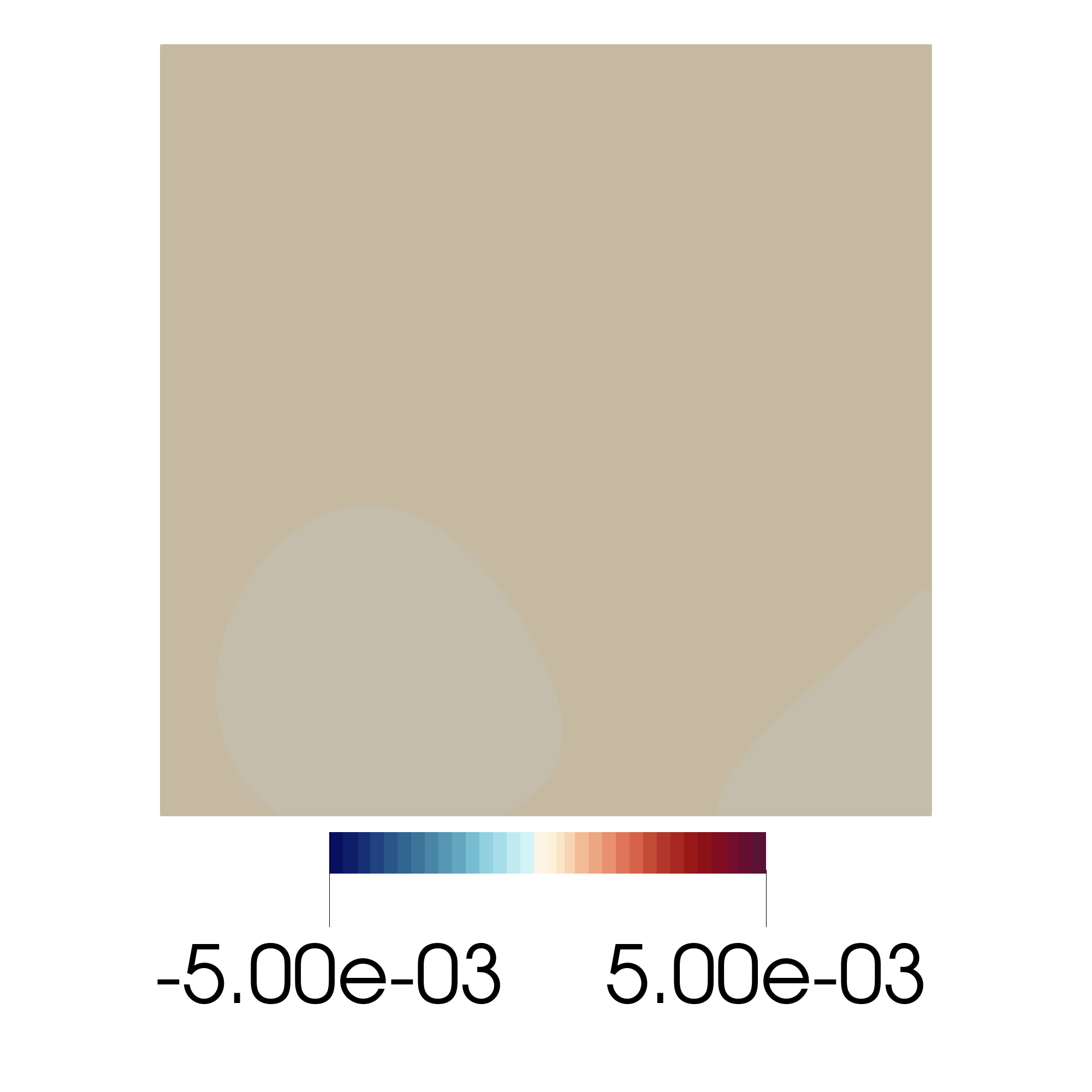}&
    \includegraphics[width=\figsize\linewidth]{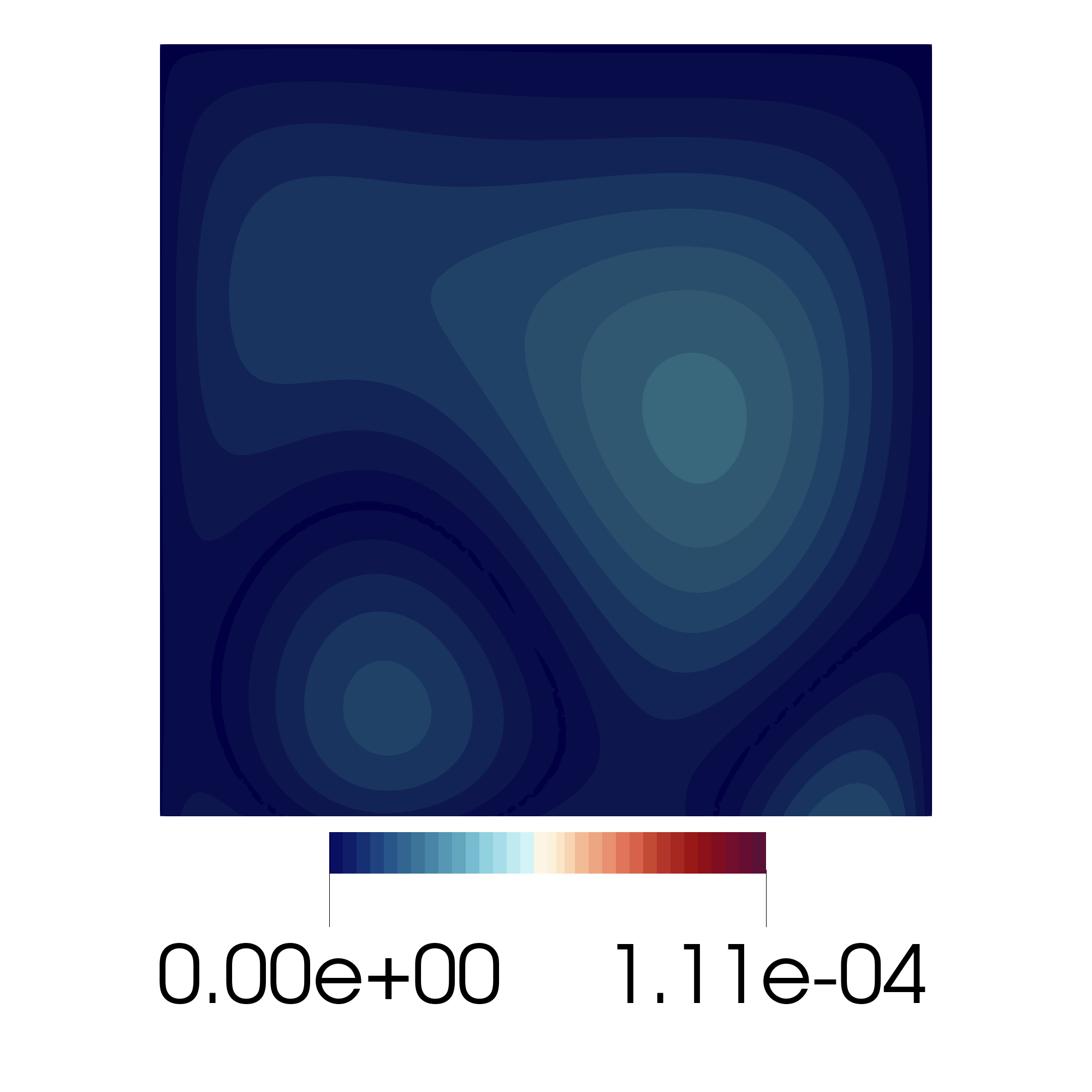}\\[-1ex]
\end{tabular}
    \caption{Predictions of the displacement component $u_y$, stresses $\sigma_{yy}$ and $\sigma_{xy}$ with corresponding absolute errors 
    $E_{abs}^{u_y}$, $E_{abs}^{\sigma_{yy}}$, $E_{abs}^{\sigma_{xy}}$ . The contact constraints 
    $\pazocal{L}_{\mathrm{KKT}}$ are enforced via three different methods: {sign}-based method, {Sigmoid}-based method and
    the \textit{Fischer-Burmeister} NCP function. Absolute error $E_{abs}^{*}=\text{abs}(\tilde{*} - *)$.} %
    \label{fig:single_block_results}
\end{figure}

As illustrated in Fig. \ref{fig:single_block_results}, all of the investigated methods correctly capture that $u_y$ is linearly distributed and 
close to zero at the bottom due to the soft enforcement of contact constraints.
The maximum absolute error for the displacement component $u_y$, denoted by $E_{abs,max}^{u_y}$, is larger for the sign-based formulation compared 
to the two other methods (see Table \ref{tab:rectangle_contact}). Meanwhile, the normal stress component, $\sigma_{yy}$, is close to -0.1 and the shear stress component $\sigma_{xy}$ is close to zero. 
Since the traction boundary condition in y-direction on the top surface and the shear stress boundary conditions are enforced as hard constraints, 
absolute errors in the corresponding regions are zero for all cases.
The sign-based formulation also performs worst in terms of the errors $E_{abs,max}^{\sigma_{yy}}$ and $E_{abs,max}^{\sigma_{xy}}$. 
Moreover, $L_2$ relative errors show that the \textit{Fischer-Burmeister} NCP function performs 
best, i.e. with errors being up to one order of magnitude smaller than for the sign-based and Sigmoid-based variants.  
Additionally, it is observed that all investigated methods require similar computing time for training and prediction.
Overall, it can be concluded that the \textit{Fischer-Burmeister} NCP function yields the best results in terms of accuracy and
computing time.

\begin{table}[thbp]
    \centering
    \begin{tabular}{@{}lcccccccc@{}}
        & \multicolumn{1}{c}{\begin{tabular}[c]{@{}c@{}}hidden \\ layers\end{tabular}} 
        & \multicolumn{1}{c}{\begin{tabular}[c]{@{}c@{}}training \\ time (s)\end{tabular}}
        & \multicolumn{1}{c}{\begin{tabular}[c]{@{}c@{}}prediction \\ time (s)\end{tabular}}
        & $\relErrorU$ (\%)
        & $\relErrorS$ (\%)
        & $E_{abs,max}^{u_y}$ 
        & $E_{abs,max}^{\sigma_{yy}}$ 
        & $E_{abs,max}^{\sigma_{xy}}$
        \\ \midrule
        sign-based
        & 5x50
        & 20.15 
        & 0.388 
        & 0.382 
        & 0.154  
        & 2.66e-4
        & 4.34e-4
        & 1.11e-4\\
        Sigmoid-based
        & 5x50
        & 20.14 
        & 0.389 
        & 0.090      
        & 0.094 
        & 8.23e-5
        & 1.51e-4
        & 6.30e-5\\ 
        \textit{Fischer-Burmeister}
        & 5x50 
        & 18.20 
        & 0.383 
        & 0.024      
        & 0.031 
        & 2.71e-5
        & 6.10e-5
        & 2.55e-5
    \end{tabular}
    \caption{The structure of hidden layers, performance measurements, and errors 
    for the contact example between an elastic block and a rigid surface. Three different methods to enforce
    the KKT constraints are compared.}
    \label{tab:rectangle_contact}
\end{table}

\subsection{Hertzian contact problem}\label{sec_hertz}
\begin{figure}[thbp]
    \centering
    \includegraphics[width=0.65\linewidth]{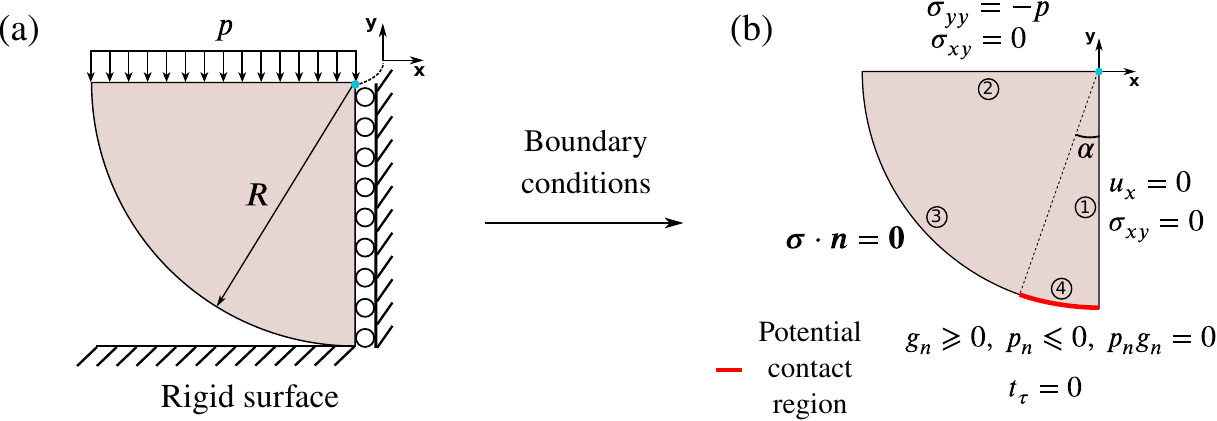}
    \caption{The Hertzian contact problem between an elastic half-cylinder and a rigid flat surface 
    (a) domain under uniform pressure on top and making use of symmetry,
    (b) accompanying boundary conditions.}
    \label{fig:hertzian_geom}
\end{figure}

\ifhidden
{}
\fi
In this example, we consider a long linear elastic
half-cylinder ($E=200$, $\nu=0.3$) lying on a rigid flat surface 
and being subjected to a uniform pressure $p=0.5$ on its top surface as shown in Fig. \ref{fig:hertzian_geom}a. 
The analytical solution for the contact pressure $p_c$ is given as \cite{kikuchi1988contact} \cite{popp2009finite}

\begin{equation}
    p_c = \frac{4 R p}{\pi b^2}\sqrt{b^2 - x^2} \quad \text{with} \  b=2\sqrt{\frac{2R^2 p (1-\nu^2)}{E\pi}}.
\end{equation}

\afterequation
Here, $b$ is the width of the contact zone, and $R$ is the radius of the cylinder. 
For the chosen set of parameters ($p=0.5$, $E=200$, $\nu=0.3$, $R=1$), $b$ can be calculated as 0.076. 
An analytical solution is available only for the contact pressure. Reference solutions for displacement and stress fields 
within the half-cylinder are obtained with well-established FEM algorithms. 
The FEM simulations are performed with our in-house multi-physics research code BACI \cite{baci}.   

The following output transformation is applied to enforce displacement and traction BCs as hard constraints on the edges numbered 1 and 2
(see Fig. \ref{fig:hertzian_geom}b)

\begin{equation}
    \pinnUComponent{x} = \frac{-x}{E} \pazocal{N}_{u_x}, \quad
    \pinnUComponent{y} = \frac{1}{E} \pazocal{N}_{u_y}, \quad
    \pinnSComponent{yy} = -p + (-y) \pazocal{N}_{\sigma_{yy}}, \quad 
    \pinnSComponent{xy} = xy \pazocal{N}_{\sigma_{xy}}.
\end{equation}

\afterequation
As for the Lam\'e problem, we scale the displacement field with the inverse of the Young's modulus $1/E$.
Additionally, only one non-zero helper function $g(x)=-p$ is required for $\pinnSComponent{yy}$.
The traction BC on edge 3 and the contact constraints on edge 4 are enforced as soft constraints. 
To define the potential contact area, we set $\alpha=15^{\circ}$ (corresponding to $b$=0.259), which extends 
well beyond the actual contact area.
Moreover, KKT constraints are enforced using the \textit{Fischer-Burmeister} method. 

In the following, we investigate four distinct PINN application cases for the Hertzian contact problem:
\textit{Case 1}: PINNs as pure forward model / PDE solver, 
\textit{Case 2}: PINNs as data-enhanced forward model, 
\textit{Case 3}: PINNs as inverse solver for parameter identification, and
\textit{Case 4}: PINNs as a fast-to-evaluate surrogate model.
We refer to Table \ref{tab:hertz_table} for information on network architecture, and training and test points.
To facilitate reproduction of the results, the table also reports the weights of individual loss terms. 


\begin{table}[thbp]
    \centering
    \def\arraystretch{1.3}
    \begin{tabular}{@{}lccccc@{}}
    & \multicolumn{1}{c}{\begin{tabular}[c]{@{}c@{}}hidden \\ layers\end{tabular}} 
    & \multicolumn{1}{c}{\begin{tabular}[c]{@{}c@{}}no. of\\ input neurons\end{tabular}} 
    & \multicolumn{1}{c}{\begin{tabular}[c]{@{}c@{}}no. of training \\ points\end{tabular}} 
    & \multicolumn{1}{c}{\begin{tabular}[c]{@{}c@{}}no. of test\\ points\end{tabular}}
    & loss weights
    \\ \midrule
    \textbf{case 1} 
    & 5x50     
    & 2         
    & 47935 
    & 26185
    & $w^{\mathrm{KKT}}$=$10^3$
    \\
    \multirow{2}{*}{\textbf{case 2,3}}
    & \multirow{2}{*}{5x50}   
    & \multirow{2}{*}{2}         
    & \multirow{2}{*}{47935}   
    & \multirow{2}{*}{26185} 
    & \multirow{2}{*}{\def\arraystretch{1.0}\begin{tabular}[c]{@{}c@{}}$w^{\mathrm{\mathrm{KKT}}}$=$10^3$, $w^{\mathrm{EXPs}}_1$=$10^4$, $w^{\mathrm{EXPs}}_2$=$10^4$, \\ $w^{\mathrm{EXPs}}_3$=$10^{-1}$, $w^{\mathrm{EXPs}}_4$=$10^{-1}$, $w^{\mathrm{EXPs}}_5$=$10^{-1}$  \end{tabular}} 
    \\
    &     
    &         
    &      
    &
    &
    \\ 
    \textbf{case 4} 
    & 8x75      
    & 3           
    & 1946 
    & 604    
    & $w^{\mathrm{KKT}}$=$10^4$\  
    \\
    \end{tabular}
    \caption{The structure of hidden layers, number of input neurons, training and test points, and loss weights 
    are given for different cases of the Hertzian contact problem. The remaining loss weights are set as $1$.
    See Appendix \ref{append_2d_loss} for a detailed explanation of loss weights.}
    \label{tab:hertz_table}
\end{table}

\subsubsection{Case 1: PINNs as pure forward model / PDE solver}\label{sec_case_1}
\ifhidden
{}
\fi

\begin{figure}[thbp]
    \def\figsize{0.18}
    \newcommand{\trimsize}{15.3}
    \def\hskiplocal{\hskip 0.1cm}
    \centering
    \def\arraystretch{0.5}%
    \begin{tabular}{c@{ }c@{\hskiplocal}c@{\hskiplocal}c@{\hskiplocal}c@{\hskiplocal}c}
    &$u_x$ &$u_y$ &$\sigma_{xx}$ &$\sigma_{yy}$ &$\sigma_{xy}$  \\
    \rownameform{\hspace*{-1.5cm}\textbf{PINN}}&
    \includegraphics[trim={0cm \trimsize cm 0cm 0cm},clip, width=\figsize\linewidth]{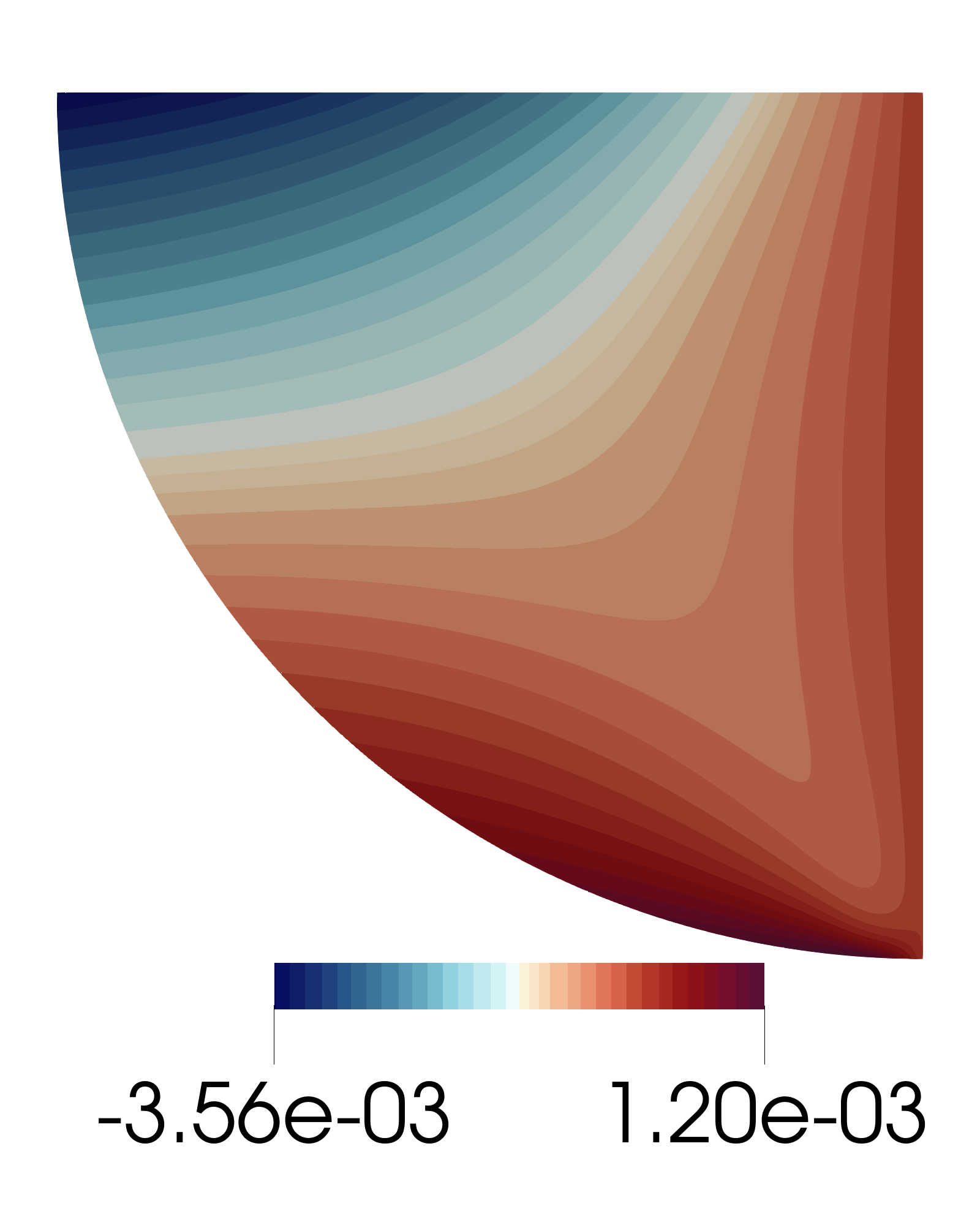}&
    \includegraphics[trim={0cm \trimsize cm 0cm 0cm},clip, width=\figsize\linewidth]{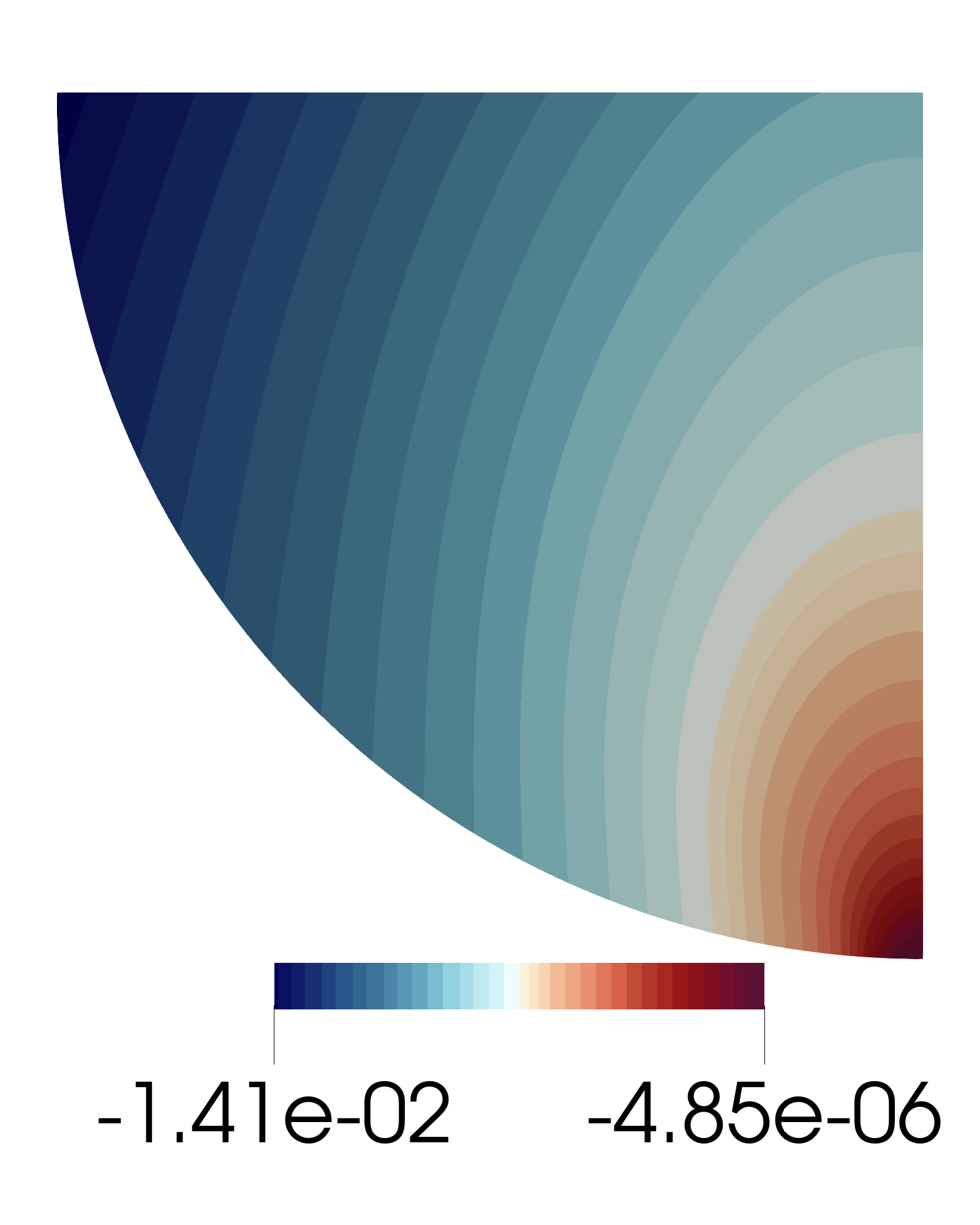}&
    \includegraphics[trim={0cm \trimsize cm 0cm 0cm},clip, width=\figsize\linewidth]{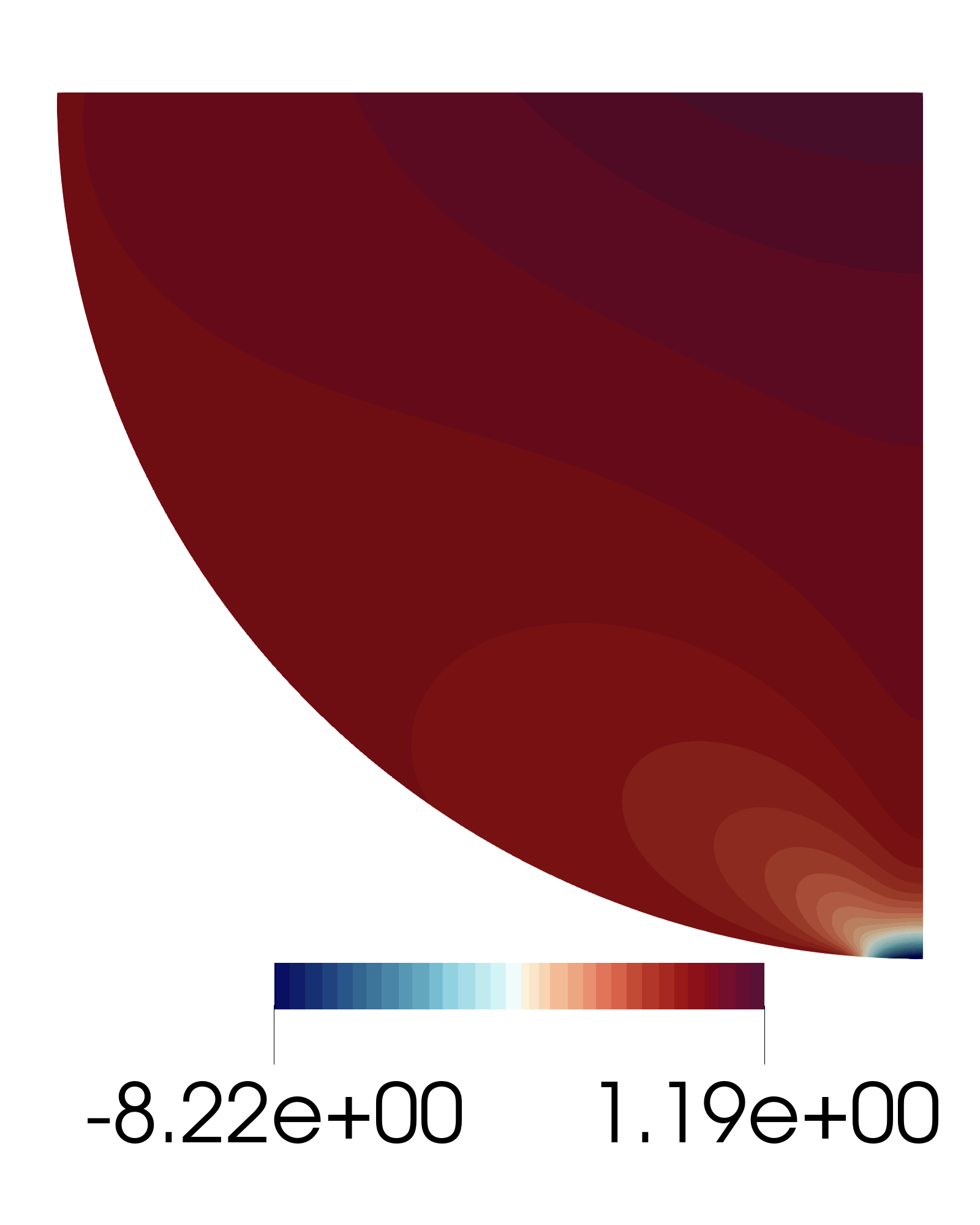}&
    \includegraphics[trim={0cm \trimsize cm 0cm 0cm},clip, width=\figsize\linewidth]{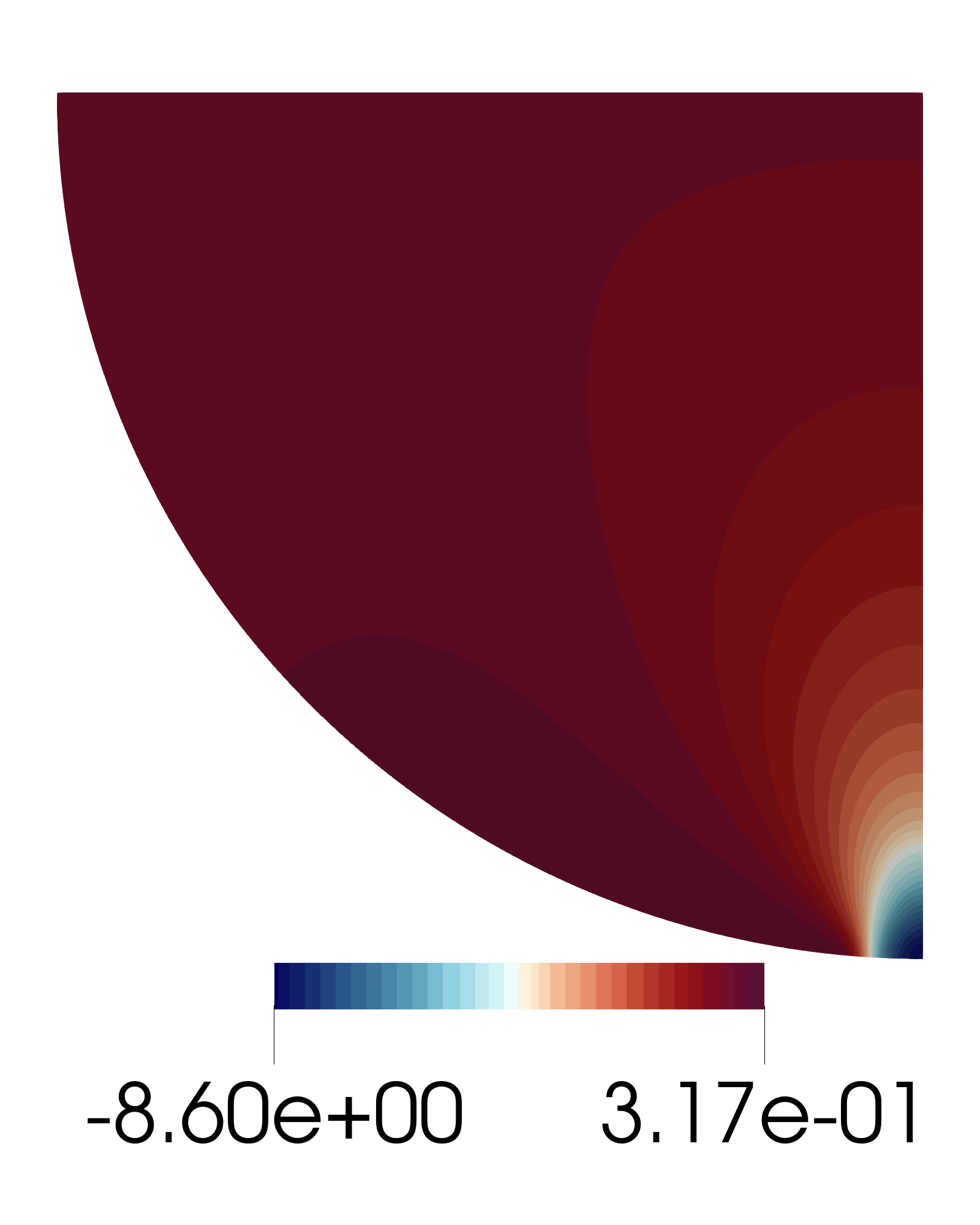}&
    \includegraphics[trim={0cm \trimsize cm 0cm 0cm},clip, width=\figsize\linewidth]{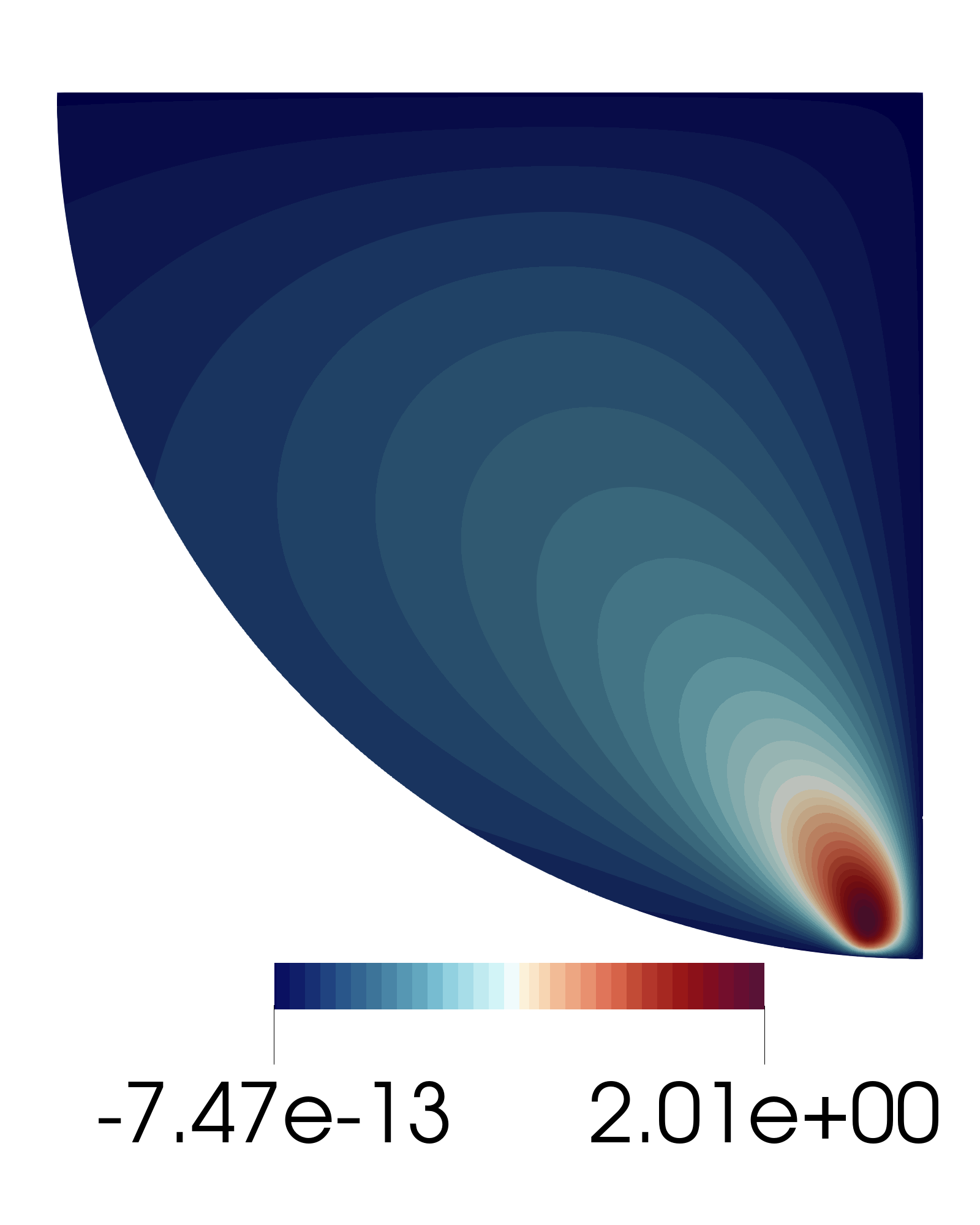}
    \\[-1ex]
    \rownameform{\textbf{FEM}}&
    \includegraphics[width=\figsize\linewidth]{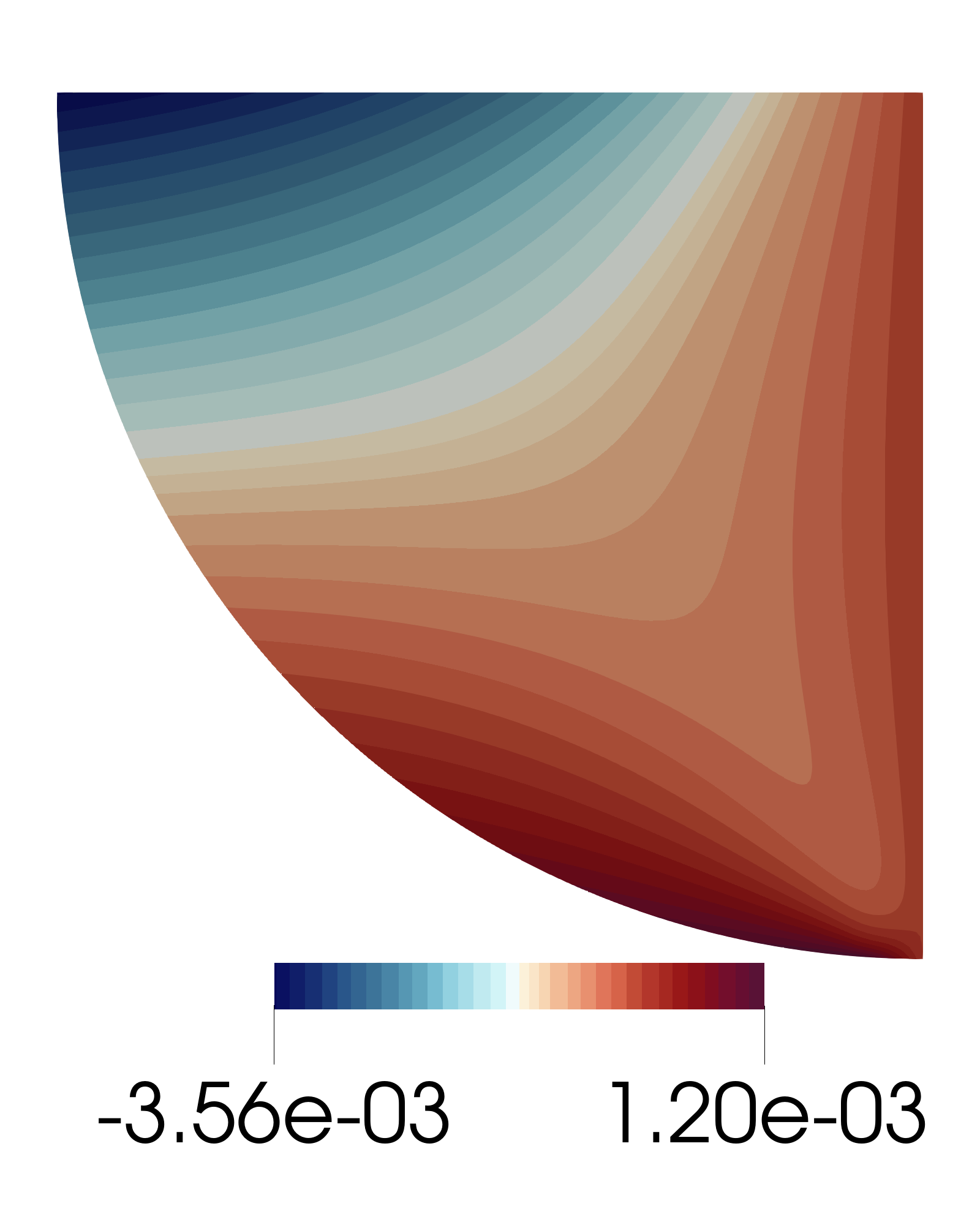}&
    \includegraphics[width=\figsize\linewidth]{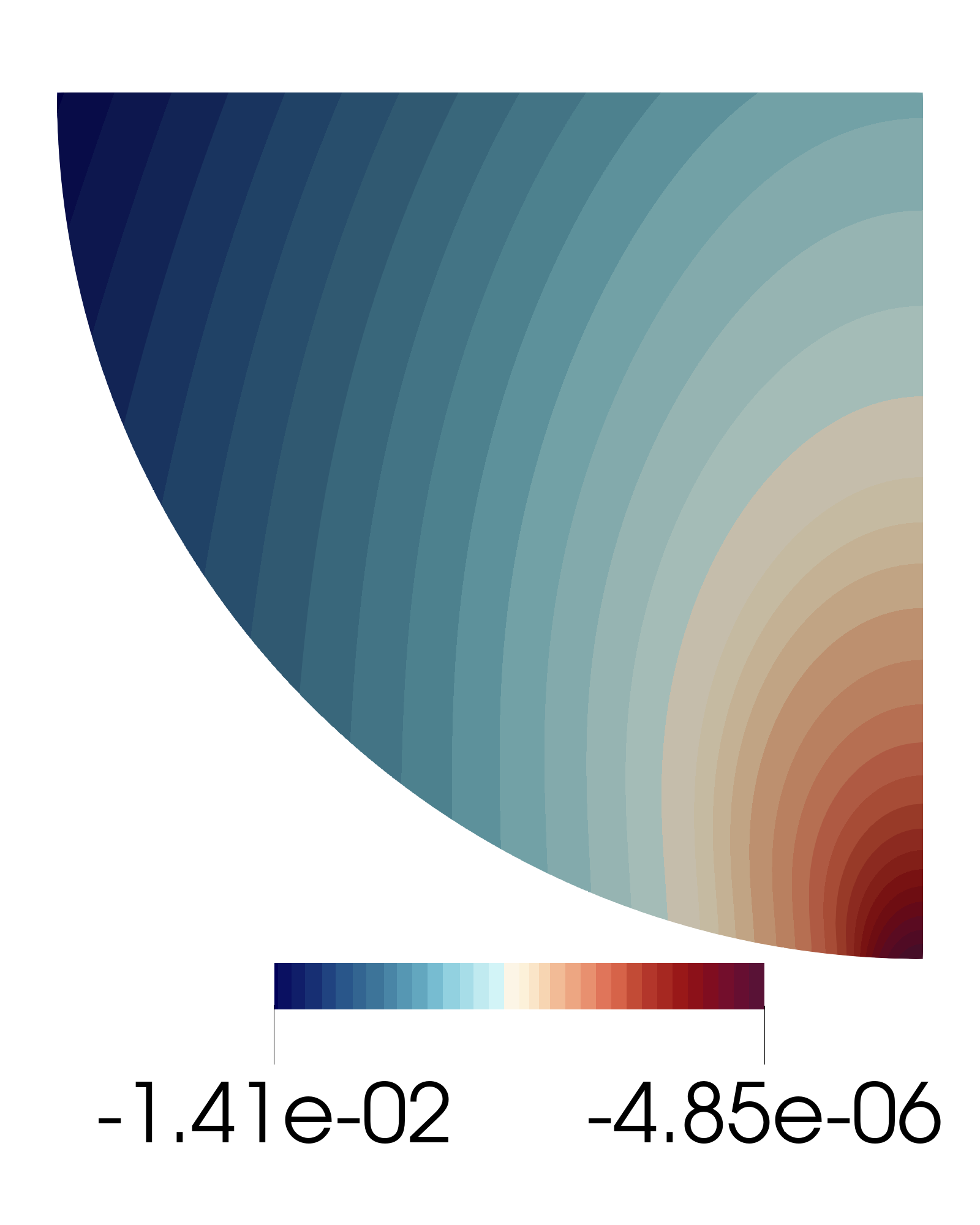}&
    \includegraphics[width=\figsize\linewidth]{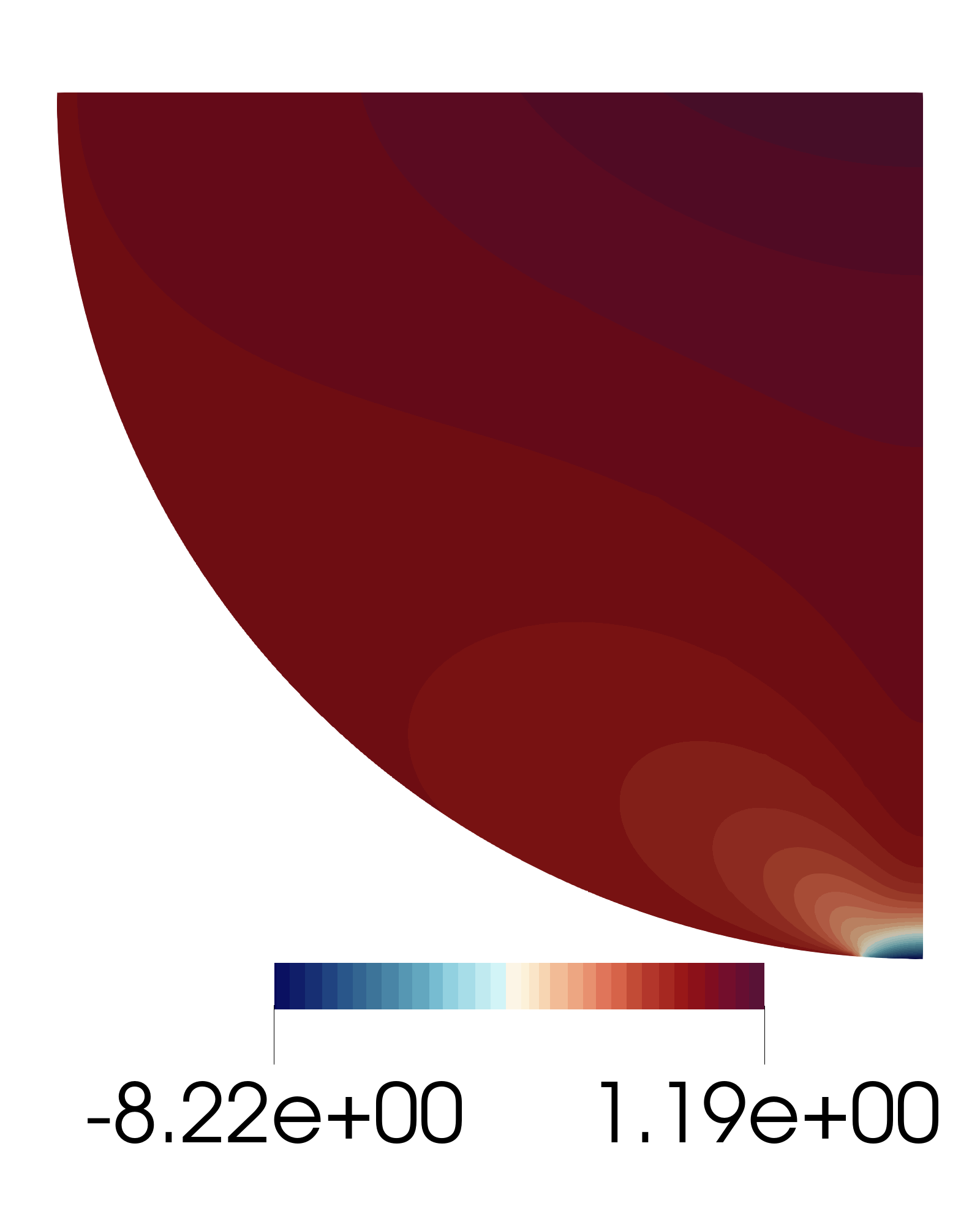}&
    \includegraphics[width=\figsize\linewidth]{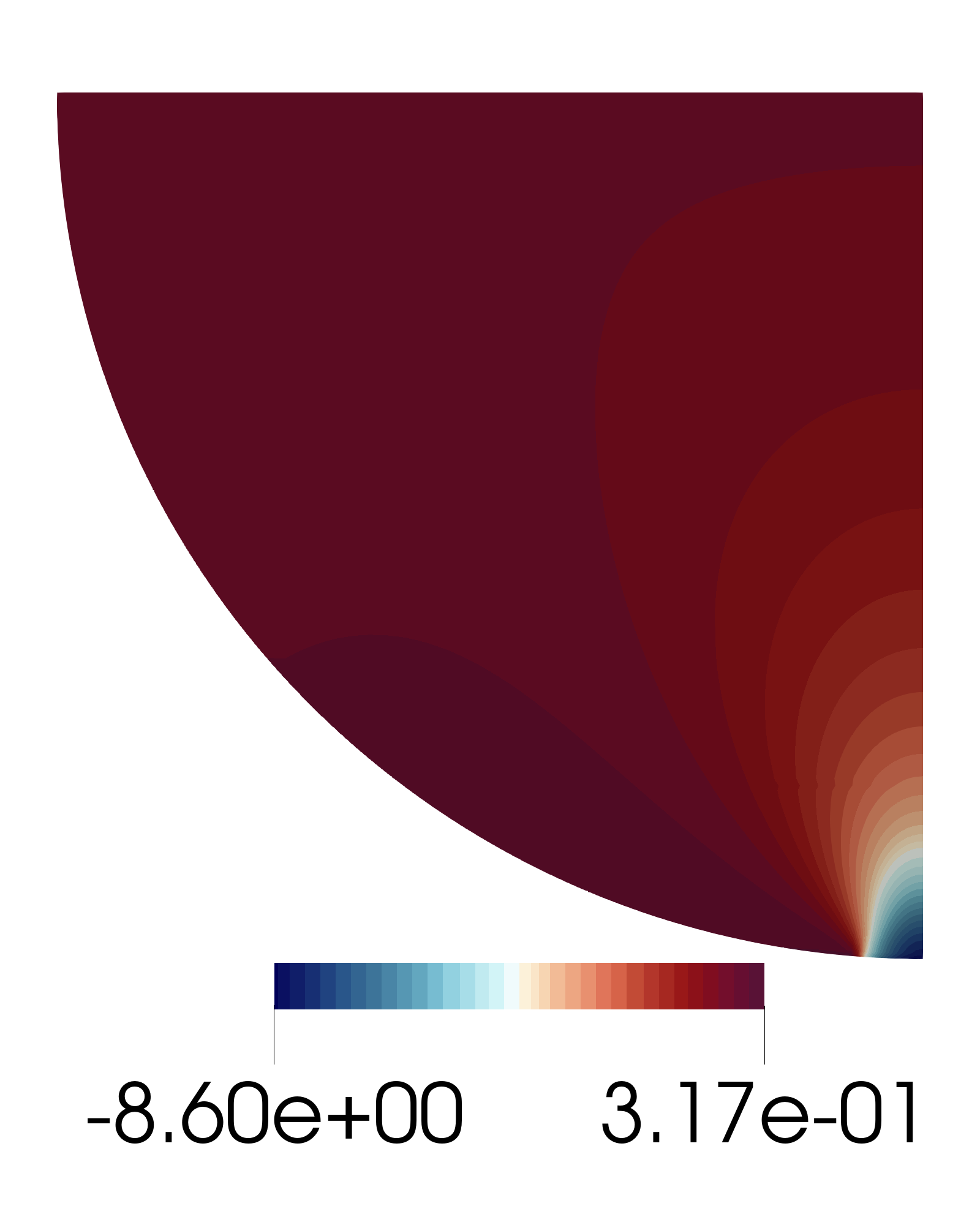}&
    \includegraphics[width=\figsize\linewidth]{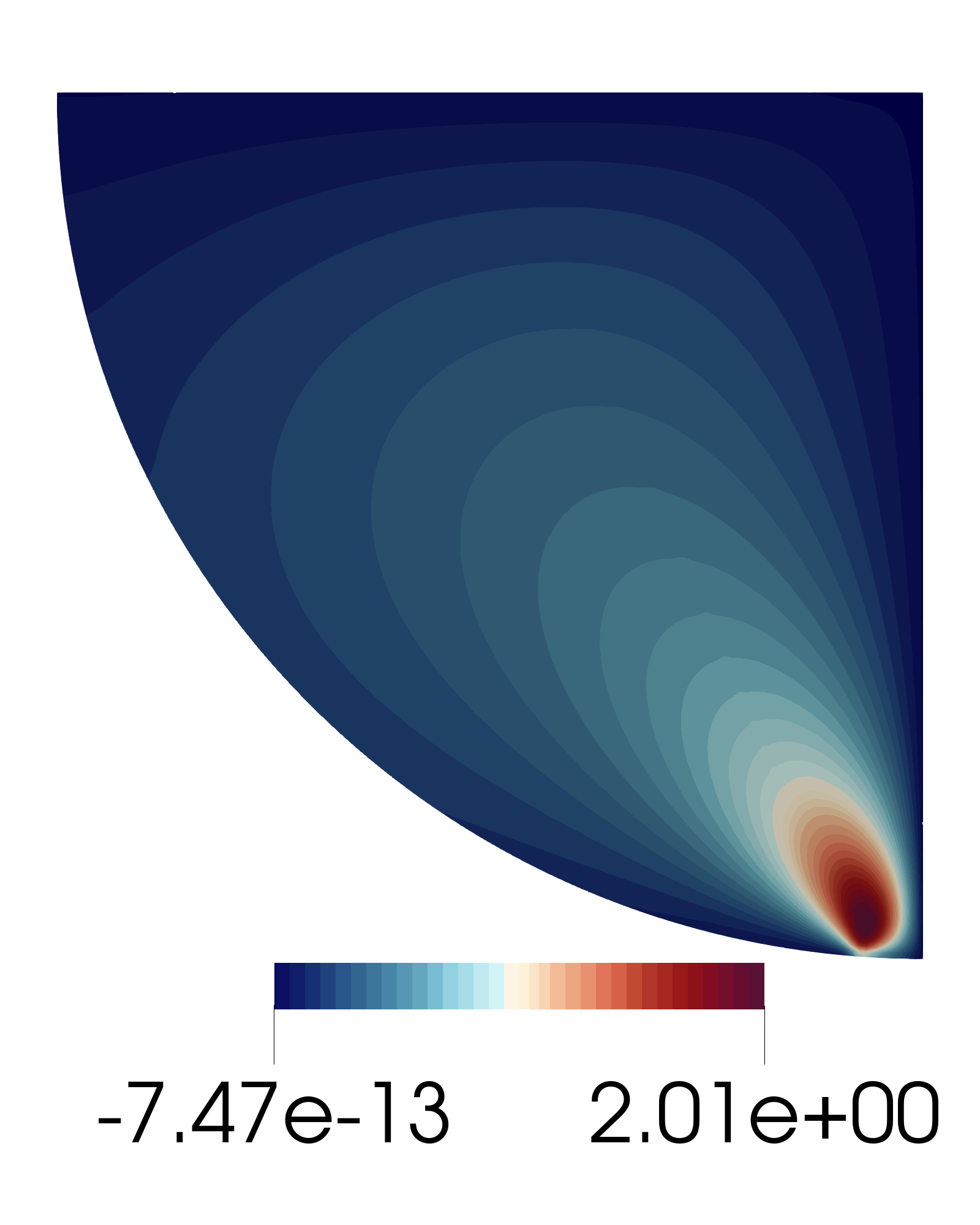}
    \\[-1ex]
    \end{tabular}
    \caption{PINNs as pure forward solver for the Hertzian contact problem. Comparison of stress and displacement components obtained by PINN and FEM.} %
    \label{fig:hertz_pinn_fem}
\end{figure}

\begin{figure}[thbp]
    \centering
    \includegraphics[width=0.375\linewidth]{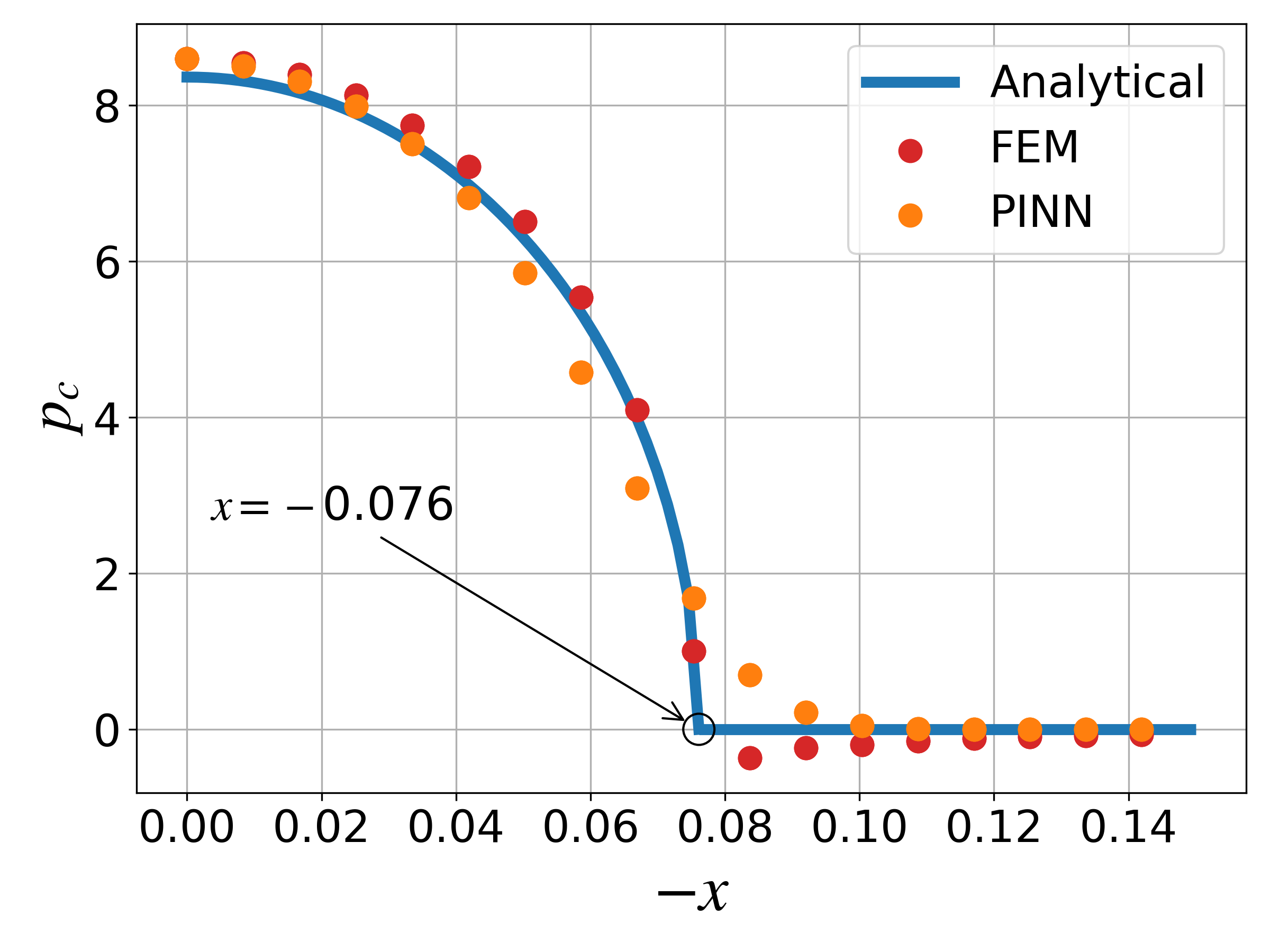}  
    \caption{Comparison of contact pressure distributions $p_c(x)$ obtained by analytical solution, PINN and FEM for case 1.
    }
    \label{fig:hertz_pinn_fem_analytical_pure}
\end{figure}

In the first use case, we deploy PINNs as a pure forward solver for contact problems to validate our approach. 
Training takes a total of 4049.8 \textit{s} and the prediction time is 0.091 \textit{s}. 
Displacement and stress components obtained through PINN and FEM are compared with contour plots in Figure \ref{fig:hertz_pinn_fem}.
Errors for displacement and stress fields are quantified using the relative ${L}_2$ norm.
As summarized in Table \ref{tab:case_1_2_results}, we obtain a relative error $E_{{L}_2}^{\boldsymbol{u}}=2.24\%$ for the displacement field 
and a relative error $E_{{L}_2}^{\boldsymbol{\sigma}}=3.74\%$ for the stress field.
Furthermore, the contact pressure distributions obtained via analytical solution, PINN and FEM are compared in Fig. \ref{fig:hertz_pinn_fem_analytical_pure}.
Since the zero traction boundary condition and the KKT constraints on the curved surface, numbered as edge 3 in Fig. \ref{fig:hertzian_geom}, are not 
enforced as hard but soft constraints, the PINN result can only resolve the kink at $x=\pm 0.076$ 
in an approximate manner depending on the chosen set of training points. 
Consequently, the contact pressure $p_c$ 
reduces to zero smoothly and slightly violates the zero traction boundary condition in the non-contact zone. 
Accordingly, the rather large error values $E_{{L}_2}^{\boldsymbol{\sigma}}$ are mostly related to this violation of the zero traction 
boundary conditions and KKT constraints close to the kink. 
Readers familiar with FEM modeling of contact problems will notice that a quite similar phenomenon occurs for 
mesh-based numerical methods where the transition between contact and non-contact zones cannot be perfectly resolved either 
and might even cause spurious oscillations in the case of higher-order interpolation functions \cite{Popp2018a} \cite{SEITZ2016259} \cite{popp2012siam}.

\begin{table}[thbp]
    \centering
    \begin{tabular}{lccccc}
        & $E_{{L}_2}^{\boldsymbol{u}}$ (\%)  
        & $E_{{L}_2}^{\boldsymbol{\sigma}}$ (\%)  
        & training time (s) 
        & prediction time (s)
        & $\int_{\Omega_c} \tilde{p}_c$
        \\ \midrule
        \textbf{case 1} 
        & 2.24 
        & 3.74 
        & 4049.8 
        & 0.091
        & 0.4993
        \\
        \textbf{case 2} 
        & 0.11 
        & 2.79 
        & 7126.7 
        & 0.092
        & 0.4978
    \end{tabular}
    \caption{Comparison of relative $L_2$ errors, training and prediction time for case 1 and case 2. 
    The term $\int_{\Omega_c} \tilde{p}_c$ denotes an integral of the predicted contact pressure $\tilde{p}_c$ over the potential contact boundary.}
    \label{tab:case_1_2_results}
\end{table}

\subsubsection{Case 2: PINNs as data-enhanced forward model}\label{sec_case_2}
\ifhidden
{}
\fi

\begin{figure}[thbp]
    \centering
    \includegraphics[trim={11cm 0cm 0cm 0cm},clip, width=0.575\linewidth]{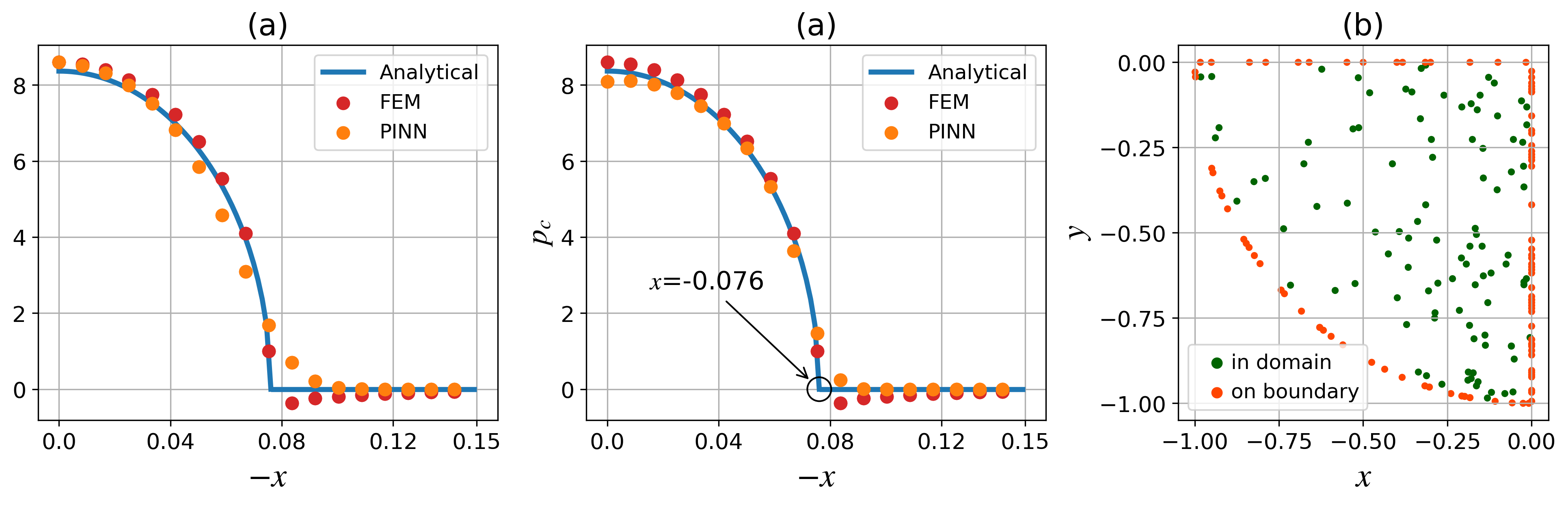}  
    \caption{Comparison of contact pressure distributions $p_c(x)$ obtained by analytical solution, PINN and FEM for case 2.
    (a) PINNs as a data-enhanced forward solver, 
    (b) distribution of data points generated through FEM simulations in the domain and on the boundary.
    }
    \label{fig:hertz_pinn_fem_analytical_data}
\end{figure}

One of the key features of PINNs is the capability of easily incorporating external data, such as measurement or simulation data, 
into the overall loss function.  
In this section, we enhance our PINN model with "artificial" measurement data obtained through FEM simulations, 
namely, data points for displacement and stress fields, to achieve better accuracy. 
The incorporated FEM data points are randomly selected, with 100 being selected within the domain and 100 being chosen along the boundary as 
depicted in Fig. \ref{fig:hertz_pinn_fem_analytical_data}b. 

A comparison of the contact pressure $p_c$ in the case of data enhancement 
with analytical solution and FEM reference solution is given in Fig. \ref{fig:hertz_pinn_fem_analytical_data}a. 
While the PINN accuracy is significantly improved upon close to the kink at $x=\pm 0.076$, the data-enhanced model underestimates the normal contact pressure around 
the origin. 
The relative $L_2$ errors confirm this assessment. While $E_{{L}_2}^{\boldsymbol{\sigma}}$ is only slightly reduced, 
$E_{{L}_2}^{\boldsymbol{u}}$ benefits dramatically from data enhancement. 
As provided in Table \ref{tab:case_1_2_results}, the integrated contact pressure is close to the applied load, $p=0.5$,
due to the conservation of momentum. To fulfill the momentum equation, overshooting after the kink is balanced by undershooting around the origin. 
Similarly, in case 1, overshooting after the kink is balanced by the undershooting from around $x=0.04$ to the kink. Moreover, 
we observe that results significantly rely on selecting appropriate additional data loss weights $w^{\mathrm{EXPs}}$ 
(reported in Table \ref{tab:hertz_table}). 
So far we identified the parameters through manual adjustments. In general, we recommend 
a hyperparameter analysis to determine optimal loss weights. 
Additionally, the data-enhanced model requires more training time compared to the pure forward model, 
which can be explained by the need to evaluate additional loss terms. However, after training is finished, the more accurate prediction takes
essentially the same time as in case 1. 

\subsubsection{Case 3: PINNs as inverse solver for parameter identification}\label{sec_case_3}
\ifhidden
{}
\fi
An interesting approach to solve an inverse problem is to simply add the unknown parameter to the set of network trainable parameters $\boldsymbol{\theta}$,
and it can then be identified with the help of additional loss terms based on the difference 
between predictions and observations (see Eq. \ref{eq:inverse_theta}). 
In the following, we exemplarily identify the applied external pressure $p$ acting on the half-cylinder
using FEM results as "artificial" measurement data. 
\begin{equation}
    \label{eq:inverse_theta}
    \boldsymbol{\Theta}^*=\underset{\boldsymbol{\Theta}}{\arg \min } \ \mathcal{L}_{\mathrm{C}}(\boldsymbol{\Theta}) 
    \quad \text{where,} \quad \boldsymbol{\Theta}=(p,\boldsymbol{\theta}).
\end{equation}

\begin{figure}[!thbp]
    \centering
    \includegraphics[width=0.90\linewidth]{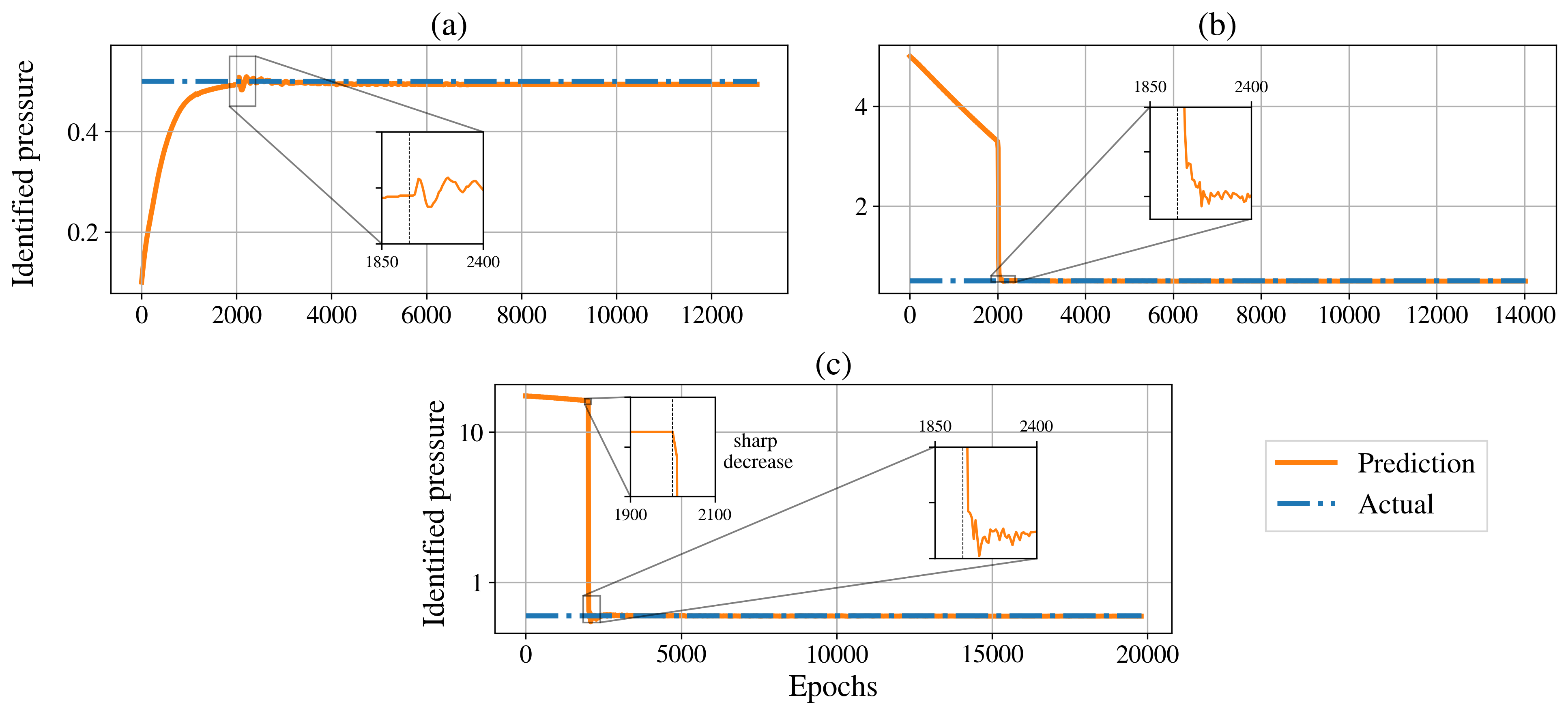}  
    \caption{Identification of the applied external pressure on the half-cylinder in case of different initial guesses
    (a) $p_o=0.1$, (b) $p_o=5$, (c) $p_o=20$.}
    \label{fig:hertz_inverse}
\end{figure}

Fig. \ref{fig:hertz_inverse} shows the convergence behavior of the identified pressure $\tilde{p}$ compared to the actual pressure $p$
through the number of epochs for different initial guesses $\tilde{p}_o$. 
First, we start training with a "good" initial guess  $\tilde{p}_o=0.1$ being quite close to 
the actual pressure $p=0.5$, and then we increase it to $\tilde{p}_o=20$ to measure how sensitive the PINN is to the choice of the initial guess. 
As depicted in Fig. \ref{fig:hertz_inverse}c, convergence can be achieved even when a relatively large and therefore unphysical initial guess is made. 
There is a steep increase in convergence rate after 2000 epochs, 
which can be explained by the fact that we switch to the \textit{L-BFGS-B} optimizer there. As mentioned in Section \ref{sec_lame}, 
we deploy \textit{Adam} for 2000 epochs to avoid the optimization process getting stuck in local minima.   
Note, however, that switching from \textit{Adam} to \textit{L-BFGS-B} introduces oscillations in the transition region so that
transition has to be handled with care. 
The relative error, denoted as $E$, is used to measure the difference between the identified and actual 
pressure values, resulting in the same value of 1.2\% for all three initial guesses $p_o=0.1$, $p_o=5$, and $p_o=20$ 
as provided in Table \ref{tab:hertz_inverse_results}. 
Additionally, a larger initial guess requires more computing time and epochs as expected.   

\begin{table}[thbp]
    \centering
    \begin{tabular}{lccccc}
        &  \multicolumn{1}{c}{\begin{tabular}[c]{@{}c@{}}initial guess \\ $\tilde{p}_o$\end{tabular}}    
        &  \multicolumn{1}{c}{\begin{tabular}[c]{@{}c@{}}identified \\ $\tilde{p}_f$\end{tabular}}                                          
        & \multicolumn{1}{c}{\begin{tabular}[c]{@{}c@{}}relative error \\ $E(\%)$ \end{tabular}}  
        & training time (s) 
        & no. of epochs \\
        \midrule
        \multirow{3}{*}{\textbf{case 3}} & 0.1 & 0.494 & 1.2 & 3472.8 & 12958 \\
                                         & 5   & 0.494  & 1.2 & 3762.2 & 14019 \\
                                         & 20  & 0.494 & 1.2 & 5394.1 & 19793
    \end{tabular}
    \caption{Comparison of relative errors, training time and number of epochs for the inverse Hertzian problem 
    used to identify the applied external pressure. The relative error is defined as $E^{*}=\text{abs}\left((\tilde{*} - *)/*\right)$.}
    \label{tab:hertz_inverse_results}
\end{table}

\subsubsection{Case 4: PINNs as fast-to-evaluate surrogate model}
\ifhidden
{}
\fi

\begin{figure}[!thbp]
    \centering
    \includegraphics[width=0.675\linewidth]{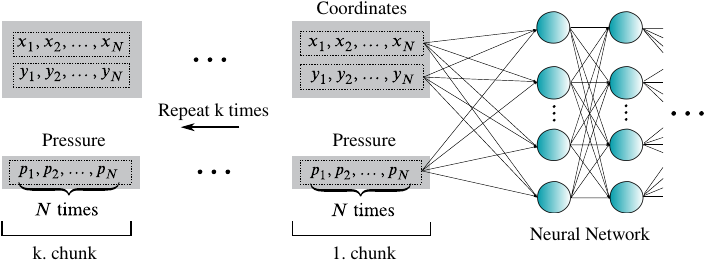}  
    \caption{An illustration of the procedure to include the applied external pressure as a neural network input.}    
    \label{fig:hertz_surrogate}
\end{figure}

In the last use case, the load (applied external pressure) is considered as an additional network input, i.e. 
$\netInput = (\boldsymbol{x}, p)$ (compare Eq. \ref{eq:approx}), to construct a fast-to-evaluate 
surrogate model that is capable of predicting displacement and stress fields for different pressure values. 
Since a single network is deployed, the length of each network input
must be the same. 
We sample the three-dimensional input space with $N = 1946$ training points ($N = N_{rp} + N_{bp}$ (Eq. \ref{eq:pinn_elasticity_detail})). 
While the spatial coordinates $\boldsymbol{x}$ are selected from a two-dimensional mesh (as in Section \ref{sec_case_1}, \ref{sec_case_2}, and \ref{sec_case_3}), 
the third (pressure) component of the input vector is drawn randomly from a uniform distribution over the 
considered pressure range.
To improve the accuracy of the prediction, the $N$ sampling points are repeated $k$ times. 
In the context of network training, we refer to one instance of the $N$ distinct sampling points as one \textit{chunk}, 
so that $k$ is the number of chunks as depicted in Fig. \ref{fig:hertz_surrogate}.
Indeed, this process increases the computing time and complexity of the model, since the input size increases from $(n,2)$ 
to $(n\cdot k, 3)$. However, such a method can lead to better accuracy since the network is trained with a larger data set 
and also it enables batch training to reduce computational effort \cite{geron2022hands}.  
For our specific example, we sample pressure values from a range of $[0.2, 1.0]$ and
we consider only two different numbers of chunks, namely $k=1$ and $k=5$.  

As shown in Figure \ref{fig:hertz_surrogate_results}, employing a single chunk is insufficient to accurately 
capture the influence of the applied external pressure $p$ on the contact pressure distribution $p_c(\pdeInput,p)$. 
However, increasing the number of chunks from $k=1$ to $k=5$ 
increases the accuracy.
Table \ref{tab:hertz_surrogate_results} provides the relative $L_2$ errors for the contact pressure distribution $p_c(\pdeInput,p)$
between the analytical solution and the predictions of surrogate PINN models, 
and it can be seen that increasing the number of chunks indeed results in improved accuracy of the surrogate model. 
Nonetheless, the overall error level of 10\%-16\% even in the case $k=5$ is still too high from an engineering perspective and will
be subject to further investigations. This example is only intended as a very first proof of concept, and we have already identified
several algorithmic modifications that could possibly increase the accuracy of the surrogate model.  
For example, we use fixed loss weights even though pressure values in the input layer vary. 
Thus, adaptive loss weights should be implemented since the PINN
accuracy highly depends on choosing appropriate loss weights. 
Additionally, the accuracy of PINNs could be further improved 
by increasing the number of chunks.

\begin{table}[!htbp]
    \centering
    \def\arraystretch{1.3}
    \def\hcoldist{0.5cm}
    \begin{tabular}{@{}cl@{\hskip \hcoldist}c@{\hskip \hcoldist}c@{\hskip \hcoldist}c}
    & 
    & 
    \multicolumn{3}{c}{\textbf{Applied pressure}} 
    \\ \cmidrule{3-5}                                                 
    &   
    & $p=0.45$ 
    & $p=0.98$ 
    & $p=1.5$ 
    \\ \cmidrule{2-5}
    \multirow{2}{*}{\textbf{chunk size}}
    & $k=1$
    & 22.81\% & 17.23\% & 16.93\%\\
    & $k=5$   
    & 16.30\% & 11.47\% & 10.87\%
    \end{tabular}
    \caption{Comparison of relative ${L}_2$ errors for the contact pressure distribution $p_c(\pdeInput,p)$ between the analytical solution and the predictions of our surrogate PINN models.
    Predictions are based on unseen pressure values.} 
    \label{tab:hertz_surrogate_results}
\end{table}

\begin{figure}[thbp]
    \def\figsize{0.7}
    \newcommand{\trimsize}{1.75}
    \newcommand{\trimsizetop}{1.}
    \def\hskiplocal{\hskip 0.1cm}
    \centering
    \def\arraystretch{0.5}%
    \begin{tabular}{c@{ }c}
    \includegraphics[trim={0cm \trimsize cm 0cm 0cm},clip, width=\figsize\linewidth]{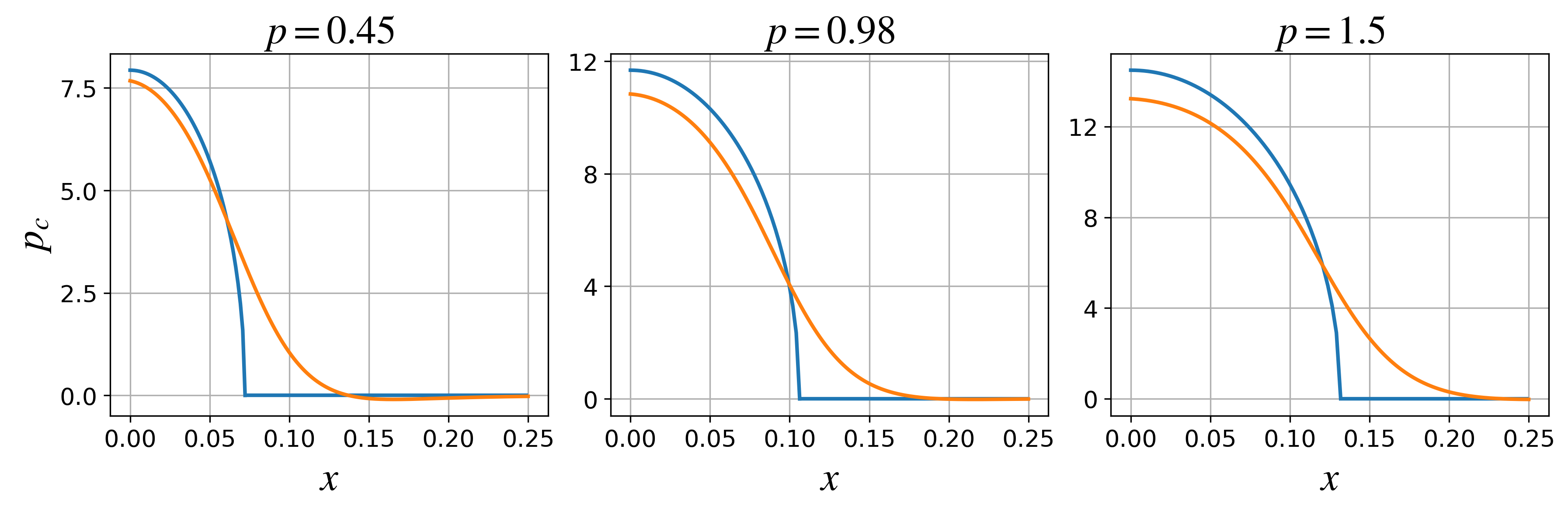}&
    \rownameform{\hspace*{-0.8cm} $k=1$}
    \\
    \includegraphics[trim={0cm 0cm 0cm \trimsizetop cm},clip, width=\figsize\linewidth]{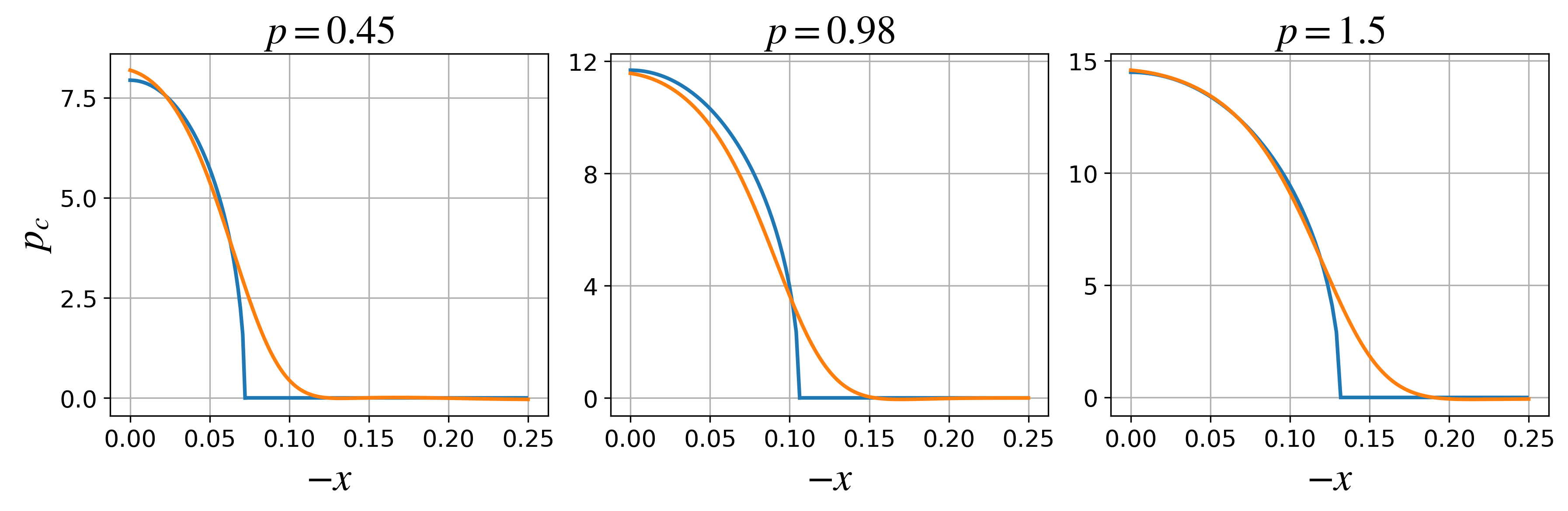}&
    \rownameform{\hspace*{0.8cm} $k=5$}
    \\
    \end{tabular}
    \caption{Comparison of the contact pressure distribution obtained through a PINN-based surrogate model
     to the analytical solution. The PINNs are trained using different numbers of chunks.} 
    \label{fig:hertz_surrogate_results}
\end{figure}

\section{Conclusion}\label{sec5}
In this study, we have presented an extension of physics-informed neural networks (PINNs) for solving forward and 
inverse problems of contact mechanics
under the assumption of linear elasticity. The framework has been tested on several benchmark examples
with different use cases, e.g. the Hertzian contact problem, and has been 
validated by existing analytical solutions or numerical simulations using the finite element method (FEM). 
As an alternative way of soft constraint enforcement as compared to existing methods, 
a nonlinear complementarity problem (NCP) function, namely \textit{Fischer-Burmeister}, 
is explored and exploited to enforce the inequality constraints inherent to contact problems. 
This aspect has not been investigated in the context of PINNs so far to the best of the authors' knowledge. Besides using 
PINNs as pure forward PDE solver,
we show that PINNs can serve as a hybrid model enhanced by experimental and/or simulation data to identify unknown parameters of contact
problems, e.g. the applied external pressure. We even go one step further and deploy PINNs as fast-to-evaluate surrogate models,
and could at least obtain a first proof of concept up to a certain level of accuracy.
   
A question that has emerged recently  
is whether data-driven approaches such as PINNs will replace classical numerical methods such as FEM in the near future. 
Within this study, we only considered benchmark examples that have been developed 
and solved decades ago using the FEM. Even for these simple examples, we came to the conclusion that deploying PINNs 
as forward solvers for contact mechanics can not compete with FEM
in terms of computational performance and accuracy.
Therefore, we doubt the applicability of PINNs to complex engineering problems without data enhancement. However, 
PINNs can be a good candidate for solving data-enhanced forward problems and especially inverse problems due to the easy integration of additional data.
Similarly, PINNs can break the curse of dimensionality of parametric models, so that more complex surrogate models can be generated.
Also, it is observed that minimizing multiple loss functions simultaneously is one of the most 
significant challenges in training PINNs, 
and current optimization algorithms are not tailored to addressing this challenge.
Therefore, using multi-objective optimization algorithms that are particularly designed for PINNs has the potential to be a 
gamechanger in improving their overall performance and accuracy.    
We believe that hybrid strategies can be a promising option to construct mixed models 
to benefit from the advantages of both classical and data-driven approaches.  

This study reveals several possibilities for further exploration and investigation. 
Although the proposed PINN formulation for benchmark examples demonstrates acceptable results, 
further applications, particularly on complex domains including three-dimensional problems should be analyzed. 
As an alternative strategy to scaling network outputs, 
a non-dimensionalized contact formulation can be implemented.
Different NCP functions other than the \textit{Fischer-Burmeister} function can be further investigated. 
The inverse solver has been applied to identify the applied external pressure, but it can be extended to also 
predict internal material parameters.  
Additionally, a hyperparameter optimization study can be performed to tune loss weights, network architecture 
and optimizer parameters. Last but not least, related techniques such as variational PINNs might overcome the limitations of collocation
inherent to PINNs and instead provide a sound variational framework.

\section{Acknowledgement}

This research paper is funded by dtec.bw - Digitalization and Technology Research Center of the 
Bundeswehr. dtec.bw is funded by the European Union - NextGenerationEU.
The authors gratefully acknowledge the computing resources provided by the Data Science \& Computing Lab at the University of the Bundeswehr Munich.

\clearpage

\nocite{*}
\bibliography{wileyNJD-AMA}%

\appendix
\section{Loss formulation for 2D linear isotropic elasticity under plane strain conditions}\label{append_2d_loss}

The elasticity tensor $\mathbb{C}$ under plane strain conditions can be expressed in terms of Lam\'e constants
$\lambda$ and $\mu$ as 
\begin{equation}
    \arraycolsep=8pt
    \label{eq:2d_elasticity_tensor}
    \mathbb{C}=\left[\begin{array}{ccc}
    2 \mu+\lambda & \lambda & 0 \\
    \lambda & 2 \mu+\lambda & 0 \\
    0 & 0 & \mu
    \end{array}\right].
\end{equation}

\afterequation 
\noindent Inserting Eq. \ref{eq:2d_elasticity_tensor} into Eq. \ref{eq:pinn_elasticity_detail} 
we can obtain the total loss for 2D linear isotropic elasticity under the plane strain condition as
\begin{equation}
    \label{eq:loss_elasticity_2d}
    \begin{aligned}
        \pazocal{L}_{\mathrm{E}} = &\ \pazocal{L}_{\mathrm{PDEs}} + \pazocal{L}_{\mathrm{DBCs}} + \pazocal{L}_{\mathrm{NBCs}} + \pazocal{L}_{\mathrm{EXPs}}\\
        = &\ w^{(\mathrm{PDEs})}_1 \bigl|\pinnSComponent{xx, x}+ \pinnSComponent{xy, y} + \bodyforcetwodim_x\bigr|_{\Omega} +
        w^{(\mathrm{PDEs})}_2\bigl|\pinnSComponent{yx, x}+\pinnSComponent{yy, y} + \bodyforcetwodim_y \bigr|_{\Omega} + \\
        &\ w^{(\mathrm{PDEs})}_3\bigl|\pinnSComponent{xx}-(\lambda + 2\mu)\pinnEComponent{xx}-\lambda\pinnEComponent{yy}\bigr|_{\Omega} + 
        w^{(\mathrm{PDEs})}_4\bigl|\pinnSComponent{yy}-\lambda\pinnEComponent{xx}-(\lambda + 2\mu)\pinnEComponent{yy}\bigr|_{\Omega} + 
        w^{(\mathrm{PDEs})}_5\bigl|\pinnSComponent{xy}-2\mu \pinnEComponent{xy}\bigr|_{\Omega} + \\
        &\ w^{(\mathrm{DBCs})}_1\bigl|\pinnUComponent{x} - \hat u_x\bigr|_{\partial \Omega_D} + w^{(\mathrm{DBCs})}_2\bigl|\pinnUComponent{y} - \hat u_y\bigr|_{\partial \Omega_D} + \\
        &\ w^{(\mathrm{NBCs})}_1\bigl|\pinnSComponent{xx} n_x + \pinnSComponent{xy} n_y - \hat t_x\bigr|_{\partial \Omega_N} + w^{(\mathrm{NBCs})}_2\bigl|\pinnSComponent{yx} n_x + \pinnSComponent{yy} n_y - \hat t_y\bigr|_{\partial \Omega_N} + \\
        &\ w^{(\mathrm{EXPs})}_1\bigl|\pinnUComponent{x}- u^*_x\bigr|_{\Omega_e} + w^{(\mathrm{EXPs})}_2\bigl|\pinnUComponent{y} - u^*_y\bigr|_{\Omega_e} + w^{(\mathrm{exps})}_3\bigl|\pinnSComponent{xx} - \sigma^*_{xx}\bigr|_{\Omega_e}
        + w^{(\mathrm{EXPs})}_4\bigl|\pinnSComponent{yy} - \sigma^*_{yy}\bigr|_{\Omega_e} + w^{(\mathrm{EXPs})}_5\bigl|\pinnSComponent{xy} - \sigma^*_{xy}\bigr|_{\Omega_e}
    \end{aligned}
\end{equation}

\afterequation
including kinematics
\begin{equation}
    \pinnEComponent{xx} =\pinnUComponent{x,x}, \quad \pinnEComponent{yy} =\pinnUComponent{y,y}, \quad  \pinnEComponent{xy} =\frac{1}{2}(\pinnUComponent{x,y}+\pinnUComponent{y,x}),
\end{equation}

\afterequation
where 
\begin{equation}
    \begin{aligned}
        \pinnUComponent{x} \approx (\mathcal{N}_{u_x}({x,y}))^{'}, & \quad  \pinnUComponent{y} \approx (\mathcal{N}_{u_y}({x,y}))^{'},\\
        \pinnSComponent{xx} \approx (\mathcal{N}_{\sigma_{xx}}({x,y}))^{'}, \quad
        \pinnSComponent{yy} \approx (\mathcal{N}_{\sigma_{yy}}&({x,y}))^{'}, \quad
        \{\pinnSComponent{xy}=\pinnSComponent{yx}\} \approx (\mathcal{N}_{\sigma_{xy}}({x,y}))^{'}.
    \end{aligned}
\end{equation}

\afterequation
The out-of-plane stress component $\sigma_{zz}$ is not considered as network output since it can be calculated in the post-processing. 
\section{Additional error comparisons}\label{append_error_comparison}
The vector-based $L_2$ error between the approximated PINN solution $\tilde{f}$ and the analytical solution $f$, 
is denoted as, 
\begin{equation}
    \vectorbasedError{f} := \sqrt{\sum_{j=1}^{N_{\text{test}}}
    \left(\tilde{f}(\boldsymbol{x}^j)-f(\boldsymbol{x}^j)\right)^2}.
\end{equation}

\noindent The corresponding integral-based $L_2$ error between the approximated PINN solution $\tilde{f}$ 
and the analytical solution $f$, is denoted as, 
\begin{equation}
    \integralbasedError{f} := \sqrt{\frac{N_{\text{test}}}{\int_{\Omega} d\boldsymbol{x}} \int_{\Omega} 
    \left(\tilde{f}(\boldsymbol{x})-f(\boldsymbol{x})\right)^2 d\boldsymbol{x}}, 
\end{equation}

where $\int_{\Omega} d\boldsymbol{x}$ represents the area for Table \ref{tab:error_comparion_lame}
and Table \ref{tab:error_comparion_patch}, while it represents the arc length for Table \ref{tab:error_comparion_hertz},
since the contact pressure $p_c$ is integrated over the potential contact boundary $\partial \Omega_c$. 
We refer to sections \ref{sec_lame}, \ref{sec_patch}, \ref{sec_hertz} for the respective number of test points $N_{\text{test}}$. 
Note that in this first study, we used the vector-based error measure that is frequently used for PINNs 
\cite{lu2021physics} \cite{xu2022transfer} and easily implemented. 
For future studies, we suggest more expressive integral-based error estimates.

\begin{table}[!htbp]
    \centering
    \def\arraystretch{1.3}
    \def\hcoldist{1cm}
    \begin{tabular}{cccc@{\hskip \hcoldist}cccc}
        \multicolumn{4}{c}{\hspace{-1cm}\textbf{Vector-based}} 
        & \multicolumn{4}{c}{\textbf{Integral-based}} \\ \cmidrule(lr{1cm}){1-4} \cmidrule{5-8}
        $\relErrorUComponent{r}$ 
        & $\relErrorSComponent{rr}$    
        & $\relErrorSComponent{\theta \theta}$       
        & $\relErrorSComponent{r \theta}$    
        & $\integralbasedError{u_r}$ 
        & $\integralbasedError{\sigma_{rr}}$    
        & $\integralbasedError{\sigma_{\theta \theta}}$       
        & $\integralbasedError{\sigma_{r \theta}}$
        \\ \cmidrule(lr{1cm}){1-4} \cmidrule{5-8}
        1.06e-5
        & 0.012    
        & 0.035    
        & 0.008  
        
        & 8.80-e6  
        & 0.010    
        & 0.030    
        & 0.006  
    \end{tabular}
    \caption{Comparison of the vector-based and integral-based $L_2$ error for
    for the Lam\'e problem of elasticity. The evaluated quantities are 
    displacement and stress components in polar coordinates
    (see Section \ref{sec_lame} for the example setup).}
    \label{tab:error_comparion_lame}
\end{table}
\aftertable

\begin{table}[!htbp]
    \centering
    \def\arraystretch{1.3}
    \def\hcoldist{0.2cm}
    \begin{tabular}{c ccccc@{\hskip \hcoldist}ccccc}
        &
        \multicolumn{5}{c}{\textbf{Vector-based}} 
        & \multicolumn{5}{c}{\textbf{Integral-based}} \\ \cmidrule(lr{0.5cm}){2-6} \cmidrule{7-11}
        
        & $\relErrorUComponent{x}$    
        & $\relErrorUComponent{y}$       
        & $\relErrorSComponent{xx}$   
        & $\relErrorSComponent{yy}$  
        & $\relErrorSComponent{xy}$  
        & $\integralbasedError{u_{x}}$    
        & $\integralbasedError{u_{y}}$       
        & $\integralbasedError{\sigma_{xx}}$ 
        & $\integralbasedError{\sigma_{yy}}$ 
        & $\integralbasedError{\sigma_{xy}}$ 
        \\ \midrule
        \textbf{sign}
        & 3.98e-3 & 1.753e-2 & 1.05e-2 & 1.166e-2 & 5.68e-3
        & 3.92e-3 & 1.750e-2 & 1.04e-2 & 1.156e-2 & 5.74e-3
        \\
        \textbf{Sigmoid}
        & 1.91e-3 & 3.79e-3 & 5.89e-3 & 7.419e-3 & 3.86e-3
        & 1.85e-3 & 3.75e-3 & 5.81e-3 & 7.421e-3 & 3.91e-3
        \\ 
        \multicolumn{1}{c}{\def\arraystretch{1}\begin{tabular}[c]{@{}c@{}}\textbf{\textit{Fischer-}} \\ \textbf{\textit{Burmeister}}\end{tabular}} 
        & 8.46e-4 & 8.04e-4 & 2.28e-3 & 2.10e-3 & 1.20e-3
        & 8.20e-4 & 7.80e-4 & 2.26e-3 & 2.08e-3 & 1.22e-3
    \end{tabular}
    \caption{Comparison of the vector-based and integral-based $L_2$ error for 
    the contact problem between an elastic block and the rigid domain. Evaluated quantities are 
    displacement and stress component in cartesian coordinates.
    Three different methods are provided to enforce KKT constraints (see Section \ref{sec_patch}
    for the example setup).}
    \label{tab:error_comparion_patch}
\end{table}
\aftertable

\begin{table}[!htbp]
    \centering
    \def\arraystretch{1.3}
    \def\hcoldist{0.2cm}
    \begin{tabular}{c c@{\hskip \hcoldist}c}
        &
        \multicolumn{1}{c}{\textbf{Vector-based}} 
        & \multicolumn{1}{c}{\textbf{Integral-based}} \\ \cmidrule(lr){2-2} \cmidrule{3-3}
        
        & $E_{{L}_2}^{{p}_{c}}$    
        & $\integralbasedError{p_{c}}$ 
        \\ \midrule
        \textbf{case 1}
        & 2.22 & 2.48
        \\
        \textbf{case 2}
        & 0.89 & 1.27
    \end{tabular}
    \caption{Comparison of the vector-based and integral-based $L_2$ error for the contact pressure $p_c$
    of the Hertzian contact problem for different cases (see Section \ref{sec_case_1} and Section \ref{sec_case_2}
    for the example setup).}
    \label{tab:error_comparion_hertz}
\end{table}

\end{document}